\documentclass[12pt,draftcls,onecolumn]{IEEEtran}
\IEEEoverridecommandlockouts                              %
\usepackage{epsfig} %

\usepackage{amsmath,graphicx,amsfonts,amssymb,epsfig,subfig,mathrsfs,mathtools}
\allowdisplaybreaks
\usepackage{ntheorem}
\usepackage{cuted}
\usepackage{arydshln}
\usepackage{blkarray}

\usepackage{multirow}
\usepackage{color}
\usepackage{multicol}
\usepackage{lipsum}
\usepackage{rotating}
\usepackage{graphicx}%
\usepackage{algorithm}
\usepackage{booktabs}

\usepackage[utf8]{inputenc}
\usepackage{xspace}
\usepackage{algorithmicx}
\usepackage{algpseudocode}
\usepackage{cite}
\usepackage{cancel}
\usepackage{textcomp}
 \definecolor{lin}{RGB}{240,0,0}
 \definecolor{paleblue}{RGB}{0,9,255}
\usepackage[colorlinks=true,linkcolor=lin,citecolor=paleblue,pdfa]{hyperref}

\usepackage{epstopdf}
\newcommand{\map}[3]{#1: #2 \rightarrow #3}
\newcommand{\setdef}[2]{\{#1 \; | \; #2\}}
\newcommand{\st}{\ensuremath{\operatorname{s.t.}}}
\newcommand{\real}{\ensuremath{\mathbb{R}}}
\newcommand{\prob}{\ensuremath{\mathbb{P}}}
\newcommand{\realnonnegative}{\ensuremath{\mathbb{R}}_{\ge 0}}

\newcommand{\until}[1]{\{0,1,\dots, #1\}}
\newcommand{\untilone}[1]{\{1,\dots, #1\}}
\newcommand{\subscr}[2]{#1_{\textup{#2}}}
\newcommand{\supscr}[2]{#1^{\textup{#2}}}
\newcommand{\vect}[1]{\mathbf{#1}}
\newcommand{\vectorones}[1]{\vect{1}_{#1}}
\newcommand{\vectorzeros}[1]{\vect{0}_{#1}}
\newcommand{\Norm}[1]{\|#1\|}
\newcommand{\trans}[1]{{#1}^\top}
\newcommand{\obj}{\ensuremath{\operatorname{obj}}}
\newcommand{\Prob}{\ensuremath{\operatorname{Prob}}}
\newcommand{\sol}{\ensuremath{\operatorname{sol}}}

\newcommand{\untilinterval}[2]{\{{#1},\dots, {#2}\}}
\newcommand{\UB}{\ensuremath{\operatorname{UB}}}
\newcommand{\LB}{\ensuremath{\operatorname{LB}}}
\newcommand{\graph}{\mathcal{G}}
\newcommand{\nodes}{\mathcal{V}}
\newcommand{\edges}{\mathcal{E}}

\newcommand{\len}{\operatorname{len}}
\newcommand{\lane}{\operatorname{lane}}

\DeclareMathOperator*{\argmax}{argmax}

\makeatletter
\newtheoremstyle{breaknote}%
  {\item{\theorem@headerfont
          ##1\ ##2\theorem@separator}\hskip\labelsep\relax}%
  {\item{\theorem@headerfont
          ##1\ ##2\ (##3)\theorem@separator}\hskip\labelsep\relax}
\newcommand{\leqnomode}{\tagsleft@true}
\newcommand{\reqnomode}{\tagsleft@false}
\makeatother
\reqnomode
\theoremstyle{breaknote}
\newtheorem{assumption}{Assumption}[section]
\newtheorem{theorem}{Theorem}[section]

\newtheorem{lemma}{Lemma}[section]
\newtheorem{proposition}{Proposition}[section]
\newtheorem{remark}{Remark}[section]

\newcommand{\qed}{~\hfill{\small $\blacksquare$}}

\newcommand{\reviseone}[1]{{#1}}

\allowdisplaybreaks
\renewcommand{\baselinestretch}{.977}

\title{Data-driven Predictive Control for a Class of Uncertain Control-Affine Systems}

\author{Dan Li$^{1}$, Dariush Fooladivanda$^{1}$ and Sonia
  Mart{\'\i}nez$^{1}$ \thanks{*This research was developed with
    funding from ONR contract N00014-19-1-2471 and DARPA (Lagrange) contract N660011824027. The
    views, opinions and/or findings expressed are those of the author
    and should not be interpreted as representing the official views
    or policies of the Department of Defense or the U.S. Government.}
  \thanks{$^{1}$D. Li, D. Fooladivanda and S. Mart{\'\i}nez are with
    the Department of Mechanical and Aerospace Engineering, University
    of California San Diego, La Jolla, CA 92092, USA {\tt\small
      lidan@ucsd.edu; dfooladivanda@ucsd.edu;soniamd@ucsd.edu}}%
}
\begin{document}

\maketitle

\begin{abstract}
  This paper studies a data-driven predictive control for a class of
  control-affine systems which is subject to uncertainty. With the
  accessibility to finite sample measurements of the uncertain
  variables, we aim to find controls which are feasible and provide
  superior performance guarantees with high probability.  This results
  into the formulation of a stochastic optimization problem (P), which
  is intractable due to the unknown distribution of the uncertainty
  variables. By developing a distributionally robust optimization
  framework, we present an equivalent and yet tractable reformulation
  of (P). Further, we propose an efficient algorithm that provides
  online suboptimal data-driven solutions and guarantees performance
  with high probability. To illustrate the effectiveness of the
  proposed approach, we consider a highway speed-limit control
  problem. We then develop a set of data-driven speed controls that
  allow us to prevent traffic congestion with high
  probability. Finally, we employ the resulting control method on a
  traffic simulator to illustrate the effectiveness of this approach
  numerically.
\end{abstract}

\section{Introduction}
Uncertainty is ubiquitous in real-world complex systems, including
financial markets~\cite{DA-MG-GG-SK:14}, human-robot mixed
systems~\cite{ST-WB-DF:05,ND-BM:12,BDL-MK-JPH:10} or transportation
networks~\cite{AN-ELG-VD:02,AG-YG-TL-KJH-FB:13,DL-DF-SM:19-ecc}.
Thanks to new advances in computation, communication, and innovative
infrastructure, large amounts of data have become widely accessible,
which can help reduce uncertainty in system design and
control~\cite{ST-WB-DF:05}.  Stochastic \textit{Model Predictive
  Control} (MPC) is a general framework that can handle
broad types of system uncertainty in a tractable
manner~\cite{MM-JHL:99,AM:16,DM:16,JR-DM-MD:17,HJF-HGB-MD:08}.  In general, the
effectiveness of MPC depends on the particular problem of interest, roughly
classified according to three criterion: 1) the class of systems to be
controlled, e.g., linear~\cite{MF-LG-RS:16,MP-SG-JL:12} or nonlinear
systems~\cite{AM-SS-RF-RDB:14,MAS-RRB:17,JAP-AM:19,LM-GDN-RS-FA:03}; 2) the way uncertainty is
handled, e.g., stochastic~\cite{GCC-LF:12,GS-LF-CF-MM:14,DQM-ECK-EV-PF:11}
   or bounded uncertainty~\cite{MC-BK-SVR-QC:10};
  and 3) the solution method, e.g., quadratic
programming~\cite{JR-DM-MD:17}, stochastic
programming~\cite{DP-AB-TA:05,AM-SS-RF-RDB:14}, or nonconvex
optimization~\cite{JM-ALV-JL:07,AM:16}. In practice, system
performance guarantees typically require large amounts of data
processing which, for online settings such as MPC, may become
specially challenging.

Motivated by the need of developing finite-data-driven predictive
control methods, we consider the application of a
\textit{Distributionally Robust Optimization} (DRO) framework to a
class of MPC problems.  DRO has attracted recent attention due to its
finite-sample performance guarantees~\cite{RG-AJK:16,PME-DK:17}, problem tractability via Wasserstein ambiguity sets~\cite{DB-JC-SM:20,AH-IY:20,IY:20},
its distributed formulation~\cite{AC-JC:17-allerton,AC-JC:20-tac} and
online implementation~\cite{DL-SM:19-tac-online}. These characteristics
provide a novel mechanism to deal with uncertainty in MPC, while
allowing for the tractability of the associated nonconvex optimization
problems. In practice, the efficacy of the DRO framework
depends on the explicit system structure.
Here, our proposed approach applies to a problem class whose dynamics can be  nonlinear, but affine both in control and states, and constraints can be linear or bi-linear.

For illustration purposes, we apply our solution method to toy examples
related to highway speed-limit control, focusing our discussion on its
performance with respect to uncertainty.
As the highway congestion significantly affects system operations~\cite{AAK-PV:10},
a variety of congestion mitigation strategies have been proposed,
including those based on optimization~\cite{SJ-KS:18,AJ-IP-MP-BDS:17,JY-RRB:05},
logic-based control~\cite{CS-GEA-AM_BC:16}, extremum seeking
control~\cite{HY-SK-TRO-MK:19}, and autonomous-vehicle scheduling~\cite{CW-AMB-AM:18}.
More recently, \textit{Speed-Limit Control} has been proposed as an
effective mechanism in
transportation~\cite{FS-IM-MS-MM:17,AYY-MR-MAD:17}.
In particular,
optimization-based speed-limit control has successfully demonstrated
the containment of traffic congestion. For
example,~\cite{AH-BDS-JH:03} proposed a discrete macroscopic
second-order model, METANET, and used it in optimization problems for
speed limit. However, %
at that time, the robustness analysis with respect to the system
uncertainty was absent. Later in~\cite{HY-HA-etl:17}, a linear model
predictive control of speed limit was developed for congestion
mitigation.
This work exploited the celebrated \textit{Cell Transmission Model}
(CTM)~\cite{CFD:94} or its extension for inhomogeneous traffic, the
\textit{Link Transmission Model} (LTM), to capture the deterministic
distribution of traffic densities along a highway, and the dynamic
properties of highways are characterized for congestion reduction,
with relative low online computational costs. In
addition,~\cite{SL-AS-JF-AN-EC-HH-BDS:16} proposed a scenario-based
optimization to account for the bounded uncertainty of transportation
systems.  However, an explicit treatment of the system uncertainty,
such as unknown drivers actions, vehicle arrival and departures, and
as well as random events that happen on highways, are yet to be fully
explored.

With increasing accessibility of real-time traffic
data~\cite{DBW-SB-OPT-BP-AMB:10,JCH-DBW-RH-XJB-QJ-AMB:10} for
uncertainty reduction, congestion mitigation and traffic control under
uncertainty can become practical.
As a first step into developing a novel data-driven traffic control
methods, we consider a problem related to highway transportation,
formulated by way of an extended CTM and controlled by speed limit,
and customize our proposed framework with the following question in
mind: \textit{Can we find an efficient approach for the computation of
  data-driven controls with guarantees on congestion elimination?}
{We would like to note that, while the particular problem
  we look at is inspired by traffic-congestion mitigation, the main
  emphasis of this work is on the solution approach. }
\noindent \textit{Statement of Contributions.} This work presents the
following contributions: 1) We first provide a general framework for
data-driven predictive control under uncertainty. To demonstrate our
approach explicitly, we consider a problem inspired by highway
speed-limit control which extends the CTM to account for random events
(due to drivers actions) as well as vehicle arrival and
departures. This provides an analytical framework and a stochastic
optimization problem formulation for the computation of speed limit
using the available flow measurements (in
Section~\ref{sec:ProbStat}). 2) We propose an optimization-based
data-driven control approach that extends DRO to account for system
dynamics. The resulting approach guarantees congestion elimination
with high probability, using only finite flow measurements (in
Section~\ref{sec:Reform}). 3) As the proposed data-driven optimization
problem is infinite dimensional and intractable, we propose an
equivalent reformulation that reduces the proposed problem into a
finite-dimensional optimization problem (in
Section~\ref{sec:Reform}). 4) Yet this problem is non-convex and
difficult to solve, so we find an equivalent reformulation via a
binary representation technique, and
propose a computationally efficient algorithm to provide online
sub-optimal speed limits that ensure congestion elimination and
guarantee of highway throughput with high probability (in
Section~\ref{sec:Alg}). 5) We propose in Section~\ref{sec:tools} an
optimization tool to analyze the performance of the proposed algorithm
offline. This tool is developed via a second-order cone relaxation
technique and we show that, under mild conditions, the resulting
Mixed-Integer Second-Order Cone Problem (MISOCP) is exact, and can be
handled by commercial solvers.  6) We numerically demonstrate the
effectiveness of our data-driven approach with performance guarantees,
in Section~\ref{sec:Sim} and~\ref{sec:App}. We claim that the proposed
approach is suitable for problems which are subject to control-affine
constraints.
A preliminary study of this work appeared in~\cite{DL-DF-SM:19-ecc},
however the differences are the following: This work 1) considers a general optimization framework,
2) demonstrates a detailed
procedure of the proposed integer-solution search algorithm, 3) provides an exact MISOCP relaxation as a tool for solution analysis, and 4) studies the effectiveness of the proposed approach in application. \\
\emph{Notation:} See the footnote\footnote{
  Let $\real^{m\times n}$ denote the $m \times n$-dimensional real
  vector space, and let $\vectorones{m}$ and $\vectorzeros{m}$ denote
  the column vectors $\trans{(1,\cdots,1)} \in \real^m$ and
  $\trans{(0,\cdots,0)} \in \real^m$, respectively. For any vector $x
  \in \real^m$, let us denote $x \geq \vectorzeros{m}$ if all the
  entries are nonnegative. Any letter $x$ may have appended the
  following indices and arguments: it may have the subscript $x_e$
  with $e \in \mathbb{N}$, the argument $x_e(t)$ with $t \in
  \mathbb{N}$, and finally a superscript $l\in \mathbb{N}$ as in
  $x^{(l)}_e(t)$. Given a finite number of elements $x^{(l)}_e(t) \in
  \real$ where $e, t, l\in\mathbb{N}$, we define vectors
  $x^{(l)}(t):=(x^{(l)}_1(t), x^{(l)}_2(t),\ldots)$,
  $x^{(l)}:=(x^{(l)}(1), x^{(l)}(2),\ldots)$, and $x:=(x^{(1)},
  x^{(2)},\ldots)$. Let $\left\langle x,\;y\right\rangle$ and $x \circ
  y$ denote the inner and component-wise products of vectors $x,y \in
  \real^m$, respectively. The component-wise square of vector $x \in
  \real^m$ is denoted by $x^2:=x \circ x$. In addition, let $x \otimes
  y$ denote the Kronecker product of vectors $x, y$ with arbitrary
  dimension. Let $\Norm{{x}}$ denote the $1$-norm
  of ${x} \in \real^m$, and let
  $\Norm{{x}}_{\star}:=\sup_{\Norm{z}\leq 1}\left\langle
    z,\;x\right\rangle$ denote the corresponding dual norm. Notice
  that $\|x\|_{**} = \|x\|$. \\
  \noindent \textbf{Optimization theory:} Consider a bounded function
  $\map{f}{X}{\real}$ where $X \subseteq
  \real^n$.
  The function $f$ is lower semi-continuous on $X$ if $f(x) \leq
  \liminf_{y\rightarrow x}f(y)$ for all $x \in X$. Similarly, the
  function $f$ is lower semi-continuous on $X$ if and only if its
  sublevel sets $\setdef{x\in X}{f(x)\leq \gamma}$ are closed for each
  $\gamma \in \real$. We let $\map{f^{\star}}{X}{\real\cup \{+\infty \}}$ denote the
  convex conjugate of $f$, which is defined as
  $f^{\star}(x):=\sup_{y\in X} \left\langle x,\;y\right\rangle
  -f(y)$. Further, the infimal convolution of two functions $f$ and $g$
  on $X$ is defined as $(f\square g)(x):= \inf_{y\in X} f(x-y)+g(y)$. If
  $f$ and $g$ are bounded, convex, and lower
  semi-continuous functions on
  $X$, %
  we will have $\left(f+g\right)^{\star}=(f^{\star}\square g^{\star})$.

  Consider a subset $A\subset X$, and let $\map{\chi_{A}}{X}{\real \cup \{+\infty\}}$ denote the characteristic function
  of $A$, i.e., $\chi_{A}(x)$ is equal to $0$ if and only if $x\in A$
  and $+\infty$ otherwise. In addition, let $\map{\sigma_{A}}{X}{\real}$
  denote the support function of $A$, which is defined as
  $\sigma_{A}(x):=\sup_{y\in A}\left\langle x,\;y\right\rangle$. Notice
  that $\chi_{A}$ is lower semi-continuous if and only if $A$ is closed,
  and that $\sigma_{A}(x)=[\chi_{A}]^{\star}(x)$ for all $x\in X$.
   } for basic notations of this work.

\section{Problem Formulation}\label{sec:ProbStat}
We first present our data-driven predictive control approach in a
general setting with the main goal of control design employing a finite
data set. We then focus on %
an %
one-way highway
speed-limit control as an application. This problem leverages a
traffic model based on the
Lighthill-Whitham-Richards (LWR)
discretization~\cite{MV-AH-BDS-JH:03,AAK-PV:12,CD:97,IY:07}. Finally,
we adapt the proposed predictive control approach for our control
problem, resulting into a stochastic optimization problem that can be
used to find speed limits with performance guarantees.

\subsection{General Optimization Framework for Stochastic MPC}
The goal of the proposed framework is to address system uncertainty
explicitly and find a control law which satisfies system constraints
and optimizes the expected objective in uncertainty. To achieve this,
let us denote by $t \in \mathbb{N}$, $x(t) \in \real^n$ and $u(t)\in
\real^m$ the time, system state and control at time $t$,
respectively. We consider the control law $u$ to be solutions of the
following $T$-long receding-horizon stochastic predictive control
problem: {\leqnomode
\begin{align} \small
  \label{eq:P0} \tag{$\mathbf{P}$}~&\sup\limits_{u} \;
  \mathbb{E}_{\prob(u)}\left[ \ell(u , x) \right],  \\
  \st \; & x \sim \prob(u) \textrm{ characterized by} \nonumber \\
  & x(t+1)=F(x(t),u(t),\xi(t)), \; t=0,1,\ldots, T-1, \nonumber \\
  &  x \in \mathcal{Z}(u), \; u \in \mathcal{U}, \; \xi \sim \prob_{\xi}, \nonumber
\end{align}}where $u$ is a concatenated variable of
$u(0),u(1),\ldots,u(T-1)$ and $x$ is that of those $x(t)$. Notice that
$x$ is a stochastic process and we denote by $\prob(u)$ the
distribution of $x$ given $u$. The objective is to maximize the
expectation of a given reward function $\map{\ell}{\real^{mT} \times
  \real^{nT}}{\real}$ taken under $\prob(u)$. We denote
by %
$\map{F}{\real^n \times \real^m \times \real^d}{\real^n}$ the given
system dynamics where $\xi$ represents the uncertainty. We assume that
the process $x$ is constrained in a given set $\mathcal{Z}(u)$ and so
does $u$ in $\mathcal{U}$. We denote by $\prob_{\xi}$ the unknown
distribution of the uncertainty process.

Due to the unknown $\prob_{\xi}$, Problem~\eqref{eq:P0} cannot be
solved exactly. To find a control law that solves~\eqref{eq:P0}, we
propose a distributionally robust optimization (DRO)
framework. %
By means of this, we will employ a finite set of realizations or
samples of the random variables $\xi$ to approximate the unknown
distribution and compute a set of feasible $u$. The main appeal of
this robust method is that it can provide \textit{out-of-sample}
probabilistic guarantees of performance~\cite{PME-DK:17}. In
particular, at each $t$, let us assume that $N$ samples of $x(0)$ and
$\xi:=(\xi(0),\ldots,\xi(T-1))$ are accessible. Under some mild
conditions on these samples, we will show in Section~\ref{sec:Reform}
that a tractable optimization-based function $J(u)$ can be constructed
using those state and uncertainty samples. The function $J(u)$ is a
surrogate objective of~\eqref{eq:P0} which accounts for the system
dynamics as well as constraints on the states. In addition, given a
confidence value $\beta \in (0,1)$, we provide the
\textit{out-of-sample} performance guarantee in the sense that the
probability of the true objective function in~\eqref{eq:P0} being
greater than ${J}({u})$ is greater than $1-\beta$. In other words, we
have
\begin{equation*} \small
{\Prob^N}\left( \mathbb{E}_{\prob(u)}\left[ \ell(u , x) \right]  \geq
{J}({u}) \right)\geq 1- \beta,
\end{equation*}
where $\Prob^N$ is the product probability over $N$ sample
trajectories of the system. Thus, with high probability, the choice of
samples to approximate the problem will provide a minimum lower value
for the original problem. As Problem~\eqref{eq:P0} needs to be solved
in a moving horizon fashion, the surrogate functions $J(u)$ is
optimized similarly, as in Algorithm~\ref{Alg:ddpc}. Notice that, the functions $J(u)$ depend on samples as wells as the system structure, and the solution to $J(u)$ is nontrivial.
\begin{algorithm}
\caption{Data-driven predictive control with guarantees} \label{Alg:ddpc}
\begin{algorithmic}[1]
  \State Initialize $t=0$
  \While {True}
  \State Take $N$ %
  measurements $x^{(l)}(t)$ and $\xi^{(l)}$, $l=1,\ldots,N$
  \State Adapt an approach to optimize $J(u)$ over $u$
  \State Apply performance-guaranteed $u$ to the system
  \State $t \leftarrow t+1$
  \EndWhile
\end{algorithmic}
\end{algorithm}

To enable the proposed tractable optimization method for $J(u)$ as in Section~\ref{sec:Alg}, we assume the following system structures.
\begin{assumption}[Control-and-state-affine systems with bi-linear constraints] \label{assump:structure}
The system dynamics $F(x,u,\xi)$ is continuous, affine in ${x}$ and affine in $u$. In addition, the set $\mathcal{Z}(u)$ is composed of constraints which are linear and bi-linear in $(x,u)$. And the set $\mathcal{U}$ is a finite set.
\end{assumption}
Assumption~\ref{assump:structure} covers a wide class of dynamical systems, including linear systems and bi-linear systems. Furthermore, with a minor modification of the proposed reformulation techniques in Section~\ref{sec:equivReform}, the proposed solution method can be extended to control-affine systems. We leave the extension as the future work. Next, we consider a traffic speed-limit control
problem and design a set of speed controls using the proposed framework.

\subsection{Highway Traffic Model}
\begin{figure}[tbp]%
\centering
\includegraphics[width=0.35\textwidth]{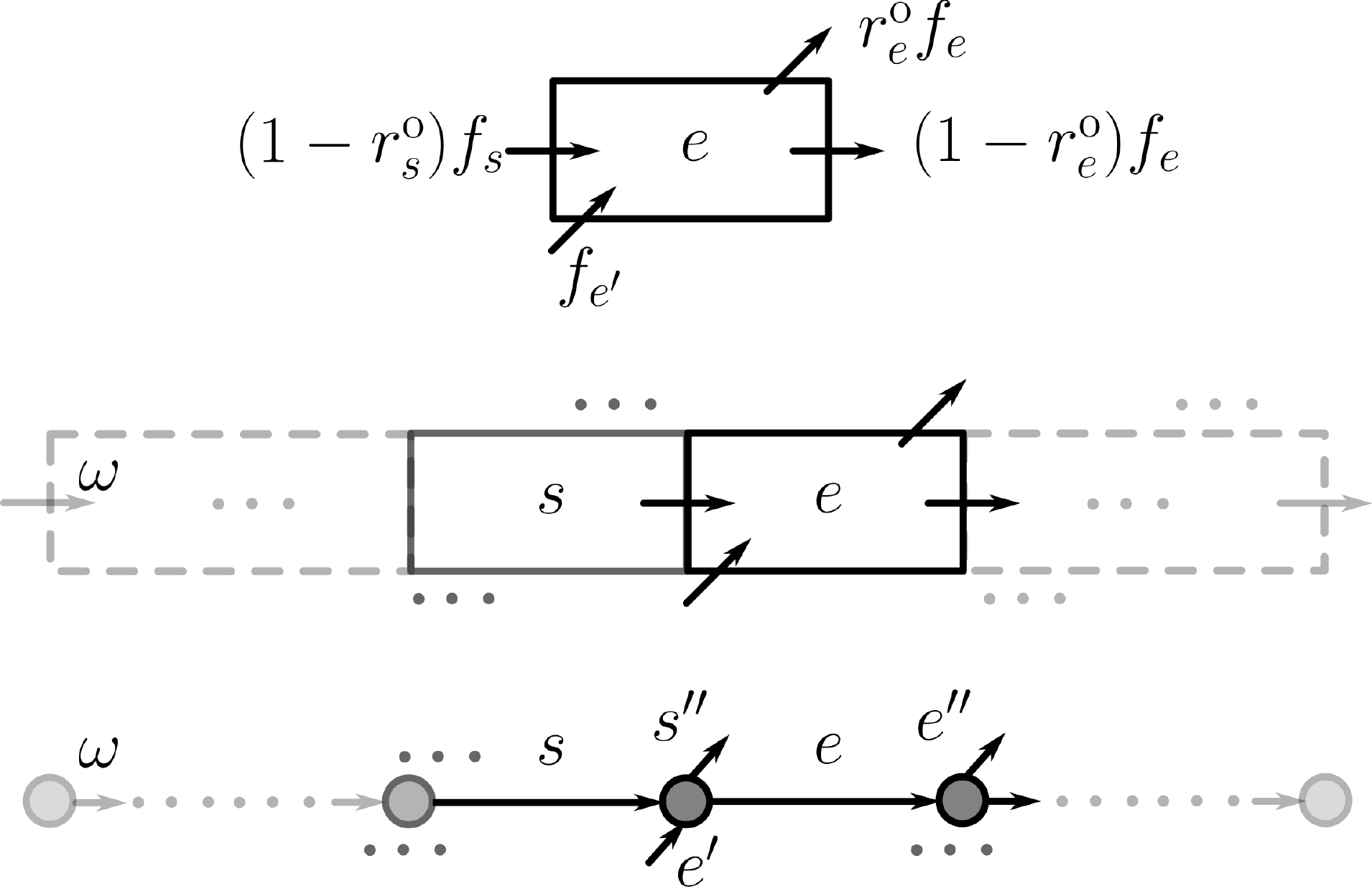}%
\caption{\small Segment of highways and its graph representation. The
  symbol $e$ indicates a highway segment, and $s$ is its preceding
  segment. The on-ramp and off-ramp of $e$ are denoted by $e^{\prime}$
  and $e^{\prime \prime}$, respectively. The variable $\omega$, $f$
  and $\supscr{r}{o}_e$ represents the mainstream inflow, edge flows and
  fraction of outflows, respectively.}%
\label{fig:highway}%
\end{figure}

\reviseone{In order to
  introduce our discrete traffic model, we start by describing our
   time and space discretization, respectively. } \reviseone{Let
  $\delta$ be a sufficiently small time-discretization step},
then %
we let $\mathcal{T}=\until{T-1}$ denote the set of time slots, where
\reviseone{we identify time by index $t \equiv t\delta$. }%
\reviseone{Further, let us} consider a one-way highway of length $L$,
and divide \reviseone{the highway} into $n$
segments. %
The topology of the highway can be described by a directed graph
$\graph=(\nodes,\edges)$, where $\nodes=\{0,1,\ldots,n\}$ and each
\reviseone{node} $v\in \nodes$ corresponds to the \reviseone{junction}
between two consecutive road segments. The fictitious node
$0$ represents the \reviseone{mainstream} inflow into the highway.
The graph edge set is $\edges=\{(0,1),\cdots,(n-1,n)\}$ with each
element $e=(v,v+1)\in \edges$ corresponding to the road segment
between nodes $v$ and $v+1$ for all $v\in\{0,1,\ldots,n-1\}$. To
illustrate this, \reviseone{the reader is referred to
  Fig. \ref{fig:highway}.}  \reviseone{Nodes $0$ and $n$ are the
  source and sink nodes of the graph $G$, respectively}. Further, node
$v \in \nodes \setminus \{0 \}$ is called an arrival node if there
exists an on-ramp at node $v$. Similarly, node $v\in \nodes \setminus
\{n \}$ is called a departure node if there exists an off-ramp at node
$v$. Let $\nodes_{A}$ and $\nodes_{D}$ denote the set of arrival and
departure nodes, respectively. By convention, we set node $0 \notin
\nodes_{A}$ and node $n \notin \nodes_{D}$. For each edge
$e=(v,v+1)\in \edges$ with starting node $v$, let $e^{\prime}$ denote
the on-ramp of edge $e$ if node $v \in \nodes_{A}$ and let $e^{\prime
  \prime}$ denote the off-ramp of edge $e$ if node $v+1 \in
\nodes_{D}$.  \reviseone{For each $e \in \edges$, we denote by
  $\len_e$ and $\lane_e$ the segment length and the number of lanes on
  $e$, respectively. In particular, the highway length
  $L=\sum_{e\in\edges}\len_e$.}

\reviseone{At each time $t\in\mathcal{T}$, macroscopic
  models of traffic use aggregated variables to describe the behavior
  of traffic. We start presenting constraints that relate the state
  variable $\rho_e(t)$, the average traffic density in segment $e$ and
  period $t$, with the average traffic flow $f_e(t)$ on segment $e$,
  as well as the average traffic  velocity $s_e(t)$ on segment $e$.
}

\reviseone{A first relationship is given by the traffic-flow
  definition
\begin{equation*}
  f_e(t)=s_e(t) \rho_e(t) , %
\end{equation*}
where $s_e(t)$ is taken to be a function of speed limit as
  follows.} Let $u_e(t)$ denote the speed limit set on
  $e$ %
  and \reviseone{assume that the majority of drivers comply with
  it. Then,
\begin{equation*}
  s_e(t)= \min \{ \overline{s}_e(\rho_e(t)), u_e(t)  \},
\end{equation*}
where $\overline{s}_e(\rho_e(t))$ is the maximal admissible %
speed on $e$ %
when the traffic density is $\rho_e(t)$. %
In practice, the non-negative values of $\overline{s}_e$
monotonically decrease as $\rho_e(t)$ increases, reflecting the
safe-driving behavior.
}

\reviseone{Further, the fundamental diagram is a plot of the nonlinear
  relationship between traffic flow and traffic density at a location
  of a highway.  %
  In this plot, the density at which the traffic flow attains its
  maximum value, $\overline{f}_e$, is called the \textit{critical
    density}, $\rho^{c}_e$. The density at which the traffic
  flow is zero is the traffic jam density
  $\overline{\rho}_e$.}
\reviseone{The traffic flow is an increasing function of $\rho_e$ on
  $(0, \rho^{c}_{e})$, and a strictly decreasing function of $\rho_e$
  on $(\rho^{c}_{e}, \overline{\rho}_e)$. }

\begin{figure}[tbp]%
\centering
\includegraphics[width=0.3\textwidth]{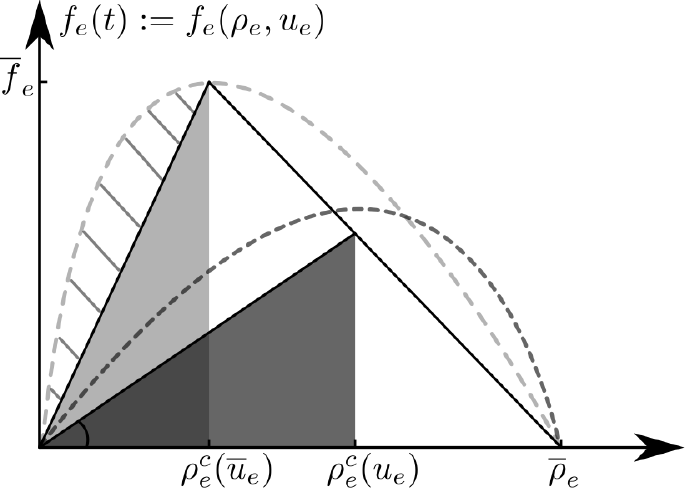}%
\caption{\small {Flow rate as a function of traffic density of edge
    $e$ for two speed limits $\overline{u}_e$ and $u_e$ such that
    $\overline{u}_e \ge u_e$. The two nonlinear curves are the
    fundamental diagrams corresponding to each speed limit, while the
    straight lines are piece-wise linear approximations of these. The
    slope of linear approximations at the origin represents the speed
    limits $\overline{u}_e$ and $u_e$, respectively. The light and
    dark shaded region guarantees no congestion of edge $e$ under
    speed limit $\overline{u}_e$ and $u_e$, respectively. (That is,
    under speed limit $u_e$, a density such that $\rho_e(t) \leq
    \rho^{c}_{e}({u}_{e})$, for all $t$, guarantees no congestion.)
    \reviseone{The slashed area reflects the compliance of drivers
      with speed limit $\overline{u}_e$. \reviseone{A higher
        compliance of drivers results in a smaller area.}}  }}%
\label{fig:fd}%
\end{figure}

\reviseone{How speed limits affect the fundamental diagram
  has been a subject of debate~\cite{FS-IM-MS-MM:17,AYY-MR-MAD:17}.
    Following~\cite{FS-IM-MS-MM:17}, we will assume that in the
  presence of speed limits, the flow rate and traffic density
  still hold a similar relationship, however the
  critical density will be a function of the velocity limit $u_e(t)$.
  Let $\overline{u}_e$ be the free flow velocity corresponding to a
  maximum value of flow under no speed limits.\footnote{\reviseone{Given drivers behavior
      $\overline{s}_e$, the value $\overline{u}_e=\max_{\rho_e}
      \overline{s}_e(\rho_e) $, $\overline{\rho}_e$ admits
      $\overline{s}_e(\overline{\rho}_e)=0$ and $\overline{f}_e$
      maximizes $f_e$ over any $u_e$ and $\rho_e \in (0,
      \overline{\rho}_e)$. On the other hand, the function
      $\overline{s}_e$ is highly dependent on}
  the physical structure of \reviseone{the segment $e$} as well as
  random events, such as accidents and temporary lane closures (see,
  e.g.,~\cite{AAK-PV:12,KCD-LY-XW-YW-HS-MC-LY-CQ-VS:15}).} A comparison of two
  fundamental diagrams with constant speed limits $\overline{u}_e$ and
  $u_e$ are shown in Fig.~\ref{fig:fd}. Notice that the reduction of
  the speed limit from $\overline{u}_e$ to $u_e$ increases the
  critical density and decreases the maximal flow rate on edge $e$.}

\reviseone{To model the fundamental diagram in the presence of speed
  limits, consider edge $e\in \edges$, and let
  $\rho^{c}_{e}(u_e(t))$ denote the critical density of
  edge $e$ at speed limit $u_e(t)$.} \reviseone{That is,
  the critical density $\rho^{c}_{e}({u}_e(t))$ determines
  the traffic density at which the maximum edge flow
    $f_e(t)$ is achievable}. \reviseone{Given speed limit
  $u_e(t)$, %
  we will work with an approximation of the fundamental diagram of
  edge $e\in \edges$ given as follows }
\begin{equation}
  f_{e}(t)=
	\begin{cases}
          \reviseone{u_e(t)} \rho_e(t) , & \; {\rm{if}} \; \rho_{e}(t) \leq \rho^{c}_{e}(\reviseone{u_e(t)}), \\
          \tau_e\overline{u}_e\left( \overline{\rho}_e - \rho_{e}(t)
          \right), & \; {\rm{otherwise},}
	\end{cases}
  \label{eq:approx}
\end{equation}
with\footnote{The parameter $\tau_e\overline{u}_e$ is known as the
  backward wave speed (see, e.g.,~\cite{MJL-GBW:55}). }
\begin{equation*}
  \rho^{c}_{e}(\reviseone{u_e(t)}):=\left({\tau_e\overline{\rho}_{e}\overline{u}_e}\right)/
  \left({\tau_e \overline{u}_{e} + \reviseone{u_e(t)} }\right),
\end{equation*}
where the parameter $\tau_e:=
{\overline{f}_e}/{\left(\overline{u}_e\overline{\rho}_e
    -\overline{f}_e\right)}$. To illustrate this, we refer the reader
to Fig. \ref{fig:fd}.

Each segment $e$ is \reviseone{congested}
when its density is higher
than its critical density (see,
e.g.,~\cite{GG-RH-AAK-PV-JK:08,AYY-MR-MAD:17}). Since critical
densities are determined by speed limits, we \reviseone{select the speed limits} such that
the following constraint is satisfied at all time slots $t \in
\mathcal{T}$
\begin{equation}
 0 \leq \rho_e(t) \leq \rho^{c}_{e}(\reviseone{u_e(t)}),~\forall e\in\edges.
  \label{eq:RhoC}
\end{equation}
The constraint above ensures that the highway is not congested
regardless of the flow rates.  Using the constraint above and the
fundamental diagram approximation in~\eqref{eq:approx}, we obtain
\begin{equation}\label{eq:fund}
f_e(t)=\reviseone{u_e(t)}\rho_e(t),~\forall e\in\edges,~t \in \mathcal{T}.
\end{equation}

\reviseone{Finally, based on the physical principle of conservation of
  mass, the discretized LWR model provides a set of difference
  equations for each road segment $e$ which enable us to analyze the
  dynamics of traffic flows on highways~\cite{CS-AM:15,IY:07}.  For
  each segment $e=(v,v+1) \in \edges$, the dynamics of
  $\rho_e$ are determined by the \textit{demand} flow
  $\supscr{f}{D}_e(t)$ from $v$ and the \textit{supply} flow
  $\supscr{f}{S}_e(t)$ to $v+1$ as follows:
\begin{equation*}
  \rho_e(t+1)= \rho_e(t)+ h_e(\supscr{f}{D}_e(t)- \supscr{f}{S}_e(t)),
  \quad \forall t \in\mathcal{T}, %
\end{equation*}
where
$h_e:=\delta/\len_e$.
For numerical stability, $\delta$  must be
selected such that %
$h_e  \leq 1/ \max_{t\in\mathcal{T}} \{u_e(t)\}$, $\forall e\in \edges$~\cite{SJ-KS:18}.
}

\reviseone{The supply flow $\supscr{f}{S}_e(t)$ denotes the maximal flow that can be transferred through edge $e$, and is given by %
\begin{equation*}
  \supscr{f}{S}_e(t)=
	\begin{cases}
          u_e(t) \rho_e(t) , & \; {\rm{if}} \; \rho_{e}(t) \leq \rho^{c}_{e}(u_e(t)), \\
          u_e(t) \rho^{c}_{e}(u_e(t)), & \; {\rm{otherwise},}
	\end{cases}
\end{equation*}
Notice that $\supscr{f}{S}_e(t) = f_e(t)$ when the highway segment $e$ is not congested, i.e., the constraint in \eqref{eq:RhoC} is satisfied.}  %
\reviseone{Consider junction $v \in \nodes_{D}$, i.e., the junction $v$ is a departure node of the preceding edge of $e$ or that of some edge $s$. Since $v \in \nodes_{D}$, a fraction of the
{supply $\supscr{f}{S}_{s}(t)$}
will depart the highway through an off-ramp edge $s^{\prime \prime}$ and
the rest will enter into succeeding segment $e$. Let
$\supscr{r}{o}_s(t) \in [0,1)$ denote the fraction of the supply $\supscr{f}{S}_{s}(t)$
that departs the system. Hence, the flow through the off-ramp edge
{$s^{\prime \prime}$ is $f_{s^{\prime
    \prime}}(t):=\supscr{r}{o}_s(t)\supscr{f}{S}_{s}(t)$}.} Notice that the fraction
$\supscr{r}{o}_s(t)$ is determined by the drivers' actions. Therefore,
$\supscr{r}{o}_s$ is random, and its value is unknown to the system
operator in advance. Each random variable $\supscr{r}{o}_s(t)$ will
have a nonempty support set denoted by
$\mathcal{Z}_{\supscr{r}{o}_s(t)} \subset \realnonnegative$.

\reviseone{The traffic demand $\supscr{f}{D}_e(t)$ depends on the supply of its preceding edge $s\in\edges$ as well as the ramp flows on their connected junction $v$. At each $v \in \nodes_{A}$, a fraction of the traffic demand $\supscr{f}{D}_e(t)$ is originated from the on-ramp edge $e^{\prime}$. Let $\supscr{r}{in}_e(t) \in [0,1)$ denote the fraction of the traffic demand $\supscr{f}{D}_e(t)$ originated from the on-ramp edge $e^{\prime}$. Hence, the on-ramp traffic flow is given by  $f_{e^{\prime}}(t):=\supscr{r}{in}_e(t)\supscr{f}{D}_e(t)$. Notice that the ratio $\supscr{r}{in}_e(t)$ is an exogenous parameter that depends on the traffic flow at the on-ramp edge $e^{\prime}$.} Hence, each ratio $\supscr{r}{in}_e(t)$ can be
modeled as a random variable with nonempty support
$\mathcal{Z}_{\supscr{r}{in}_e(t)} \subset \realnonnegative$\footnote{\reviseone{For $v \notin \nodes_{D}$ or $v \notin \nodes_{A}$}, the value of
$\supscr{r}{o}_s(t)$ or $\supscr{r}{in}_e(t)$ is zero, respectively.}.
\reviseone{Then, by the conservation of flows, at each time slot $t$, the traffic demand $\supscr{f}{D}_e(t)$ must satisfy the following constraint:
\begin{equation*}
\begin{aligned}
\supscr{f}{D}_e(t)&=\supscr{f}{S}_s(t) - f_{s^{\prime
    \prime}}(t) + f_{e^{\prime}}(t) , \quad \forall \; t \in \mathcal{T}. %
\end{aligned}
\end{equation*}
At edge $e=(0,1)$, we have  $\supscr{f}{D}_{e}(t):=\omega(t)$ which is a random mainstream with support $\mathcal{Z}_{\omega(t)} \subset \realnonnegative$.
At each edge $e$, the demand $\supscr{f}{D}_e(t)$ must be %
admissible to %
edge $e$, %
i.e.,
\begin{equation}
  \supscr{f}{D}_e(t) \leq
    \min\{ \overline{f}_e, \; \tau_e\overline{u}_e\left( \overline{\rho}_e - \rho_{e}(t) \right)  \} .
  \label{eq:demand}
\end{equation}
Notice that, %
the constraints~\eqref{eq:demand} allows for transient speed of $\supscr{f}{D}_e(t)$ higher than speed limit $u_e(t)$, as long as the mean speed $s_e(t)$ complies with $u_e(t)$.}

Let $\rho(0)=(\rho_1(0),\ldots,\rho_n(0))$
denote the traffic density of the highway with support
$\mathcal{Z}_{\rho(0)} \subset \realnonnegative^n$.
Using the constraints above, the traffic density dynamics at each time
slot $t \in \mathcal{T}$ %
are \reviseone{given by}
\begin{equation}
\begin{aligned}
  \rho_e(t+1)&= \rho_e(t)+ \reviseone{h_e} \frac{1-\supscr{r}{o}_s(t)}{1-\supscr{r}{in}_e(t)}f_s(t)- \reviseone{h_e} f_e(t), \\
	& \hspace{3cm}  \quad \forall e \in \edges \setminus \{(0,1) \}, \\
	  \rho_e(t+1)&= \rho_e(t)+ \reviseone{h_e}( \omega(t) - f_e(t) ), \;  e =(0,1) . \\
 \end{aligned} \label{eq:CTM}
\end{equation}
Recall that $f_e(t)=\reviseone{u_e(t)}\rho_e(t)$. %
\reviseone{For each $e \in \edges \setminus \{ (0,1) \}$ and  $t\in\mathcal{T}$, %
the constraint in~\eqref{eq:demand} can be written as
\begin{equation}
    \frac{1-\supscr{r}{o}_s(t)}{1-\supscr{r}{in}_e(t)}f_s(t) \leq  \min\{ \overline{f}_e, \; \tau_e\overline{u}_e\left( \overline{\rho}_e - \rho_{e}(t) \right)  \} ,
  \label{eq:demand_new}
\end{equation}
where $s$ is the preceding edge of edge $e$.
}

Our goal is to design \reviseone{a set of} speed limits for \reviseone{drivers}. To achieve this goal,
we consider a finite set of speed limits, and then approximate the
fundamental diagram of each segment $e$ with a finite set of piece-wise
linear functions, as shown in Fig.~\ref{fig:fd}. Let $\Gamma$ be a
\reviseone{finite set of feasible speed limits} for the \reviseone{highway segments}. More precisely,
the speed limit of each edge $e \in \edges$ must satisfy
\begin{equation}
\begin{aligned}
 \reviseone{u_e(t)} \in \Gamma:=& \{\gamma^{(1)}, \ldots, \gamma^{(m)}\}, \quad \reviseone{t\in \mathcal{T}}. %
\end{aligned}	\label{eq:ue}
\end{equation}
The set of real values $\Gamma$ is determined by the physical
structure of the highway as well as its maximal free flow speed and
traffic jam density. \reviseone{As mentioned earlier, random
events, such as traffic incidents or lane closure}, can change these values.

\subsection{Problem Formulation for the Traffic Control Problem}
We aim at maximizing \reviseone{the expected flow rate of highway
  segments while reducing congestion via speed
  limits.} %
To compute \reviseone{the} average flow, let $\prob_{\varpi}$ denote
the distribution of the concatenated random variable $\varpi:=(\omega,
\rho(0),\supscr{r}{in},\supscr{r}{o})$. \reviseone{Given the
  parameters $\{\overline{f}_e\}_{e\in\edges}$,
  $\{\overline{\rho}_e\}_{e\in\edges}$, and
  $\Gamma$,} %
the problem of computing speed limits which
\reviseone{maximize the expected flow}, can be formulated as follows:
{\leqnomode
\begin{align}
  \label{eq:P} \tag{$\mathbf{P}$}~ \max\limits_{ \substack{u, \rho }}
  \quad & {\mathbb{E}_{\prob_\varpi}\left[ \frac{1}{T} \sum_{e\in\edges,t\in\mathcal{T} }{\rho_e(t) \reviseone{u_e(t)}} \right] } ,  \\
  \st \quad & \eqref{eq:RhoC},\; \eqref{eq:fund}, \; \eqref{eq:CTM},
  \; \eqref{eq:demand_new}, \; \eqref{eq:ue},\; \nonumber
\end{align}}where %
\reviseone{$\rho$ and $u$
  are %
  the concatenated variables of
  $\{\rho_e(t)\}_{e\in\edges,t\in\mathcal{T}
  }$\footnote{$\rho:=(\rho_1(0),\rho_2(0), \ldots,\rho_n(0),\rho_1(1),
    \ldots,\rho_n(T-1))$.} and
  $\{u_e(t)\}_{e\in\edges,t\in\mathcal{T}}$,
  respectively.}

Problem~\eqref{eq:P} is \reviseone{nonconvex}.  In addition, the
probability distribution $\prob_\varpi$ is unknown, i.e., it is
impossible to compute \reviseone{the expected flow (i.e., the
  objective function)} exactly. Our goal is to compute a
  set of speed limits that are feasible to
  Problem~\eqref{eq:P} and guarantee a minimum achievable expected
  flow in the presence of uncertainty on $\prob_\varpi$.  To achieve
  this goal, we adapt the proposed %
  framework to compute the desired speed limits.
  In this way, given an optimal speed limit $u$ and a set
    of $N$ samples, we obtain $J(u)$, an achievable average flow
    rate---let us call it certificate. In particular, $J(u)$ is a
    function of the $N$ random samples. This certificate is a minimum
    with confidence $\beta \in (0,1)$ in the sense that the
    probability of the true objective function being greater than
    ${J}({u})$ is greater than $1-\beta$. In other words, let
    $\Prob^N$ denote the product probability distribution over $N$
    samples. Under certain conditions on $\prob_{\varpi}$, the proposed %
    approach
    guarantees that the following out-of-sample performance constraint
    is satisfied: %
\begin{equation} \small
  {\Prob^N}\left({\mathbb{E}_{\prob_\varpi}
    \left[ \frac{1}{T}  \sum_{e\in\edges,t\in\mathcal{T}}{\rho_e(t) \reviseone{u_e(t)}}\right]} \geq
    {J}({u}) \right)\geq 1- \beta.
\label{eq:perf}
\end{equation}
The probabilistic guarantee~\eqref{eq:perf} enables us to evaluate the performance of a feasible solution $u$ to Problem ~\eqref{eq:P} via only finite samples of $\varpi$ and a lower bound $J(u)$.
We call a solution $u$ with the probabilistic  guarantee~\eqref{eq:perf}  \textit{a data-driven Speed-Limit control}.

\begin{figure}[tbp]%
\centering
\includegraphics[width=0.45\textwidth]{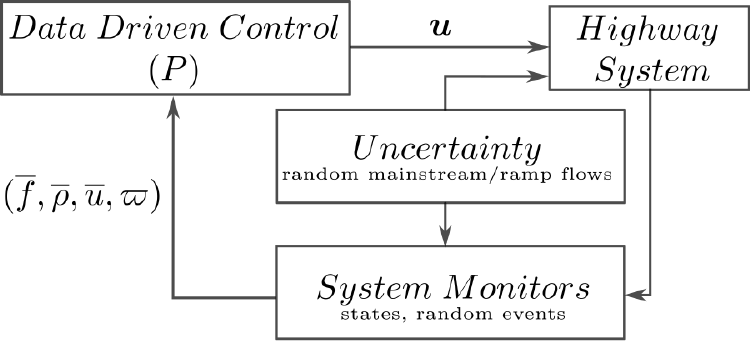}%
\caption{\small \reviseone{Data-driven speed-limit control on highways under random ramp flows and system events.}}%
\label{fig:control_framework}%
\end{figure}

\begin{remark}[Implementation of data-driven control]
{\rm \reviseone{With online-accessible samples of $\varpi$ and parameters $(\overline{f},\overline{\rho},\overline{u})$ which are related to real-time highway random events, the data-driven control $u$ can be achieved via online solutions to~\eqref{eq:P} in a moving horizon fashion. Fig.~\ref{fig:control_framework} demonstrates how to implement the proposed approach in real-world applications.
}}\end{remark}
\begin{remark}[Admissible operation zone]
{\rm In this work, we propose a set of speed-limit controls that prevents congestion with probabilistic guarantees. Note that the existence of such controls is highly dependent on the feasibility of Problem (P). When traffic demands are higher than the highway capacity in a sufficiently long time, congestion is inevitable. Therefore, there is no hope to prevent congestion via the speed-limit control in the presence of high traffic demands. In such scenarios, we set speed limits to a set of predefined values.}
\end{remark}

\section{Performance Guaranteed Speed-Limit Design}\label{sec:Reform}
Our goal is to compute a set of speed limits with certain
out-of-sample performance guarantees. To achieve this goal, we follow
a four-step procedure. First, we reformulate Problem~\eqref{eq:P} into
an equivalent problem (call it Problem~\eqref{eq:P1}). Second, we
propagate the admissible sample trajectories using the $N$ measurements of
$\varpi$. Third, we adopt a distributionally robust optimization
approach to~\eqref{eq:P1}. The first three steps enable us to
obtain a distributionally robust optimization framework for computing
speed limits with guarantees equivalent to~\eqref{eq:perf}. Finally, we obtain a
tractable problem reformulation.

\noindent \textbf{Step 1: (Equivalent Reformulation of~\eqref{eq:P})}
Traffic densities $\rho_e(t)$, for all $e\in \edges$ and $t \in
\mathcal{T}$, are random since traffic inflows and
outflows are random in each segment $e$. Using this observation, we
now consider the variable $\rho$ in Problem~\eqref{eq:P} as a random
variable, and derive an equivalent Problem~\eqref{eq:P1} via a
reformulation of the constraints in~\eqref{eq:P}.

Given a speed limit $u$ that satisfies the constraint~\eqref{eq:ue},
let $\prob(u)$ and $\mathcal{Z}(u)$ denote the probability
distribution of variable $\rho$ and the support of $\rho$,
respectively.\footnote{%
  The support $\mathcal{Z}(u)$ %
  is the smallest closed set such that $P(\rho \in
  \mathcal{Z}(u))=1$.}
Recall that in Problem~\eqref{eq:P}, constraints~\eqref{eq:RhoC}
ensure no congestion on the highway. %
Hence, the given $u$ should be
such that $\mathcal{Z}(u) \subseteq \setdef{\rho \in
  \real^{nT}}{\eqref{eq:RhoC}}$.
Without loss of generality, we consider \reviseone{the largest possible support} $\mathcal{Z}(u):=\setdef{\rho
\in \real^{nT}}{\eqref{eq:RhoC}}$. To fully characterize random
variable $\rho$, the distribution of $\prob(u)$ needs to be
determined. Using the traffic density dynamics in~\eqref{eq:CTM},
flow representation in~\eqref{eq:fund} \reviseone{and sustainability constraints in~\eqref{eq:demand_new}}, the probability distribution
$\prob(u)$ can be represented as a convolution of the distribution
$\prob_{\varpi}$.

Let ${\mathcal{M}}(\mathcal{Z}(u))$ denote the space of all
probability distributions supported on $\mathcal{Z}(u)$. \reviseone{Let us} write the objective function of Problem~\eqref{eq:P} compactly as
\begin{equation*}
  H(u;{\rho}):= \frac{1}{T}
  \sum_{e\in\edges,t\in\mathcal{T}}{\rho_e(t) \reviseone{u_e(t)}},
\label{eq:H}
\end{equation*}
and reformulate~\eqref{eq:P} as follows
{\leqnomode
\begin{align}
  \label{eq:P1} \tag{$\mathbf{P1}$}~\max\limits_{u} \quad &
  \mathbb{E}_{\prob(u)}\left[ H(u;{\rho}) \right] , \\
  \st \quad & \prob(u) \textrm{ characterized by~\eqref{eq:fund},~\eqref{eq:CTM},~\eqref{eq:demand_new} and } \prob_{\varpi}, \nonumber \\
  & \prob(u) \in {\mathcal{M}}(\mathcal{Z}(u)), \;
  \eqref{eq:ue}. \nonumber
\end{align}} Notice that Problems \eqref{eq:P} and \eqref{eq:P1} are
equivalent. \reviseone{Hence}, we obtain the performance guarantee
of~\eqref{eq:P1} by considering the induced out-of-sample performance
on $\mathbb{P}(u)$, written as $\mathbb{E}_{\prob(u)} \left[
  H(u;{\rho}) \right]$. For all problems derived later, we use
the performance guarantees equivalent to~\eqref{eq:perf}, as follows
\begin{equation}
  {\Prob^N}\left(\mathbb{E}_{\prob(u)}
    {\left[ H(u;{\rho}) \right] }\geq
    {J}({u}) \right)\geq 1- \beta.
\label{eq:perfnew}
\end{equation}
Recall that
$\Prob^N$ denotes the probability that the event $\mathbb{E}_{\prob(u)}
{\left[ H(u;{\rho}) \right] } \geq {J}({u})$ happens on the $N$
product of the sample space that defines $\rho$, the value $J(u)$ is
the certificate to be determined, and $\beta \in (0,1)$ is
the desired confidence value.

Next, we characterize the probability distribution $\prob(u)$ using
the $N$ sample measurements of $\varpi$.

\noindent \textbf{Step 2: (Admissible Sample Trajectory Propagation)} Let
$\mathcal{L}= \untilone{N}$ denote the index set for the $N$
realizations of the random variable ${\varpi}$, and let
$\{\varpi^{(l)}\}_{l\in\mathcal{L}}$, where
$\varpi^{(l)}:=(\omega^{(l)},
\rho^{(l)}(0),\supscr{r}{in,(\textit{l})},\supscr{r}{o,(\textit{l})}
)$, denote the set of independent and identically distributed (i.i.d.)
realizations of $\varpi$. Given a speed limit $u$ and measurement
$\varpi^{(l)}$, a unique traffic density trajectory ${\rho}^{(l)}$ can
be computed by using~\eqref{eq:fund},~\eqref{eq:CTM}. Notice that this trajectory is
unique since density dynamics~\eqref{eq:CTM} is linear (in $\rho$) for a given
speed limit $u$ and measurement $\varpi^{(l)}$. \reviseone{Further,
the resulting ${\rho}^{(l)}$ is an \textit{admissible} traffic trajectory if the given $u$ achieves the flow sustainability constraints~\eqref{eq:demand_new}.}
Given these
realizations $\{\varpi^{(l)}\}_{l\in\mathcal{L}}$, the admissible sample
trajectories $\{{\rho}^{(l)}\}_{l\in\mathcal{L}}$ of the random
traffic flow densities for each edge $e \in \edges$ with its precedent $s\in \edges\cup \varnothing$, are given by
\begin{equation} \small
\begin{aligned}
&  \rho^{(l)}_e(t+1)= \reviseone{h_e}
  \frac{1-\supscr{r}{o,(\textit{l})}_s(t)}{1-\supscr{r}{in,(\textit{l})}_e(t)} \reviseone{u_s(t)} {\rho}^{(l)}_s(t)- \reviseone{h_e} \reviseone{u_e(t)} {\rho}^{(l)}_e(t)  \\
  &  \hspace{2.8cm} \quad +\rho^{(l)}_e(t),  \quad \forall e \in \edges \setminus \{(0,1) \}, \\
&  \rho^{(l)}_e(t+1)= \rho^{(l)}_e(t)+ \reviseone{h_e}( \omega^{(l)}(t) - \reviseone{u_e(t)} {\rho}^{(l)}_e(t) ), \;  e =(0,1) , \\
& \reviseone{\frac{1-\supscr{r}{o,(\textit{l})}_s(t)}{1-\supscr{r}{in,(\textit{l})}_e(t)} u_s(t) {\rho}^{(l)}_s(t) \leq  \min\left\{ \overline{f}_e, \; \tau_e\overline{u}_e\left( \overline{\rho}_e - \rho^{(l)}_{e}(t) \right)  \right\}}, \\
& \hspace{5cm} \reviseone{\forall e \in \edges \setminus \{(0,1) \}} ,  \\
 \end{aligned} 	\label{eq:Rhohat}
\end{equation}
for all $t \in \mathcal{T}$ %
and $l \in
\mathcal{L}$. The following lemma establishes that
$\{{\rho}^{(l)}\}_{l\in\mathcal{L}}$ are i.i.d. samples from
$\prob(u)$.

\begin{lemma}[Independent and identically distributed sample
  generators of ${\rho}$] Given a speed limit $u$ and a set of
  i.i.d. realizations of $\varpi$, the system
  dynamics~\eqref{eq:Rhohat} generate i.i.d. admissible sample trajectories
  $\{{\rho}^{(l)}\}_{l\in\mathcal{L}}$ of $\prob(u)$.
\label{lemma:RhoGener}
\end{lemma}
The proof is provided in %
Appendix~\ref{appx:whole}, which also contains all the proofs of the later lemmas, propositions and theorems.

Let $\subscr{\mathcal{M}}{lt}(\mathcal{Z}_{\varpi}) \subset
{\mathcal{M}}(\mathcal{Z}_{\varpi})$ denote the space of all
light-tailed probability distributions supported on
$\mathcal{Z}_{\varpi}$\footnote{For any set $\mathcal{Z}$, we use
  notion $\subscr{\mathcal{M}}{lt}(\mathcal{Z})$ to denote the space
  of all light-tailed probability distributions supported on
  $\mathcal{Z}$.}. We make the following assumption on the probability
distribution $\prob_{\varpi}$:
\begin{assumption}[Light-tailed unknown distributions]
  The  distribution $\prob_{\varpi}$ satisfies $\prob_{\varpi} \in
  {\subscr{\mathcal{M}}{lt}}(\mathcal{Z}_{\varpi})$, i.e., there
  exists an exponent $a>1$ such that $b:= \mathbb{E}_{\prob_{\varpi}}
  \left[ \exp(\Norm{\varpi}^a) \right] < \infty$. \label{assump:1}
\end{assumption}
\begin{remark}[Accessible light-tailed i.i.d. samples of
  $\varpi$] \label{remark:iid_sample}{\rm\reviseone{The random
      variable $\varpi$ essentially represents the random flows and
      densities on the highway. This results into an unknown compact
      support of $\prob_{\varpi}$, indicating that $\prob_{\varpi}$
      is also light-tailed. Further, online samples of $\varpi$ can come
      from various independent system monitors, e.g., Traffic
      Performance Measurement Systems (T-PeMS), or real-time GPS
      systems, which provide online i.i.d. samples of
      $\varpi$. %
        }}\end{remark}
The following lemma establishes that the probability distribution of
traffic densities is light-tailed when $\prob_{\varpi}$ is
light-tailed.

\begin{lemma}[Light-tailed distribution of ${\rho}$] Let
  Assumption~\ref{assump:1} hold, then $\prob(u) \in
  \subscr{\mathcal{M}}{lt}(\mathcal{Z}{(u)})$.
\label{lemma:hatRho}
\end{lemma}
The above Lemmas~\ref{lemma:RhoGener} and~\ref{lemma:hatRho} %
enable the application of the %
proposed framework to~\eqref{eq:P1} in the
next step.

\noindent \textbf{Step 3: (Performance Guarantee Certificate)} Given a
speed limit ${u}$, we design a certificate ${J}({u})$ that satisfies
the performance guarantee condition~\eqref{eq:perfnew}. To achieve
this goal, the proposed %
approach consists of solving a robust
  version of the problem over a set of distributions. In particular,
we will consider a set of distributions $\mathcal{P}({u})$ that is
small, tractable and yet rich enough to contain ${\prob}(u)$ with high
probability. Then by evaluating~\eqref{eq:P1} under \reviseone{the}
worst-case distribution in $\mathcal{P}({u})$, \reviseone{the}
performance of~\eqref{eq:P1} in the form of~\eqref{eq:perfnew} can be
guaranteed.

Consider the Wasserstein ball\footnote{Let
${\mathcal{M}}(\mathcal{Z})$ denote the space of all probability
distributions supported on $\mathcal{Z}$. Then for any two
distributions $\mathbb{Q}_1$, $\mathbb{Q}_2 \in
\mathcal{M}(\mathcal{Z})$, the Wasserstein metric~\cite{KLV-RGS:58}
$\map{d_W}{\mathcal{M}(\mathcal{Z}) \times
  \mathcal{M}(\mathcal{Z})}{\realnonnegative}$ is defined by
\begin{equation*}
  d_W(\mathbb{Q}_1,\mathbb{Q}_2)
  := \min\limits_{\Pi} \int_{\mathcal{Z} \times \mathcal{Z}} \Norm{\xi_1 -\xi_2} \Pi(d \xi_1, d \xi_2),
\end{equation*}
where $\Pi$ is in a set of all distributions on $\mathcal{Z} \times
\mathcal{Z}$ with marginals $\mathbb{Q}_1$ and $\mathbb{Q}_2$. A
closed Wasserstein ball of radius $\epsilon$ centered at a
distribution $\prob \in \mathcal{M}(\mathcal{Z})$ is denoted by
$\mathbb{B}_{\epsilon}(\prob):=\setdef{\mathbb{Q} \in
  \mathcal{M}(\mathcal{Z})}{d_W(\prob,\mathbb{Q}) \le \epsilon}$.
} %
$\mathbb{B}_{\epsilon}(\hat{\prob}({u}))$ of center $\hat{\prob}(u):=
({1}/{N})\sum_{l\in\mathcal{L}} \delta_{\{{\rho}^{(l)}\}}$ and radius
$\epsilon$. Notice that the center $\hat{\prob}(u)$ is obtained using
the point mass operator $\delta$ under i.i.d. admissible
sample trajectories
$\{{\rho}^{(l)}\}_{l\in\mathcal{L}}$. \reviseone{These trajectories
  are} distributed according to ${\prob}(u)$ which is a function of
the controlled system dynamics~\eqref{eq:Rhohat} and samples of
$\varpi$. Then we propose certificates ${J}({u})$ of~\eqref{eq:P1} in
the following theorem.
\begin{theorem}[Performance guarantees
  of~\eqref{eq:P1}] \label{thm:cert} \reviseone{Let us assume that $N$
    i.i.d. samples of ${\varpi}$ are given, together with a speed
    limit $u$, and confidence value $\beta$, and that
    Assumption~\ref{assump:1} on light-tailed distributions of
    ${\varpi}$ holds. Further, let us define the set of distributions
    $\mathcal{P}({u})$ as follows:}
\begin{equation*}
  \mathcal{P}({u}) := \mathbb{B}_{\epsilon(\beta)}(\hat{\prob}({u})) \cap  \subscr{\mathcal{M}}{lt}(\mathcal{Z}{(u)}).
\end{equation*}
Then, there exists a Wasserstein
  radius $\epsilon:=\epsilon(\beta)$, depending on the confidence value
  $\beta$ and Assumption~\ref{assump:1}, such that \reviseone{${\prob}(u)$ is in the set of distributions $\mathcal{P}({u})$} as
  described in~\eqref{eq:P1} with probability at least $1-\beta$,
  i.e.,
\begin{equation*}
  {\Prob^N}\left( \prob(u) \in \mathcal{P}(u) \right)\geq 1- \beta.
\end{equation*}
Further, the following proposed certificate ${J}({u})$ of~\eqref{eq:P1} satisfies guarantee~\eqref{eq:perfnew}
\begin{equation*}
  {J}({u}):= \inf_{\mathbb{Q} \in \mathcal{P}({u})}
  {\mathbb{E}_{\mathbb{Q}} \left[H(u;{\rho})\right]}.
\end{equation*}
\end{theorem}

\begin{remark}[Effect of Assumption~\ref{assump:1} on
  ${J}({u})$] \label{remark:radius} {\rm The certificate ${J}({u})$
    highly depends on the set $\mathcal{P}({u})$ and
    the Wasserstein radius $\epsilon(\beta)$. In
    Theorem~\ref{thm:cert}, the value $\epsilon(\beta)$ is calculated
    via the parameters $a$ and $b$ in Assumption~\ref{assump:1}. As
    these parameters may not be known, one can determine
    $\epsilon(\beta)$ in a data driven fashion via \reviseone{Monte-Carlo} simulations. That is, we start by setting $\epsilon(\beta)=0$
    and gradually increase it until the performance
    guarantees~\eqref{eq:perfnew} hold with that given $\beta$ for a
    sufficiently large number of simulation runs.  }\end{remark}

Theorem~\ref{thm:cert} provides each feasible solution $u$
of~\eqref{eq:P1} with a certificate ${J}({u})$ that guarantees
performance as in~\eqref{eq:perfnew}. This motivates a tractable
reformulation of~\eqref{eq:P1} as follows.

\noindent \textbf{Step 4: (Tractable Reformulation of~\eqref{eq:P1})}
Our goal is to obtain a speed limit $u$ that maximizes the average
flow through the highway as in~\eqref{eq:P1} while ensuring that the
performance guarantee condition~\eqref{eq:perfnew} is satisfied.
Given a set of \reviseone{inflow-and-outflow-related} samples
$\{\varpi^{(l)}\}_{l\in\mathcal{L}}$, a speed limit ${u}$ that
provides the highest ${J}({u})$ \reviseone{(i.e., the best objective
  lower bound)}, can be computed by solving the following optimization
problem: {\leqnomode
\begin{align}
  \label{eq:P2} \tag{$\mathbf{P2}$}~\sup\limits_{u \; \st \;
    {\eqref{eq:ue}}} \quad & {J}({u}).
\end{align}} \reviseone{Problem~\eqref{eq:P2} is an
infinite-dimensional optimization problem, and, hence, it is hard to
solve.} %
\reviseone{The following theorem}, \reviseone{provides a
  finite-dimensional reformulation of} Problem~\eqref{eq:P2},
\reviseone{called~\eqref{eq:P3}, and shows} that
problems~\eqref{eq:P2} and~\eqref{eq:P3} are equivalent for
$({u},{J})$.

\begin{theorem}[Tractable reformulation of~\eqref{eq:P2}]
\reviseone{Consider}
  {\leqnomode
  \begin{align}
      \label{eq:P3} \tag{$\mathbf{P3}$}~~~~
      \max\limits_{u,{\rho},\lambda,\mu,\nu,\eta } \quad & -\lambda
      \epsilon(\beta) -
      \frac{1}{N}\sum_{e\in\edges,t\in\mathcal{T},l\in\mathcal{L}}{\overline{f}_e\overline{\rho}_e\eta_e^{(l)}(t)  } \nonumber \\
      & \hspace{2.9cm} \quad+\frac{1}{N}\sum_{l\in\mathcal{L}}{\left\langle \nu^{(l)} , {\rho}^{(l)}\right\rangle}, \nonumber \\
      \st \quad & %
      \reviseone{ [\overline{f} \otimes
      \vectorones{T} +
      (\overline{\rho}-\rho^{c}(\overline{u})) \otimes
      \vectorones{T} \circ u  ] \circ \eta^{(l)} }
      \nonumber \\
      & \hspace{2.9cm} -\mu^{(l)} \geq \vectorzeros{nT},
      \; \forall \; l\in\mathcal{L}, \nonumber \\
      & \nu^{(l)} =\mu^{(l)}+ \frac{1}{T}u , %
      \; \forall \; l\in\mathcal{L}, \nonumber \\
      & \Norm{\nu^{(l)}}_{\star} \leq \lambda, \; \forall
      \; l\in\mathcal{L}, \nonumber\\
      & \eta^{(l)} \geq \vectorzeros{nT}, \; \forall \; l\in\mathcal{L},\nonumber \\
      & \eqref{eq:ue},\; \eqref{eq:Rhohat},\nonumber
    \end{align}} where decision variables
  $(u,{\rho},\lambda,\mu,\nu,\eta)$ are concatenated versions of $ \reviseone{u_e(t)}$,
  ${\rho}^{(l)}_e(t)$, $\lambda$, $\mu^{(l)}_e(t)$, $\nu^{(l)}_e(t)$,
  $\eta^{(l)}_e(t) \in \real$, for all $l \in\mathcal{L}$, $t \in
  \mathcal{T}$, %
  and $e \in \edges$. The value
  $\rho^{c}(\overline{u}):=\overline{f}/\overline{u}$ is the vector of
  critical densities under the free flow, and $\overline{\rho}$ is the
  jam density vector.

  Problem~\eqref{eq:P2} is equivalent to~\eqref{eq:P3} in the sense
  that their optimal objective values coincide and the set of
  optimizers of~\eqref{eq:P2} are the projection of that
  of~\eqref{eq:P3}.  Further, for any feasible point
  $(u,{\rho},\lambda,\mu,\nu,\eta)$ of~\eqref{eq:P3}, let $\hat{J}(u)$
  denote the value of its objective function. Then the pair
  $({u},\hat{J}(u))$ gives a data-driven solution $u$ with an estimate
  of its certificate ${J}(u)$ by $\hat{J}(u)$, such that the
  performance guarantee~\eqref{eq:perfnew} holds for $(u,\hat{J}(u)
  )$.
\label{thm:P2}
\end{theorem}

\begin{remark}[Formulation~\eqref{eq:P3} depends on
  data] \label{remark:data_in_P}{\rm\reviseone{Problem~\eqref{eq:P3}
      is parameterized by the variables $(\overline{f},
      \overline{\rho},\overline{u})$, which are related to the highway
      infrastructure and random events, and the data
      $\{\varpi^{(l)}\}_{l\in\mathcal{L}}$, which are related to the
      traffic system initial state, $\{\rho^{(l)}(0)\}_{ l
          \in \mathcal{L}}$, on-ramp, and off-ramp
      flows. The decision variables include future states
        of the density $\rho^{(l)}(t)$, $t>0$, the speed limits, and
        other multipliers to make the constraints hold. In particular,
        note that the ambiguity-ball constraint is enforced via the
        $\lambda$ multiplier and maximization in $\rho$. In practice,
      highway random events such as accidents are monitored in real
      time, which provides particular information about
      $(\overline{f}, \overline{\rho},\overline{u})$. On the other
      hand, data $\{\varpi^{(l)}\}_{l\in\mathcal{L}}$ is accessible
      from various independent monitors as mentioned in the previous
      remark on data-driven control implementation. }}\end{remark}

Problem~\eqref{eq:P3} is inherently difficult to solve due to the
discrete decision variables $u$, bi-linear terms $u %
\circ \eta^{(l)}$ in the first group of constraints,
and the nonlinear admissible sample trajectories
$\{{\rho}^{(l)}\}_{l\in\mathcal{L}}$, which motivates our next
section.

\section{\reviseone{Equivalent} Reformulation}
\label{sec:equivReform} Our goal is
to compute exact solutions to Problem~\eqref{eq:P3}. To achieve this
goal,
we focus on the
feasibility set of Problem~\eqref{eq:P3}, and transform \reviseone{a group of} its non-convex
quadratic terms, \reviseone{which are comprised of a continuous variable and a binary variable,} into a set of mixed-integer linear constraints. We
call the new \reviseone{equivalent} formulation, Problem~\eqref{eq:P4}.
\noindent \textbf{Binary Representation of Speed Limits:}
Let
$\mathcal{O}:=\untilone{m}$ be the index set of the speed limit set
$\Gamma$. For each edge $e\in\edges$, \reviseone{time slot $t \in \mathcal{T}$} and speed limit value
$\gamma^{(i)} \in \Gamma$, let us define the binary variable $\reviseone{x_{e,i}(t)}$
to be equal to one if $\reviseone{u_e(t)}=\gamma^{(i)}$; otherwise $\reviseone{x_{e,i}(t)}=0$. We
will then have $\reviseone{u_e(t)}=\sum_{i\in \mathcal{O}} \gamma^{(i)} \reviseone{x_{e,i}(t)}$ for
each edge $e \in\edges$. Using this representation, we can reformulate
the speed limit constraints~\eqref{eq:ue} as follows

\begin{equation}
	\begin{aligned}
          &\gamma^{(1)} \leq  \sum_{i \in \mathcal{O}} \gamma^{(i)}
          \reviseone{x_{e,i}(t)} \leq \gamma^{(m)}, \; \forall e \in \edges, \; \forall \reviseone{t \in \mathcal{T}}, \\
          &\sum_{i \in
            \mathcal{O}}\reviseone{x_{e,i}(t)}=1,
          \forall e \in \edges, \; \forall \reviseone{t \in \mathcal{T}}, \\
          &\reviseone{x_{e,i}(t)} \in \{0,\; 1\}, \; \forall e \in \edges, \;
          i \in \mathcal{O}, \; \reviseone{t \in \mathcal{T}}, \\
	\end{aligned} 	\label{eq:ueNew}
\end{equation}
and we update the admissible sample trajectories formula~\eqref{eq:Rhohat} for
all $t\in\mathcal{T}$ %
and $l\in\mathcal{L}$ as follows
\begin{equation} \small
\begin{aligned}
&  \rho^{(l)}_e(t+1)= \reviseone{h_e}
 \frac{1-\supscr{r}{o,(\textit{l})}_s(t)}{1- \supscr{r}{in,(\textit{l})}_e(t)}
  \sum_{i\in \mathcal{O}} \gamma^{(i)} \reviseone{x_{s,i}(t)} {\rho}^{(l)}_s(t) +\rho^{(l)}_e(t) \\
  &  \hspace{1.4cm} \quad  - \reviseone{h_e} \sum_{i\in \mathcal{O}} \gamma^{(i)} \reviseone{x_{e,i}(t)} {\rho}^{(l)}_e(t),  \quad \forall e \in \edges \setminus \{(0,1) \}, \\
  & \rho^{(l)}_e(t+1)= \rho^{(l)}_e(t)+ \reviseone{h_e}\omega^{(l)}(t) - \reviseone{h_e} \sum_{i\in \mathcal{O}} \gamma^{(i)} \reviseone{x_{e,i}(t)} {\rho}^{(l)}_e(t) ,  \\
  & \hspace{5.8cm} \quad  e =(0,1). \\
  & \reviseone{\frac{1-\supscr{r}{o,(\textit{l})}_s(t)}{1-\supscr{r}{in,(\textit{l})}_e(t)} \sum_{i\in \mathcal{O}} \gamma^{(i)} x_{s,i}(t) {\rho}^{(l)}_s(t) \leq } \\
  & \hspace{1cm} \reviseone{\min\left\{ \overline{f}_e, \; \tau_e\overline{u}_e\left( \overline{\rho}_e - \rho^{(l)}_{e}(t) \right)  \right\}, \; \forall e \in \edges \setminus \{(0,1) \}}, \\
 \end{aligned} 	\label{eq:RhohatNew}
\end{equation}
\reviseone{We} are particularly interested in two groups of
bi-linear terms: 1) the bi-linear terms $\reviseone{x_{e,i}(t)} {\rho}^{(l)}_e(t)$ in
the admissible sample trajectories formula~\eqref{eq:RhohatNew} and 2) the
bi-linear terms $\sum_{i\in \mathcal{O}} \gamma^{(i)}
\reviseone{x_{e,i}(t)}\eta^{(l)}_e(t)$ which appear in the first set of constraints,
e.g., $u %
\circ \eta^{(l)}$.  Each of these
bi-linear terms is comprised of a continuous variable and a binary
variable. We represent each of these bi-linear terms with a set of
linear constraints using the following linearization technique.

Let us introduce variables $y^{(l)}_{e,i}(t)$ and $z^{(l)}_{e,i}(t)$
for all $e \in \edges$, $i\in \mathcal{O}$, $t\in \mathcal{T}$
and $l \in \mathcal{L}$ as follows
\begin{equation}\label{eq:bilinearNew}
\begin{aligned}
&y^{(l)}_{e,i}(t)=\reviseone{x_{e,i}(t)}\rho^{(l)}_e(t),\\
&z^{(l)}_{e,i}(t)=\reviseone{x_{e,i}(t)}\eta^{(l)}_e(t).
\end{aligned}
\end{equation}
We further make the following assumption
\begin{assumption}[Bounded dual variable $\eta$]
  There exists a positive constant $\overline{\eta}$ such that for any optimizers of~\eqref{eq:P3}, the components $\eta^{(l)}_e(t) \leq \overline{\eta}$ for all $e \in \edges$, $t
  \in \mathcal{T}$ %
  and $l
  \in\mathcal{L}$. \label{assump:2}
\end{assumption}
We achieve Assumption~\ref{assump:2} by selecting $\overline{\eta}$ large enough. This enables the following lemma to represent the non-convex
equality constraint in~\eqref{eq:bilinearNew} with a set of linear
constraints.

\begin{lemma}[Linearization technique] \label{lemma:linearization} Let
  Assumption~\ref{assump:2}
   hold.  Then for all $e
  \in \edges$, $i\in \mathcal{O}$, $t\in \mathcal{T}$ %
  and $l \in \mathcal{L}$, the non-convex equality constraint
  in~\eqref{eq:bilinearNew} can be equivalently represented with the
  following set of linear constraints
\begin{equation}
\begin{array}{l}
  0\leq z^{(l)}_{e,i}(t) \leq \overline{\eta} \reviseone{x_{e,i}(t)}, \\
  \eta^{(l)}_e(t) -\overline{\eta}(1- \reviseone{x_{e,i}(t)}) \leq z^{(l)}_{e,i}(t)
  \leq \eta^{(l)}_e(t), \\
\label{eq:z}
\end{array}
\end{equation}
\begin{equation}
	\begin{array}{l}
  0\leq y^{(l)}_{e,i}(t) \leq \overline{\rho}_e \reviseone{x_{e,i}(t)}, \\
          {\rho}^{(l)}_e(t) -
          \overline{\rho}_e(1- \reviseone{x_{e,i}(t)}) \leq y^{(l)}_{e,i}(t)
          \leq {\rho}^{(l)}_e(t) . \\
 \end{array}
\label{eq:y}
\end{equation}
\end{lemma}
In particular, we \reviseone{let} %
$y^{(l)}_{e,i}(0)=\reviseone{x_{e,i}(0)}{\rho}^{(l)}_e(0)$ for each $e \in \edges$, $i\in \mathcal{O}$ and $l \in \mathcal{L}$. Then using the new variables $y^{(l)}_{e,i}(t)$ and $z^{(l)}_{e,i}(t)$, we can now reformulate the admissible sample trajectories formula~\eqref{eq:RhohatNew} as follows

\begin{equation} \small
\begin{aligned}
&  \rho^{(l)}_e(t+1)= \reviseone{h_e} \frac{1-\supscr{r}{o,(\textit{l})}_s(t)}{1-\supscr{r}{in,(\textit{l})}_e(t)}
	\sum_{i\in \mathcal{O}} \gamma^{(i)} y^{(l)}_{s,i}(t) +\rho^{(l)}_e(t) \\
	 &  \hspace{1.4cm} \quad  - \reviseone{h_e} \sum_{i\in \mathcal{O}} \gamma^{(i)} y^{(l)}_{e,i}(t),  \quad \forall e \in \edges \setminus \{(0,1) \}, \\
&	  \rho^{(l)}_e(t+1)= \rho^{(l)}_e(t)+ \reviseone{h_e}\omega^{(l)}(t) - \reviseone{h_e} \sum_{i\in \mathcal{O}} \gamma^{(i)} y^{(l)}_{e,i}(t) ,  \\
		& \hspace{5.5cm} \quad  e =(0,1). \\
    & \reviseone{\frac{1-\supscr{r}{o,(\textit{l})}_s(t)}{1-\supscr{r}{in,(\textit{l})}_e(t)} \sum_{i\in \mathcal{O}} \gamma^{(i)} y^{(l)}_{s,i}(t) \leq } \\
    & \hspace{1cm} \reviseone{\min\left\{ \overline{f}_e, \; \tau_e\overline{u}_e\left( \overline{\rho}_e - \rho^{(l)}_{e}(t) \right)  \right\}, \; \forall e \in \edges \setminus \{(0,1) \}}, \\
 \end{aligned} 	\label{eq:RhohatNew2}
\end{equation}
Problem~\eqref{eq:P3} can now be equivalently reformulated
as the following optimization problem
{\leqnomode
\begin{align} \small
  \label{eq:P4}\tag{$\mathbf{P4}$}~ \max\limits_{\substack{x,y,z,
      {\rho}, \\ \lambda,\mu,\nu,\eta }} -\lambda \epsilon(\beta)
  - \frac{1}{N}\sum_{e,t,l}{\overline{f}_e\overline{\rho}_e\eta_e^{(l)}(t)}
    + \frac{1}{N}\sum_{e,t,l}{\nu^{(l)}_e(t){\rho}^{(l)}_e(t)},
\end{align}}
{\reqnomode
\begin{align}
  \st \quad & \sum_{i\in \mathcal{O}} \gamma^{(i)}
  (\overline{\rho}-\rho^{c}(\overline{u}))\otimes
  \vectorones{T} \circ z^{(l)}_i  -\mu^{(l)}
  \nonumber \\
  & \hspace{1.8cm} +\overline{f}\otimes \vectorones{T} \circ
  \eta^{(l)} \geq \vectorzeros{nT}, \; \forall \;
  l\in\mathcal{L}, \label{eq:dual1} \\
  & \nu^{(l)} =\mu^{(l)}+ \frac{1}{T}\sum_{i\in \mathcal{O}} \gamma^{(i)}x_{i}
  ,\; \forall \;
  l\in\mathcal{L}, \label{eq:dual2} \\
  & \Norm{\nu^{(l)}}_{\star} \leq \lambda, \; \forall
  \; l\in\mathcal{L}, \label{eq:dual3} \\
  & \vectorzeros{nT} \leq \eta^{(l)} \leq \overline{\eta},
  \; \forall \; l\in\mathcal{L},\label{eq:dual4} \\
  & { \textbf{\small speed limits}}\; \eqref{eq:ueNew},
  \;{\textbf{\small dual variable}} \; \eqref{eq:z}, \nonumber \\
  & {\textbf{\small sample trajectories}} \{
  \eqref{eq:y}, \; \eqref{eq:RhohatNew2} \}. \nonumber
\end{align}}
\begin{remark}[Performance guarantee~\eqref{eq:perfnew} in the setting
  of Problem~\eqref{eq:P4}] \label{remark:P4} {\rm Let $\hat{J}(u)$
    denote the value of the objective function of~\eqref{eq:P4} at a
    computed feasible solution
    $(x,y,z,{\rho},\lambda,\mu,\nu,\eta)$. Then, the resulting speed
    limits $u:=\sum_{i\in \mathcal{O}} \gamma^{(i)} \reviseone{x_i}
    $ provide a data-driven solution such
    that $(u,\;\hat{J}(u))$ satisfies the performance
    guarantee~\eqref{eq:perfnew} of~\eqref{eq:P1}. This result
    exploits the fact that Problem~\eqref{eq:P4} is equivalent
    to~\eqref{eq:P3} and~\eqref{eq:P3} is equivalent to~\eqref{eq:P2}
    as in Theorem~\ref{thm:P2}.  }
\end{remark}

\section{Computationally Efficient Algorithms} \label{sec:Alg}
\reviseone{We propose} a decomposition-based, \reviseone{integer-solution search algorithm} \reviseone{which} computes \reviseone{online-tractable,} high-quality feasible solutions
to \eqref{eq:P4} with performance guarantees.
Similar algorithms have
been proposed in the literature~\cite{LX-TA-BP:11,DL-XL:16}. Such
methods allow us to compute sub-optimal solutions to mix-integer
nonlinear programs efficiently.
The proposed integer-solution search
algorithm is shown in Algorithm~\ref{Alg:issa}. This algorithm
iteratively computes sub-optimal solutions to \eqref{eq:P4} until a
stopping criteria is met. At each iteration, the algorithm solves an
upper-bounding problem to~\eqref{eq:P4}, and then solves a
lower-bounding problem to~\eqref{eq:P4}. The upper-bounding
Problem~\eqref{eq:UBPk} is obtained through McCormick relaxations of
the bi-linear terms $\{\nu^{(l)} \circ {\rho}^{(l)}\}_{l
  \in\mathcal{L}}$. This upper bounding problem is a mixed-integer
linear program and its solution provides us with an upper bound on
Problem~\eqref{eq:P4} and a candidate speed limits
$x^{(k)}$. \reviseone{Notice that $x^{(k)}$ respects the sustainability constraints in~\eqref{eq:RhohatNew2}}.
We then use the computed speed limit $x^{(k)}$ to construct
a set of admissible sample trajectories $\{{\rho}^{(l,k)}\}_{l \in \mathcal{L}}$
and equivalently reduce Problem~\eqref{eq:P4} to a linear lower-bounding Problem~\eqref{eq:LBPk} for potential
feasible solutions of~\eqref{eq:P4}. If~\eqref{eq:LBPk} is feasible, then the candidate speed limits together with the objective value of~\eqref{eq:LBPk} provide %
guarantee~\eqref{eq:perfnew} for~\eqref{eq:P4}. Next, we present the
upper-bounding and lower-bounding problems in detail.

\begin{algorithm}
\caption{Integer solution search algorithm} \label{Alg:issa}
\begin{algorithmic}[1]
  \State Initialize $k=0$ \Repeat \State $k \leftarrow k+1$ \State
  Solve Problem~\eqref{eq:UBPk}, \Return $x^{(k)}$ and $\UB_k$ \State
  Generate admissible sample trajectories $\{{\rho}^{(l,k)}\}_{l \in
    \mathcal{L}}$ \State Solve Problem~\eqref{eq:LBPk}, \Return
  $\obj_{k}$ and $\LB_k$ \Until{${\UB}_k -{\LB}_k \leq \epsilon$,
    or~\eqref{eq:UBPk} is infeasible, or a satisfactory sub-optimal solution is found after certain running time $\subscr{T}{run}$}\\
  \Return data driven solution $\subscr{u}{best}:=u^{(q)}$ with
  certificate $\hat{J}(u^{(q)})$ such that $q \in
  \argmax_{p=1,\ldots,k} \{ {\obj}_p\}$
\end{algorithmic}
\end{algorithm}

\subsection{Upper-bounding Problem}
Problem~\eqref{eq:UBPk} is constructed in two stages:

\noindent \textbf{Stage 1:} We use a standard McCormick relaxation to
handle the non-convex quadratic terms
$\{\nu^{(l)}_e(t){\rho}^{(l)}_e(t)\}_{e\in\edges,t\in\mathcal{T},l\in\mathcal{L}}$
in the objective function of~\eqref{eq:P4}. Notice that the McCormick
envelope~\cite{GPM:76} provides relaxations of bi-linear terms, which
is stated in the following remark.

\begin{remark}[McCormick envelope]
{\rm  Consider two variables $x, \; y\; \in \real$ with upper and lower
  bounds, $\underline{x} \leq x \leq \overline{x}$, $\underline{y}
  \leq y \leq \overline{y}$. The McCormick envelope of the variable
  $s:=xy \in \real$ is characterized by the following constraints
\begin{equation*}
  \begin{aligned}
    s\geq \overline{x}y+ x\overline{y}-\overline{x}\overline{y},
    \quad&
    s\geq \underline{x}y+x\underline{y}-\underline{x}\underline{y}, \\
    s\leq
    \overline{x}y+x\underline{y}-\overline{x}\underline{y},\quad&s\leq
    \underline{x}y+x\overline{y}-\underline{x}\overline{y}.
  \end{aligned}
\end{equation*}
}\end{remark}

To construct a McCormick envelope for~\eqref{eq:UBPk}, let us denote
$\overline{\nu}_e:= \overline{u}_{e} \left(T^{-1} +
  \overline{\rho}_{e}\overline{\eta} \right)$ for each edge
$e\in\edges$. We have $0 \leq \nu^{(l)}_e(t) \leq
\overline{\nu}_e$, $0 \leq {\rho}^{(l)}_e(t) \leq
\overline{\rho}_e$ for all $e \in\edges$, $t \in\mathcal{T}$, and
$l \in\mathcal{L}$. Therefore, the McCormick envelope of $s^{(l)}_e(t):=
\nu^{(l)}_e(t){\rho}^{(l)}_e(t)$ is given by
\begin{equation}
	\begin{aligned}
          s^{(l)}_e(t) & \geq \overline{\nu}_e {\rho}^{(l)}_e(t)+
          \nu^{(l)}_e(t)\overline{\rho}_e-\overline{\nu}_e
          \overline{\rho}_e, \\
          s^{(l)}_e(t) & \geq 0, \\
          s^{(l)}_e(t)&\leq \overline{\nu}_e {\rho}^{(l)}_e(t),\\
          s^{(l)}_e(t)&\leq \nu^{(l)}_e(t)\overline{\rho}_e.
	\end{aligned} \label{eq:McR}
\end{equation}

\noindent \textbf{Stage 2:} We identify appropriate canonical integer
cuts to prevent~\eqref{eq:UBPk} from choosing examined candidate variable
speed limits $\{ x^{(p)} \}_{p=1}^{k-1}$. Let
$\Omega^{(p)}:=\setdef{(e,i,\reviseone{t}) \in \edges \times
  \mathcal{O} \reviseone{\times \mathcal{T}} }{\reviseone{x^{(p)}_{e,i}(t)}=1}$ denote the index set of $x$ for which
the value $\reviseone{x^{(p)}_{e,i}(t)}$ is $1$ at the previous iteration $p$. In
addition, let $c^{(p)}:=|\Omega^{(p)}|$ and $\overline{\Omega}^{(p)}:=
\left(\edges \times \mathcal{O} \reviseone{\times \mathcal{T}}\right)\setminus \Omega^{(p)}$ denote
the cardinality of the set $\Omega^{(p)}$ and the complement of
$\Omega^{(p)}$, respectively. Therefore, the canonical integer cuts of
Problem~\eqref{eq:UBPk} at iteration $k$ are given by
\begin{equation}
  \begin{aligned}
    &\sum_{ (e,i,t) \in \Omega^{(p)} } \reviseone{x_{e,i}(t)} -
    \sum_{ (e,i,t) \in \overline{\Omega}^{(p)} } \reviseone{x_{e,i}(t)} \leq c^{(p)} -1, \\
    &\hspace{4.5cm} \forall p \in \untilone{k-1}.
	\end{aligned} \label{eq:cut}
\end{equation}

Upper-bounding Problem~\eqref{eq:UBPk} can be formulated as follows
{\leqnomode
\begin{align}
  \label{eq:UBPk} \tag{{UBP}$_k$}~&\max\limits_{\substack{x,y,z,s,
      {\rho}, \\ \lambda,\mu,\nu,\eta }}
  -\lambda \epsilon(\beta) - \frac{1}{N}\sum_{e,t,l}{\left(\overline{f}_e\overline{\rho}_e\eta_e^{(l)}(t)  - s^{(l)}_e(t) \right)}, \\
  \st \quad & {\textbf{\small speed limits}}\; \eqref{eq:ueNew}, \; {\textbf{\small sample trajectories }}  \{\eqref{eq:y}, \; \eqref{eq:RhohatNew2} \} , \nonumber\\
  &  {\textbf{\small no congestion}} \;  \{\eqref{eq:z},\; \eqref{eq:dual1}, \;\eqref{eq:dual2}, \;\eqref{eq:dual3}, \;\eqref{eq:dual4}\}, \nonumber \\
  & {\textbf{\small McCormick envelope }} \;\eqref{eq:McR}, \;
  {\textbf{\small integer cuts}} \; \eqref{eq:cut}. \nonumber
\end{align}} Let $\UB_k$ denote the optimal objective value
of~\eqref{eq:UBPk}, and let $x^{(k)}$ denote the integer part of the
optimizers of~\eqref{eq:UBPk}. Then $\UB_k$ is an upper bound of the
original non-convex Problem~\eqref{eq:P4}. We use $x^{(k)}$ as a
candidate speed limit in the lower-bounding problem LBP$_k$.

\subsection{Lower-bounding Problem}
Problem~\eqref{eq:P4} can be equivalently written as {\leqnomode
  \begin{align}  ~\max\limits_{\substack{x,y,z, {\rho}, \\
        \lambda,\mu,\nu,\eta }} & -\lambda \epsilon(\beta) -
    \frac{1}{N}\sum_{e,t,l}{ \left( \overline{f}_e
        \overline{\rho}_e\eta_e^{(l)}(t)  - \nu^{(l)}_e(t){\rho}^{(l)}_e(t) \right)}, \nonumber\\
    \st \quad & (z,\lambda,\mu,\nu,\eta) \in \Phi(x), \;
    (y,{\rho})\in\Psi(x),\; x\in X. \nonumber
  \end{align}}
where
\begin{align}
&\Phi(x):=\setdef{(z,\lambda,\mu,\nu,\eta)}{\textbf{no congestion}
},\nonumber\\
&\Psi(x):=\setdef{(y,{\rho})}{\textbf{sample
    trajectories} },\nonumber\\
&X:=\setdef{x}{\textbf{speed limits}
}.\nonumber
\end{align}
The solution $x^{(k)}$ to~\eqref{eq:UBPk} at iteration $k$ provides us
with a candidate speed limit $u^{(k)}:=\sum_{i\in \mathcal{O}}
  \gamma^{(i)} \reviseone{x^{(k)}_i} %
  $ \reviseone{which respect the sustainability constraints}. For each
$l\in\mathcal{L}$ with the given $u^{(k)}$, the admissible sample trajectory
${\rho}^{(l)}$ is uniquely determined by $(\omega^{(l)},
\rho^{(l)}(0),\supscr{r}{in,(\textit{l})},\supscr{r}{out,(\textit{l})}
)$ using the uniqueness solution of the linear time-invariant
systems. Therefore, the element $(y,{\rho}) \in\Psi(x^{(k)})$ is
unique. Using the constraints set $\Psi(x^{(k)})$, we then construct
the unique admissible sample trajectories $\{{\rho}^{(l,k)}\}_{l \in
  \mathcal{L}}$. The unique admissible sample trajectories allow us to
reduce~\eqref{eq:P4} to the following lower bounding problem\footnote{\reviseone{Given $u^{(k)}$ and $\{{\rho}^{(l,k)}\}_{l \in \mathcal{L}}$, the equivalent dual of~\eqref{eq:LBPk} is %
\begin{equation*}
  \begin{aligned}
    \min_{\rho^{(l)},l\in\mathcal{L}} &  \frac{1}{NT} \sum_{e,t,l}  u^{(k)}_e(t){\rho}^{(l)}_e(t) \\
    \st \quad &  0 \leq \rho^{(l)} \leq \rho^{c}(u^{(k)}), \; \forall l \in \mathcal{L}, \\
    & \sum_{l \in \mathcal{L}} \Norm{{\rho}^{(l)} - {\rho}^{(l,k)}} \leq  \epsilon(\beta),
  \end{aligned}
\end{equation*}
which results in more efficient online solutions.}}
{\leqnomode
\begin{align}
  \label{eq:LBPk} \tag{{LBP}$_k$}~\max\limits_{\substack{z,
      \lambda,\mu,\nu,\eta }} &-\lambda \epsilon(\beta) -
  \frac{1}{N}\sum_{e,t,l}{\left( \overline{f}_e
    \overline{\rho}_e\eta_e^{(l)}(t)  - \nu^{(l)}_e(t){\rho}^{(l,k)}_e(t) \right) }, \\
  \st \quad & (z,\lambda,\mu,\nu,\eta) \in \Phi(x^{(k)}). \nonumber
\end{align}}
Note that~\eqref{eq:LBPk} is a linear program and much easier to solve
than the non-convex Problem~\eqref{eq:P4}. Let $\obj_{k}$
denote the optimal objective value of~\eqref{eq:LBPk}. If
Problem~\eqref{eq:LBPk} is solved to optimum with a finite $\obj_{k}$,
we will then obtain a feasible solution of~\eqref{eq:P4} with speed
limit $u^{(k)}:=\sum_{i\in \mathcal{O}} \gamma^{(i)} \reviseone{x^{(k)}_i} %
$ and certificate
$\hat{J}(u^{(k)}):={\obj}_{k}$; otherwise, Problem~\eqref{eq:LBPk} is
either infeasible or unbounded, i.e., $\obj_{k}=-\infty$. The lower
bound of~\eqref{eq:P4} can be calculated by $\LB_k=
\max_{p=1,\ldots,k}\{ {\obj}_{p} \}$. The stopping criterion of the
algorithm can be determined by one of the following criteria
\begin{enumerate}
\item ${\UB}_k -{\LB}_k \leq
\epsilon$,
\item \eqref{eq:UBPk} is infeasible,
\item A satisfactory
sub-optimal solution is found after certain running time
$\subscr{T}{run}$.
\end{enumerate}
In~\cite{LX-TA-BP:11}, it is shown that such algorithms convergence to
a global $\epsilon$-optimal solution after finite number of iterations
when we use the first and second stopping criteria.  The third
solution criterion allow us to find a potentially good performance-guaranteed feasible
solution within certain running time $\subscr{T}{run}$.
\begin{remark}[Online tractable solutions to~\eqref{eq:P4}] {\rm
    \reviseone{The data-driven control requires to solve a sequence
      of~\eqref{eq:P4} online. We achieve this by an online ``warm
      start'' of Algorithm~\ref{Alg:issa}, which employs an
      assimilation set $\mathcal{I}_t:=\{u^{(s)}\}_s$ that contains
      the historically-generated speed-limit candidates, where $s$
      indexes these candidates. In particular, at each time solving a
      problem~\eqref{eq:P4}, the candidates in $\mathcal{I}_t$ can be
      explored by the solution to~\eqref{eq:LBPk} of each
      $u^{(s)}\in\mathcal{I}_t$. Notice that these~\eqref{eq:LBPk} can
      be executed in parallel. Then, these examined candidates
      contribute to integer cuts in~\eqref{eq:UBPk} when executing
      Algorithm~\ref{Alg:issa}. At the
      termination %
      of the current~\eqref{eq:P4}, a new set of candidates are
      updated to $\mathcal{I}_{t+1}$ for later evaluation of~\eqref{eq:P4}.  }}
\end{remark}
\section{Analysis via Second-order Cone Problems} \label{sec:tools}
\reviseone{This section provides a tool to analyze the efficacy of the
  proposed algorithm for the nonconvex Problem~\eqref{eq:P4}. In
  particular,} we propose a second-order cone relaxation for the
non-convex quadratic terms \reviseone{$\{\nu^{(l)} \circ
  {\rho}^{(l)}\}_{l \in\mathcal{L}}$} in~\eqref{eq:P4}, and a
second-order cone relaxation for it. We then present the conditions
under which this convex relaxation is exact.
\reviseone{To enable the tool for analysis}, we assume the following:

\begin{assumption}[Highway densities are  nontrivial]
  For all $e \in \edges$, $t \in \mathcal{T}$ and $l \in\mathcal{L}$,
  we assume ${\rho}^{(l)}_e(t) \geq \epsilon(\beta)$, where the
  parameter $\epsilon(\beta)$ is the radius of the Wasserstein ball.  \label{assump:3}
\end{assumption}
\begin{remark}[On nontrivial highway densities] \label{remark:relax}
  {\rm Assumption~\ref{assump:3} depends on the radius of the
    Wasserstein ball, which is selected as in
    Remark~\ref{remark:radius}. In reality, we could select the value
    $\epsilon(\beta)$ to be sufficiently small, even if it potentially
    sacrifices confidence on performance guarantees. In any case,
    there are three cases to consider: 1) the density on each segment
    of highway is just zero, 2) there are zero density values
    ${\rho}^{(l)}_e(t)$, for some $(e,t,l)$, while the rest of
    ${\rho}^{(l)}_e(t)$ are upper bounded by a value that is smaller
    than the maximal critical density $\max_{\reviseone{u_e(t)} \in
      \Gamma}\rho^{c}_e(\reviseone{u_e(t)})$, and 3) there are some values
    ${\rho}^{(l)}_e(t)$ that go beyond the maximal critical
    density, \reviseone{e.g., ${\rho}^{(l)}_e(t) > \epsilon(\beta) + \argmax_{u \in \Gamma} \rho^{c}_e(u)$}. In the first case, no congestion would happen and there
    is no need for speed-limit control. The second case can be handled by
    tuning $\epsilon(\beta)$ to be small enough. In the third case, \reviseone{with a given small $\epsilon(\beta)$},
    there is already congestion on some segment of the highway and no
    feasible speed limit would eliminate that congestion.  }\end{remark}
Assumption~\ref{assump:3} enables us to explore properties of the optimizers
of~\eqref{eq:P4} as in the following proposition.
\begin{proposition}[Optimizers in a cone]
  Let $\sol^{\star}:=(x^{\star},y^{\star},z^{\star}, {\rho}^{\star},
  \lambda^{\star},\mu^{\star},\nu^{\star},\eta^{\star})$ be any
  optimizer of Problem~\eqref{eq:P4}. If Assumption~\ref{assump:3}
  holds, we have $\nu^{\star} \circ {\rho}^{\star} \geq
  \vectorzeros{nTN}$.
\label{prop:cone}
\end{proposition}%
Proposition~\ref{prop:cone} allows us to explore structure of the
bi-linear terms $\{\nu^{(l)} \circ {\rho}^{(l)}\}_{l \in\mathcal{L}}$
via second-order cone constraints. This is achieved by introducing
variables $\vartheta^{(l)}_{e}(t)$ and writing Problem~\eqref{eq:P4}
as follows {\leqnomode
\begin{align} \small
  \label{eq:P4A}\tag{$\mathbf{P4'}$}~ \max\limits_{\substack{x,y,z,
      {\rho}, \\ \lambda,\mu,\nu,\eta,\vartheta }} -\lambda \epsilon(\beta)
  - \frac{1}{N}\sum_{e,t,l}{\overline{f}_e\overline{\rho}_e\eta_e^{(l)}(t)}
    + \frac{1}{N}\sum_{e,t,l}{ \left(\vartheta^{(l)}_{e}(t)\right)^2},
\end{align}}
{\reqnomode
\begin{align}
  \st \quad & \vartheta^{2} \leq \nu \circ {\rho}, \label{eq:cone1} \\
  & {\textbf{\small speed limits}}\; \eqref{eq:ueNew}, \;
  {\textbf{\small sample trajectories }}  \{ \eqref{eq:y}, \; \eqref{eq:RhohatNew2} \} , \nonumber\\
  & {\textbf{\small no congestion}} \; \{\eqref{eq:z},\;
  \eqref{eq:dual1}, \;\eqref{eq:dual2}, \;\eqref{eq:dual3},
  \;\eqref{eq:dual4}\}. \nonumber
\end{align}}
For each $e\in\edges$, $t\in\mathcal{T}$ and
$l\in\mathcal{L}$, the constraint~\eqref{eq:cone1} can be equivalently
written as the following second-order cone:
\begin{equation}
\sqrt{ \left(\nu^{(l)}_{e}(t)\right)^2 +\left(\rho^{(l)}_{e}(t)\right)^2 + 2\left(\vartheta^{(l)}_{e}(t)\right)^2} \leq \nu^{(l)}_{e}(t) +\rho^{(l)}_{e}(t).
\label{eq:cone2}
\end{equation}
Problem~\eqref{eq:P4} and~\eqref{eq:P4A} are equivalent as the following:
\begin{lemma}[Equivalent optimizer sets of~\eqref{eq:P4} and~\eqref{eq:P4A}]
 \label{lemma:coneP4A}
 Problem~\eqref{eq:P4A} is equivalent to~\eqref{eq:P4} in the sense
 that their optimal objective values are the same and the set of
 optimizer of~\eqref{eq:P4} is the projection of that
 of~\eqref{eq:P4A}.  Further, any feasible solution of~\eqref{eq:P4A}
 can give us a valid performance guarantee~\eqref{eq:perfnew} with
 certificate to be objective function of~\eqref{eq:P4A} evaluated at
 that feasible solution.
\end{lemma}
Problem~\eqref{eq:P4A} is still non-convex. To approximate the
quadratic terms in the objective, we will first show that variables
$\vartheta^{(l)}_{e}(t)$ are bounded.
\begin{lemma}[Bounded variable $\vartheta$]
  \reviseone{Let} Assumption~\ref{assump:2} on bounded $\eta$ hold, then there exists large enough scalar
  $\overline{\vartheta}$ such that $|\vartheta^{(l)}_e(t)| \leq
  \overline{\vartheta}$ for all $e \in \edges$, $t \in \mathcal{T}$ %
  and $l \in\mathcal{L}$. \label{lemma:vartheta}
\end{lemma}
Lemma~\ref{lemma:vartheta} enables us to approximate each component of
$\vartheta$ by a finite set of points within its range. Let (sufficiently large) $K$
denote the number of points and let us denote the set of
these points by $Q:= \{\pi_1, \ldots, \pi_K\} \subset \real$.
We use the set $\mathcal{Q}:=\untilinterval{1}{K}$ to index these
points. For each edge $e\in\edges$, time $t\in\mathcal{T}$ and sample
$l\in\mathcal{L}$, let us define the binary variable
$q^{(l)}_{e,k}(t)$ to be equal to one if $\vartheta^{(l)}_e(t)$ is
approximated by $\pi_k$; otherwise $q^{(l)}_{e,k}(t)=0$. Then for each
$e\in\edges$, $t\in\mathcal{T}$ and $l\in\mathcal{L}$, we will then
represent $\vartheta^{(l)}_e(t)$ by the following constraints
\begin{equation}
	\begin{aligned}
          & \vartheta^{(l)}_e(t)= \sum_{k\in \mathcal{Q}} \pi_k
          q^{(l)}_{e,k}(t), \quad \sum_{k \in
            \mathcal{Q}}q^{(l)}_{e,k}(t)=1, \\
          & \qquad\qquad\qquad\qquad \qquad \quad \forall e \in \edges, \; t\in\mathcal{T}, \; l\in\mathcal{L},\\
          & q^{(l)}_{e,k}(t) \in \{0,\; 1\}, \; \forall e \in \edges,
          \; t\in\mathcal{T}, \; l\in\mathcal{L}, \; k \in
          \mathcal{Q}.
	\end{aligned} 	\label{eq:varthetanew}
\end{equation}
Using this representation, we find approximated solutions
of~\eqref{eq:P4A} by solving the following {\leqnomode
\begin{align} \small
  \label{eq:P5}\tag{$\mathbf{P5}$}~ \max\limits_{\substack{x,q,y,z,
      {\rho}, \\ \lambda,\mu,\nu,\eta,\vartheta }} -\lambda
  \epsilon(\beta)
  - &\frac{1}{N}\sum_{e,t,l}{\overline{f}_e\overline{\rho}_e\eta_e^{(l)}(t)} \\
  & \qquad + \frac{1}{N}\sum_{e,t,l,k}{ \left(\pi_k\right)^2
    q^{(l)}_{e,k}(t)}, \nonumber
\end{align}}
{\reqnomode
\begin{align}
\st \quad & {\textbf{\small level approximation}}\; \eqref{eq:varthetanew}, \; {\textbf{\small second-order cone}} \; \eqref{eq:cone2}, \nonumber\\
	&	{\textbf{\small speed limits}}\; \eqref{eq:ueNew}, \; {\textbf{\small sample trajectories }}  \{ \eqref{eq:y}, \; \eqref{eq:RhohatNew2} \} , \nonumber\\
  &  {\textbf{\small no congestion}} \;  \{\eqref{eq:z},\; \eqref{eq:dual1}, \;\eqref{eq:dual2}, \;\eqref{eq:dual3}, \;\eqref{eq:dual4}\}. \nonumber
\end{align}}Problem~\eqref{eq:P5} is an SOCMIP and can be solved to
optimum by commercial solvers, such as GUROBI and MOSEK. Note that for
any feasible solution of~\eqref{eq:P5}, it is feasible
for~\eqref{eq:P4A} and its objective value is a valid lower bound for
that of~\eqref{eq:P4A}. Thus, any feasible solution of~\eqref{eq:P5}
together with its objective value provide performance
guarantees~\eqref{eq:perfnew} via Lemma~\ref{lemma:coneP4A}.
Further, as the number of partition points $K \rightarrow \infty$,
Problem~\eqref{eq:P5} is computationally equivalent to
Problem~\eqref{eq:P4A}, and therefore the same as
Problem~\eqref{eq:P4}. Notice that online solutions to
  Problem~\eqref{eq:P5} are inaccessible and~\eqref{eq:P5} serves as a
  tool for the performance analysis of the proposed algorithm.
In the following section, we leverage this tool to analyze the performance of the algorithm in an academic example. %
\section{Simulations}\label{sec:Sim}
In this section, we demonstrate in an example how to efficiently find
a solution to~\eqref{eq:P4} that results in a robust data-driven
variable-speed limit $u$ with performance
guarantee~\eqref{eq:perfnew}. The analysis of the results are
twofold. First, we verify the effectiveness of the proposed integer
solution search algorithm by comparing it with the monolith approach
that solves an approximation to~\eqref{eq:P4},
Problem~\eqref{eq:P5}. Second, to verify the robustness of the
solution and performance guarantees in probability, we compare the
resulting distributionally robust data-driven control with the control
obtained from the sample-averaged optimization problem. Both controls
are applied on a naive highway simulator
developed via the standard
cell transmission model~\cite{CFD:94}.

\noindent \textbf{Simulation setting:} We consider an \reviseone{$8$-lane} highway with
length $L=10\rm{km}$ and we divide it into $n=5$ segments \reviseone{with equal length}. Let the
unit size of each time slot $\delta=30\sec$ and consider $T=20$ time
slots for a $10$min planning horizon.  For each edge $e \in \edges$,
we \reviseone{propose} a traffic jam density of
$\overline{\rho}_e=1050\rm{vec/km}$, a capacity of $\overline{f}_e=3.1\times
10^4\rm{vec/h}$\footnote{The unit ``vec'' stands
  for ``vehicles''. \reviseone{Notice the proposed capacity is about $50\%$ higher than the actually highway capacity in order to leverage the actual fundamental diagram for control.} } and a maximal speed limit of
$\overline{u}_e=140\rm{km/h}$. Let us consider $m=5$ different
candidate speed limits $\Gamma=\{40{\textrm{km/h}},\;
60{\textrm{km/h}},\;80{\textrm{km/h}},\;100{\textrm{km/h}},\;120{\textrm{km/h}}
\}$. On the $\supscr{4}{th}$ edge $e:=(3,4) \in \edges$, we assume an
accident happens during $\mathcal{T}$ with parameters
$\overline{f}_e=2.7\times 10^4\rm{vec/h}$. To evaluate the effect of
the proposed algorithm, samples of the random variables $\rho(0)$,
$w$, $\supscr{r}{in}$ and $\supscr{r}{o}$ are needed. In real-case
studies, samples $\{\rho^{(l)}(0)\}_{l\in\mathcal{L}}$ can be obtained
from highway sensors (loop detectors), while samples of the uncertain
mainstream flows $\{\omega^{(l)}\}_{l\in\mathcal{L}}$ and flow
fractions $\{\supscr{r}{in,(\textit{l})},\;
\supscr{r}{o,(\textit{l})}\}_{l\in\mathcal{L}}$ can be constructed
either from a database of flow data on the highway, or from the
current measurements of ramp flows with the assumption that the
stochastic process $\{\omega(t),\; \supscr{r}{in}(t),\;
\supscr{r}{o}(t)\}_{t\in\mathcal{T}}$ is trend stationary.

\noindent \textbf{Fictitious datasets:} In this simulation example,
the index set of accessible samples is given by
$\mathcal{L}=\{1,2,3\}$. For each $l\in\mathcal{L}$, let us assume
that each segment $e\in\edges$ initially operates under a free flow
condition with an initial density $\rho^{(l)}_e(0)=260\rm{vec/km}$. To
ensure significant inflows of the system, we let the samples
$\{\omega^{(l)}(t)\}_{l\in\mathcal{L},t \in \mathcal{T}}$ of the
mainstream inflow to be chosen from the uniform distribution within
interval $[2\times 10^4,2.4\times 10^4]\rm{vec/h}$. For each edge
$e\in\edges$ and time $t\in\mathcal{T}$, we further assume that
samples $\{\supscr{r}{in,(\textit{l})}(t)\}_{l\in\mathcal{L}}$ and
$\{\supscr{r}{o,(\textit{l})}(t)\}_{l\in\mathcal{L}}$ are generated
from uniform distributions within interval $[0,5\%]$ and $[0,3\%]$,
respectively. We also let the confidence value be $\beta=0.05$ and the
radius of the Wasserstein Ball $\epsilon(\beta)=0.985$ as calculated
in~\cite{DL-SM:18-cdc}.

\noindent \textbf{Effectiveness of the algorithm:} To generate
feasible solutions that can be carried out for a real time
transportation system, we allocate $\subscr{T}{run}=1$min execution
time for control design and run algorithms on a machine with two core
$1.8$GHz CPU and $8$G RAM. In this allocated $1$ minute, we consider
the speed limit design in $2$ approaches: 1) we run the proposed
Algorithm~\ref{Alg:issa} to solutions of the Problem~\eqref{eq:P4},
and 2) we run optimization solver MOSEK to solutions of the monolith
Problem~\eqref{eq:P5}. The partition number $K=5$ is selected for the
monolith approach.

We present in Table~\ref{tab:result} the comparison of the mentioned two approaches. In $1$ minute, the Algorithm~\ref{Alg:issa} executed $19$ candidate speed limits where $2$ feasible speed limits were found at time $6.7\sec$ and $28.7\sec$. We verified that
$\hat{J}(u^{(2)})=1.17\times 10^5\rm{vec/h}$ is the
highest certificate obtained, i.e., $\hat{J}(u^{(2)}) \in
\argmax_{p=1,2} \{ \hat{J}(u^{(p)}) \;| \; u^{(p)}
  \text{ is feasible} \}$, and the desired speed limits are
\begin{equation*}
  u^{(2)}=[120,100,80,80,100]{\rm{km/h}}.
\end{equation*}
Compared with the proposed algorithm, the monolith approach returned a
feasible solution with the speed limits
$\supscr{u}{mon}=[120,120,120,40,40]{\rm{km/h}}$ and an approximated
throughput $3.93\times 10^4\rm{vec/h}$.  It can be seen that 1) the
gap (difference between UB and LB) obtained from the
Algorithm~\ref{Alg:issa} is tighter than that obtained from the
monolith approach, and 2) the implementable speed limits proposed by
the Algorithm~\ref{Alg:issa} results in higher throughput than that
achieved by the monolith.

In the following subsection, we use the speed limits $u^{(2)}$ to
verify the guarantees on congestion elimination with high probability.

\begin{table} \centering \small
\begin{tabular}{@{}lll@{}}\toprule
    & Algorithm~\ref{Alg:issa}$^{a}$ & Monolith \\ \midrule
  $\#$ of feasible cand.$^{b}$ & 2 & 1 \\
  $\#$ of infeasible cand. & 17 & NA \\
  LB (vec/h)  & $1.17 \times 10^{5}$ & $3.93 \times 10^{4}$ \\
  UB (vec/h) & $1.55 \times 10^{5}$ & $1.55 \times 10^{5}$ \\
\bottomrule
\end{tabular}
 \vspace{0.5ex} \\
  \raggedright
 \footnotesize $^{a}$: Subproblems~\eqref{eq:UBPk} and~\eqref{eq:LBPk} are solved via MOSEK. \\
$^{b}$: Candidate speed limit.
 \caption{The efficiency of the proposed Algorithm~\ref{Alg:issa}.} \label{tab:result}
\end{table}

\noindent \textbf{Distributionally robust decisions:} To demonstrate
the distributional robustness of our approach, we compare the
performance of our speed limits design $u^{(2)}$ with the performance
of the speed limits developed from a sample average optimization
problem, which is also known as the dual version of the scenario-based approach, such as in~\cite{SL-AS-JF-AN-EC-HH-BDS:16}. In particular, the sample averaged version of~\eqref{eq:P}
(equivalently,~\eqref{eq:P1}) is the one substitutes the unknown
distribution $\prob(u)$ with its empirical distribution
$\hat{\prob}(u)$. The resulting tractable formulation of the sample
average problem, analogous to~\eqref{eq:P4}, is the following
{\leqnomode
\begin{align} \small
  \max\limits_{\substack{x,y,z, {\rho}, \\ \mu,\nu,\eta }} &
  - \frac{1}{N}\sum_{e,t,l}{\overline{f}_e\overline{\rho}_e\eta_e^{(l)}(t)}
    + \frac{1}{N}\sum_{e,t,l}{\nu^{(l)}_e(t){\rho}^{(l)}_e(t)}, \nonumber \\
  \st \quad & \sum_{i\in \mathcal{O}} \gamma^{(i)}
  (\overline{\rho}-\rho^{c}(\overline{u}))\otimes
  \vectorones{T} \circ z^{(l)}_i  -\mu^{(l)} \nonumber \\
  & \hspace{1.8cm} +\overline{f}\otimes \vectorones{T} \circ
  \eta^{(l)} \geq \vectorzeros{nT}, \; \forall \;
  l\in\mathcal{L},  \nonumber \\
  & \nu^{(l)} =\mu^{(l)}+ \frac{1}{T}\sum_{i\in \mathcal{O}} \gamma^{(i)}x_{i}
  ,\; \forall \;
  l\in\mathcal{L},  \nonumber \\
  & \vectorzeros{nT} \leq \eta^{(l)} \leq \overline{\eta},
  \; \forall \; l\in\mathcal{L}, \nonumber \\
  & { \textbf{\small speed limits}}\; \eqref{eq:ueNew},
  \;{\textbf{\small dual variable}} \; \eqref{eq:z}, \nonumber \\
  & {\textbf{\small sample trajectories}} \{
  \eqref{eq:y}, \; \eqref{eq:RhohatNew2} \}. \nonumber
\end{align}} Note that the difference between the previous sample
average problem and~\eqref{eq:P4} is that the former has a Wasserstein
Ball radius $\epsilon(\beta)=0$ and, thus, unlike~\eqref{eq:P4}, it
does not provide a performance guarantee on congestion. We solve the
above sample average problem to a suboptimal solution via the
Algorithm~\ref{Alg:issa} with the same setting as in solution to
$u^{(2)}$. The resulting speed limit design is the following
\begin{equation*}
  \supscr{u}{sav}=[60,60,80,60,100]{\rm{km/h}}.
\end{equation*}
To verify the performance of $u^{(2)}$ and $\supscr{u}{sav}$, we
generated $\subscr{N}{val}= 10^3$ validation samples of random
variables $\rho(0)$, $w$, $\supscr{r}{in}$ and $\supscr{r}{o}$ that
are from the distributions described in the Fictitious datasets
paragraph. The speed limit design $u^{(2)}$, $\supscr{u}{sav}$ as well
as the validation dataset are integrated into a highway simulator with
the highway parameter settings described in the Simulation setting
paragraph.

\begin{figure}[tbp]%
\centering
\includegraphics[width=0.35\textwidth]{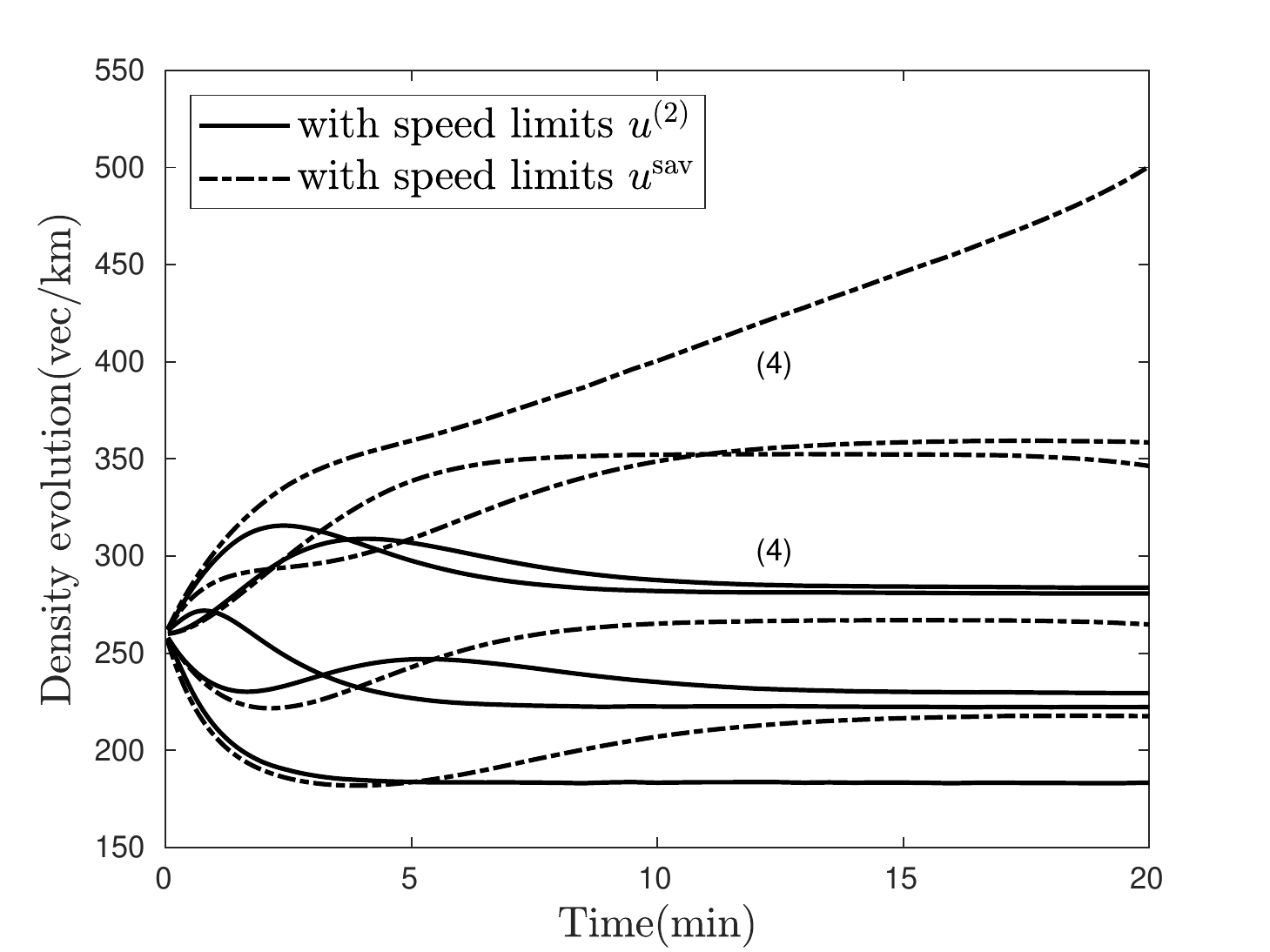}
\caption{\small Density evolution of each segment $e$, with speed
  limits $u^{(2)}$ and $\supscr{u}{sav}$. Each trajectory corresponds
  to a segment $e\in \{(1),(2),\ldots,(5) \}$. For simplicity we only
  label segment $(4)$, which happens to have an accident during the
  planning horizon.}%
\label{fig:result_density}%
\end{figure}

\begin{figure}[tbp]%
\centering
\includegraphics[width=0.35\textwidth]{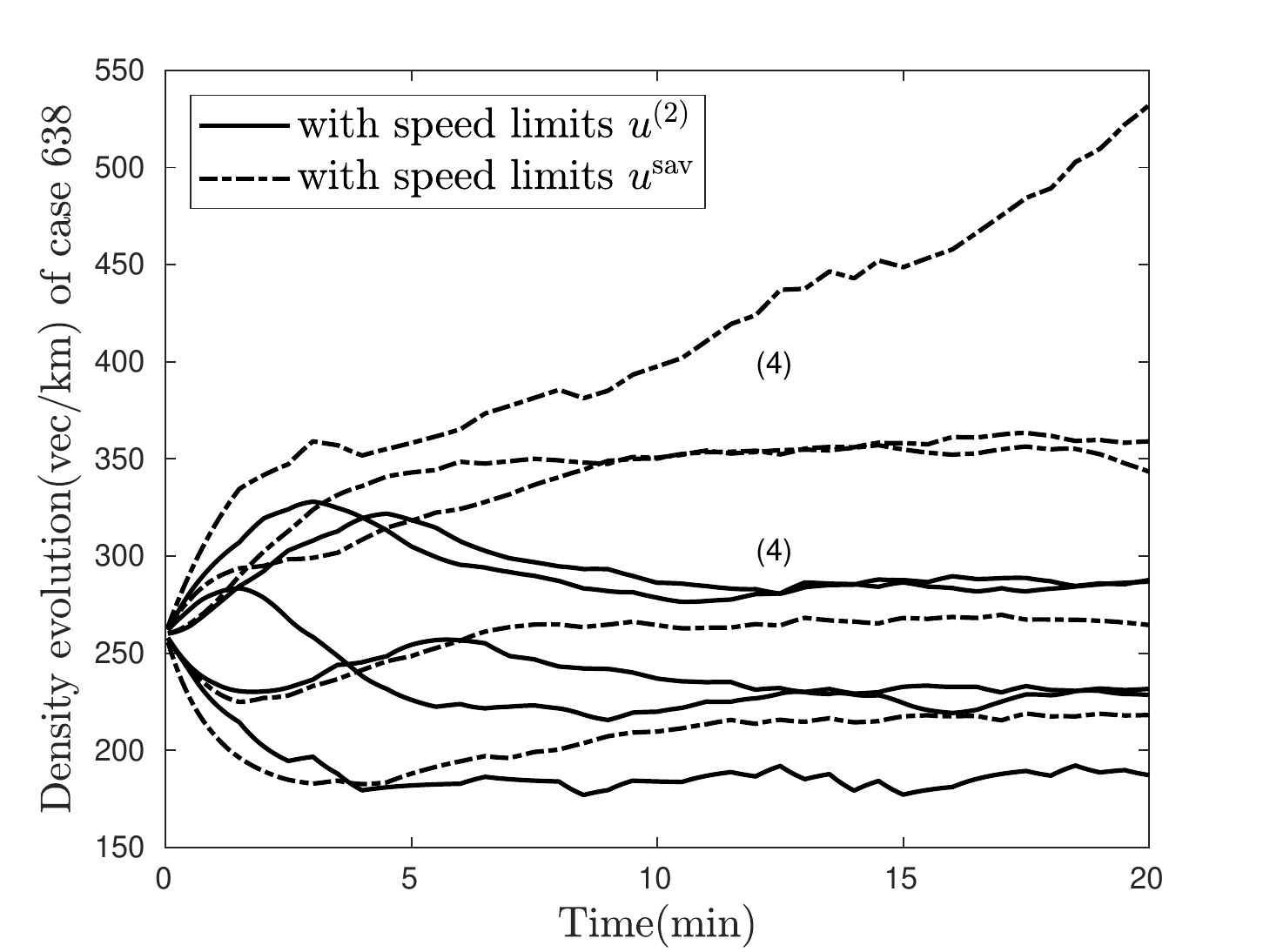}%
\caption{\small A representative density evolution of segments, with
 speed limits $u^{(2)}$ and $\supscr{u}{sav}$. The sample $638$ was arbitrarily chosen for demonstration purpose.}%
\label{fig:jamcase}%
\end{figure}

Due to space limitations we cannot showcase all admissible sample trajectories
for $10^3$ scenarios, therefore in Fig.~\ref{fig:result_density} we
show an average of the admissible sample trajectories, i.e., the function
$\frac{1}{\subscr{N}{val}}\sum_{l\in \untilone{\subscr{N}{val}}}
{\rho}_e^{(l)}(t)$ for each segment $e$, with speed limits $u^{(2)}$
and $\supscr{u}{sav}$, and present an arbitrarily chosen scenario
$638$ in Fig.~\ref{fig:jamcase}. We select the simulation time horizon
to be twice of that the planning horizon's in order to see the effect of
the design clearly.  We verified that the density evolution under
speed limits $u^{(2)}$, and, in particular, the density trajectory of
accident edge $(4)$, did not exceed the critical density values
($\rho_4^{c}(80\rm{km/h})=335\rm{vec/km}$) for each sample. Thus the
highway $\mathcal{G}$ is kept free of congestion in this planning
horizon $\mathcal{T}$ with high probability. However, the same robust
behavior can not be guaranteed under speed limits $\supscr{u}{sav}$,
as vehicles accumulate significantly on edge $(4)$ for too many
samples (contrast to its critical density
$\rho_4^{c}(60\rm{km/h})=403\rm{vec/km}$), see
Fig.~\ref{fig:result_density}.
We claim the robustness of our design compared to the design from
sample average problem, as the latter does not ensure such out-of-sample performance
guarantees.

\section{Simulation study}\label{sec:App}
\reviseone{In this section, we illustrate the proposed data-driven
  speed-limit control on a highway system located in San Diego,
  California, USA. The purpose of this simulation is to show the good
  behavior of the proposed method under uncertainty as compared with
  that of the current traffic speed limits on the highway.  While this
  is a first good indication, more work would be required to assess
  the behavior of the method over different traffic indices. However,
  our main focus is the theoretical development of the algorithm
  itself and so these questions are left for future work.
}

\begin{figure}[tbp]%
\centering
\includegraphics[width=0.45\textwidth]{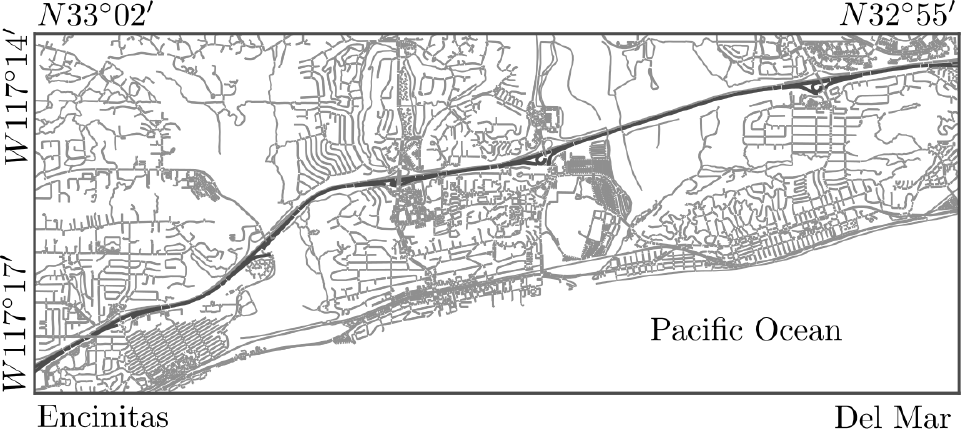}%
\caption{\small \reviseone{Highway section from Encinitas to Del Mar,
    San Diego, US.}}%
\label{fig:system_map}%
\end{figure}

\reviseone{\noindent \textbf{Highway system:} We selected a highway
  section from Encinitas to Del Mar on the I-5 San
  Diego Freeway with length $L \approx 11$ km, as shown in
  Fig.~\ref{fig:system_map}. The highway was divided into $n=26$
  segments with various lengths $\{ \len_e \}_{e\in\edges}$
  ranging from $200$ m to $2$
  km. %
  These segments have a number of lanes $\{ \lane_e \}_{e\in\edges}$
  ranging from $4$ to $8$,
  and %
  there are $7$ on-ramps and $5$ off-ramps distributed on the
  highway. We obtained real-time traffic data with $30$ seconds
  precision from the California Highway Performance Measurement System
  (PeMS), and used it to reproduce the actual traffic flow features
  via the software Simulation of Urban MObility
  (SUMO)~\cite{PAL-MB-LB-JE-YF-RH-LL-JR-PW-EW:18}, which is a
  microscopic and continuous traffic simulation
  package. %
  We calibrated the simulator using PeMS data which were collected
  between $12$pm and $2$pm on a particular day, and the speed limit on
  San Diego Freeway was $105$ km/h or $65$ mph.  }

\reviseone{\noindent \textbf{Data-driven control:} The data-driven
  control to solve~\eqref{eq:P} considered $T=80$ time slots dividing
  $4$-minutes planning horizon, with a unit size of each time slot
  $\delta=3$ seconds. The infrastructure-related parameters
  $(n,L,\len,\lane)$ were selected to be the same as that of the
  highway system, only that ramp-related ultra-short lanes were
  excluded.  We considered $m=6$ different candidate speed limits
  $\Gamma=\{ 30,\; 50,\; 70,\; 90,\; 110,\; 130 \}$
  $(\textrm{km/h})$. The control $u$ was obtained by implementing
  Algorithm~\ref{Alg:issa} that solves the equivalent~\eqref{eq:P4},
  taking a $\subscr{T}{run}=1$ minute and executing
  Algorithm~\ref{Alg:issa} every $2$ minutes. To achieve realistic
  driving instructions, we added extra speed limit constraints to
  ensure constant $u$ over each $1$ minute interval.}

\reviseone{\noindent \textbf{System monitor:} We assume the existence
  of a system monitor which provides information of the real-time
  traffic flow data $\{\varpi^{(l)}\}_{l\in\mathcal{L}}$ as
  well as the random events parameters $(\overline{f},
  \overline{\rho},\overline{u})$. Practically, these parameters can be
  calibrated in advance from PeMS historical data. In particular, at
  each time when Algorithm~\ref{Alg:issa} is to be executed, we
  consider $N=2$ accessible samples with the set $\mathcal{L}=\{
  1,2\}$. Precisely, values of $\varpi^{(1)}$ were constructed and
  propagated using real-time highway sensor measurements (loop
  detectors data in PeMS) and values of $\varpi^{(2)}$ were obtained
  as the $7$-day average data corresponding to the same time period
  ($12$pm to $2$pm). This results in a radius for the Wasserstein Ball
  $\epsilon(\beta) \approx 5$, given the confidence value
  $\beta=0.05$.  Notice that more samples can be added to reduce the
  radius if various source of measurements are accessible, e.g., such
  as real-time GPS data. On the other hand, the infrastructure and
  event-related data $\{(\overline{f}_e,
  \overline{\rho}_e,\overline{u}_e)\}_{e\in\edges}$ are determined
  theoretically, where values $\{\overline{f}_e\}_{e\in\edges}$ have
  range $[1.1,2.3] \times 10^4$ (vec/h),
  $\{\overline{\rho}_e\}_{e\in\edges}$ with values in $[0.5,1] \times
  10^3$ (vec/km/edge), and $\{\overline{u}_e\}_{e\in\edges}$ are
  assumed $200$ (km/h).  }

\begin{figure}[tp]%
\centering
\includegraphics[width=0.3\textwidth]{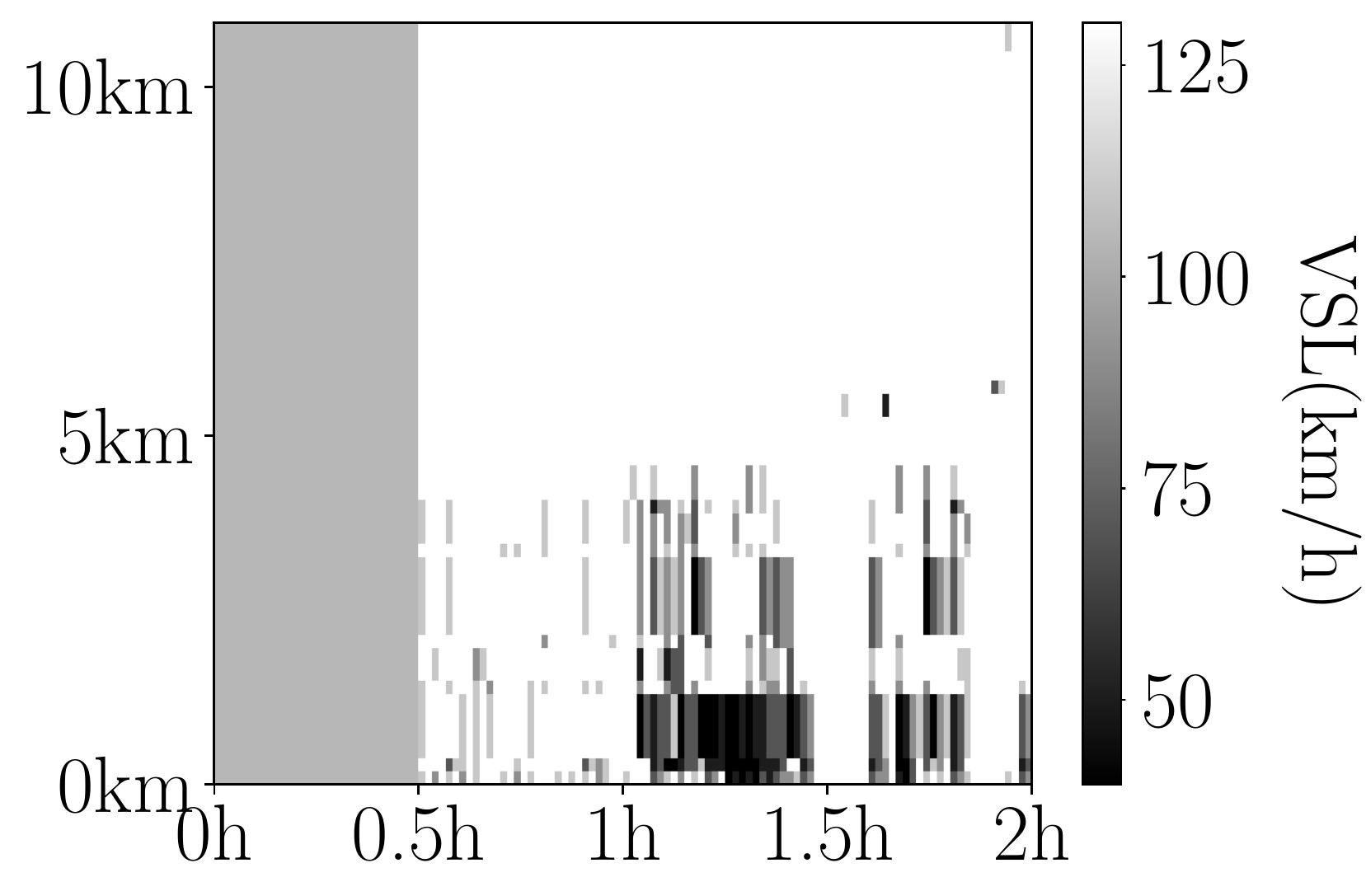}
\caption{\small Time-space profile of speed limits $u(t)$.}%
\label{fig:sd_speed_limit}%
\end{figure}

\reviseone{\noindent \textbf{Benchmarks:} We consider a $2$-hour
  scenario, from $12$pm to $2$pm, and assume a temporary lane closure
  on the $\supscr{15}{th}$ segment between $1$pm to $1:20$pm, which
  introduces a capacity and jam density drop by $35\%$ on that
  particular segment, located at $4.5$ km from the entry. Further, we
  implement the proposed data-driven control between $12:30$pm to
  $2$pm, and compare the resulting performance with that of the
  highway system without control, i.e., with constant speed limit
  $105$ km/h. Notice that, if congestion is inevitable, namely, the
  data-driven control problem~\eqref{eq:P} is infeasible under the
  admissible regime, we implement the default speed limit ($105$ km/h)
  instead.  }
  \begin{figure}[tp]%
  \centering
  \includegraphics[width=0.49\textwidth]{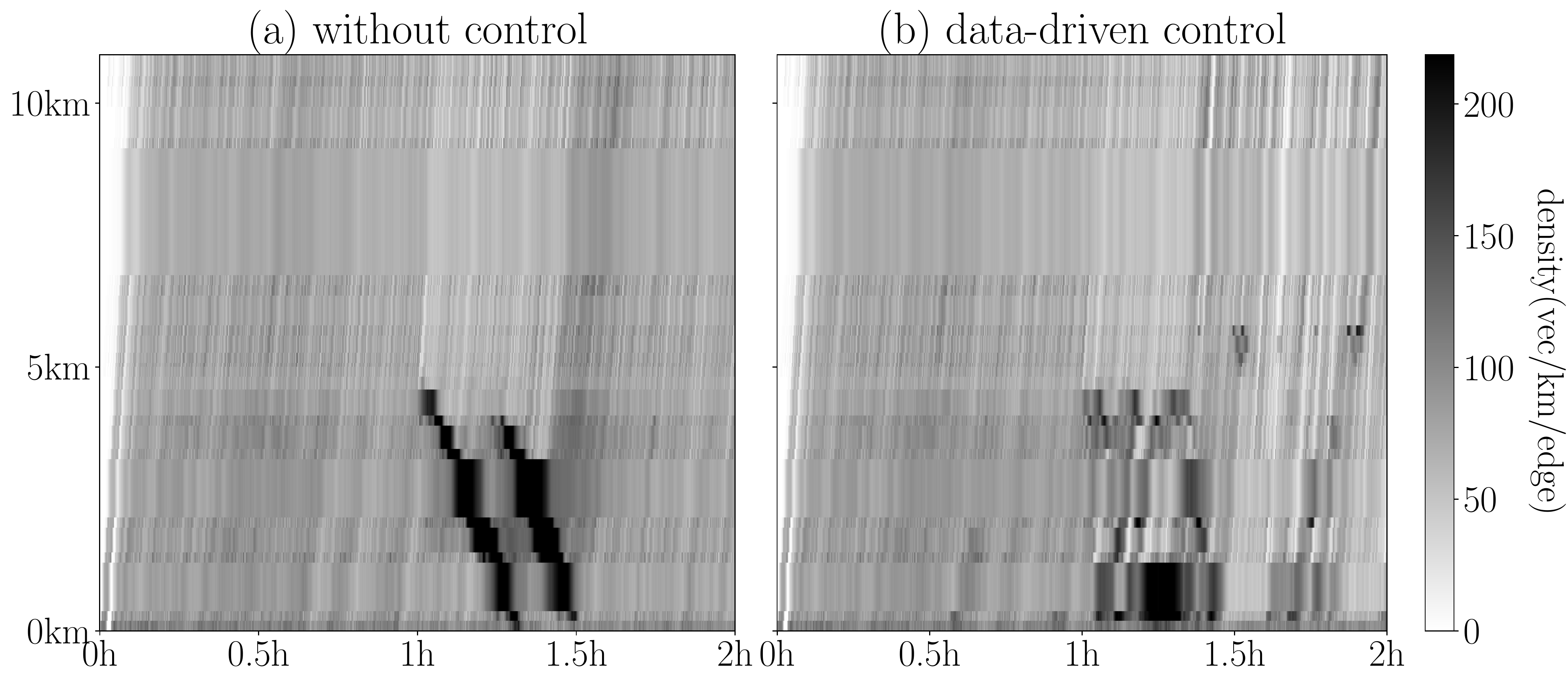}
  \caption{\small Time-space profile of traffic average density.}%
  \label{fig:sd_density}%
  \end{figure}
  \begin{figure}[tp]%
  \centering
  \includegraphics[width=0.49\textwidth]{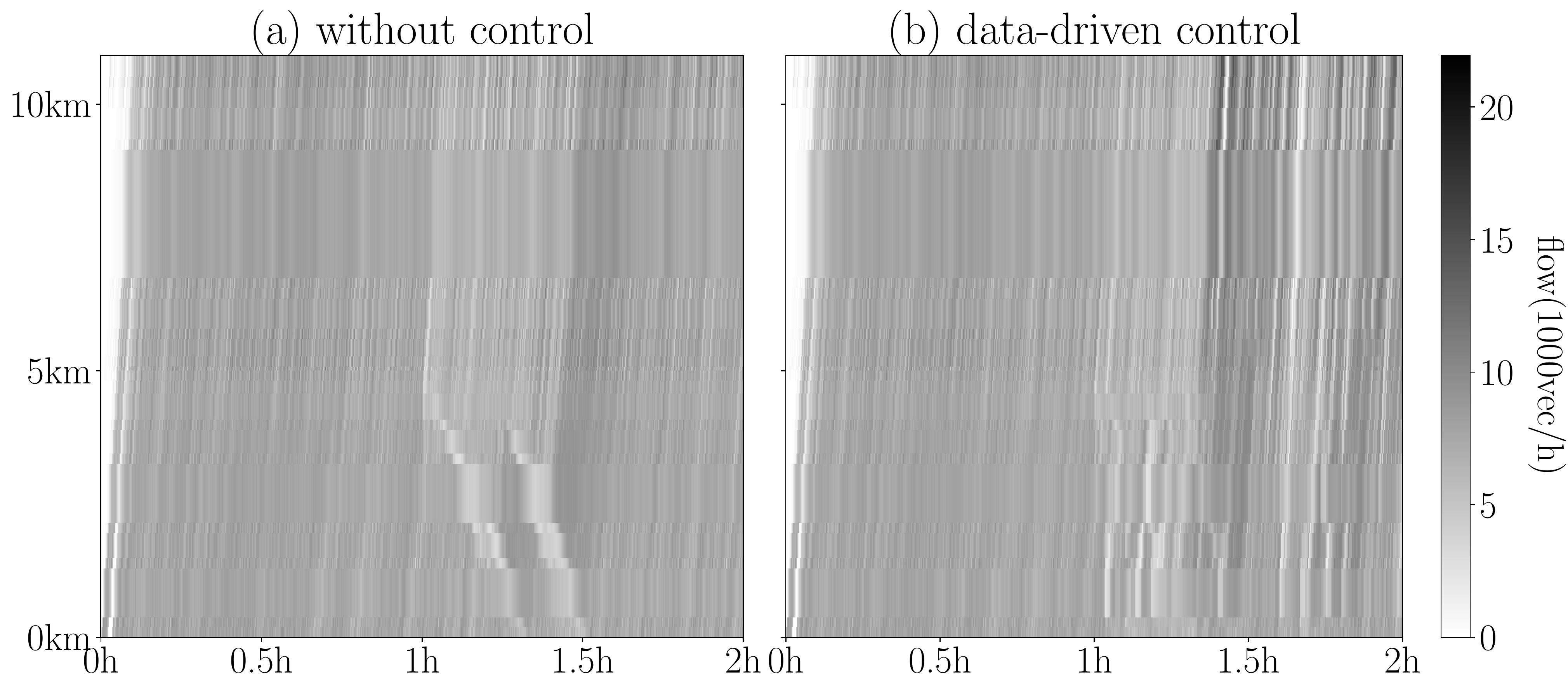}
  \caption{\small Time-space profile of traffic average flow.}%
  \label{fig:sd_flow}%
  \end{figure}
  \begin{figure}[tp]%
  \centering
  \includegraphics[width=0.49\textwidth]{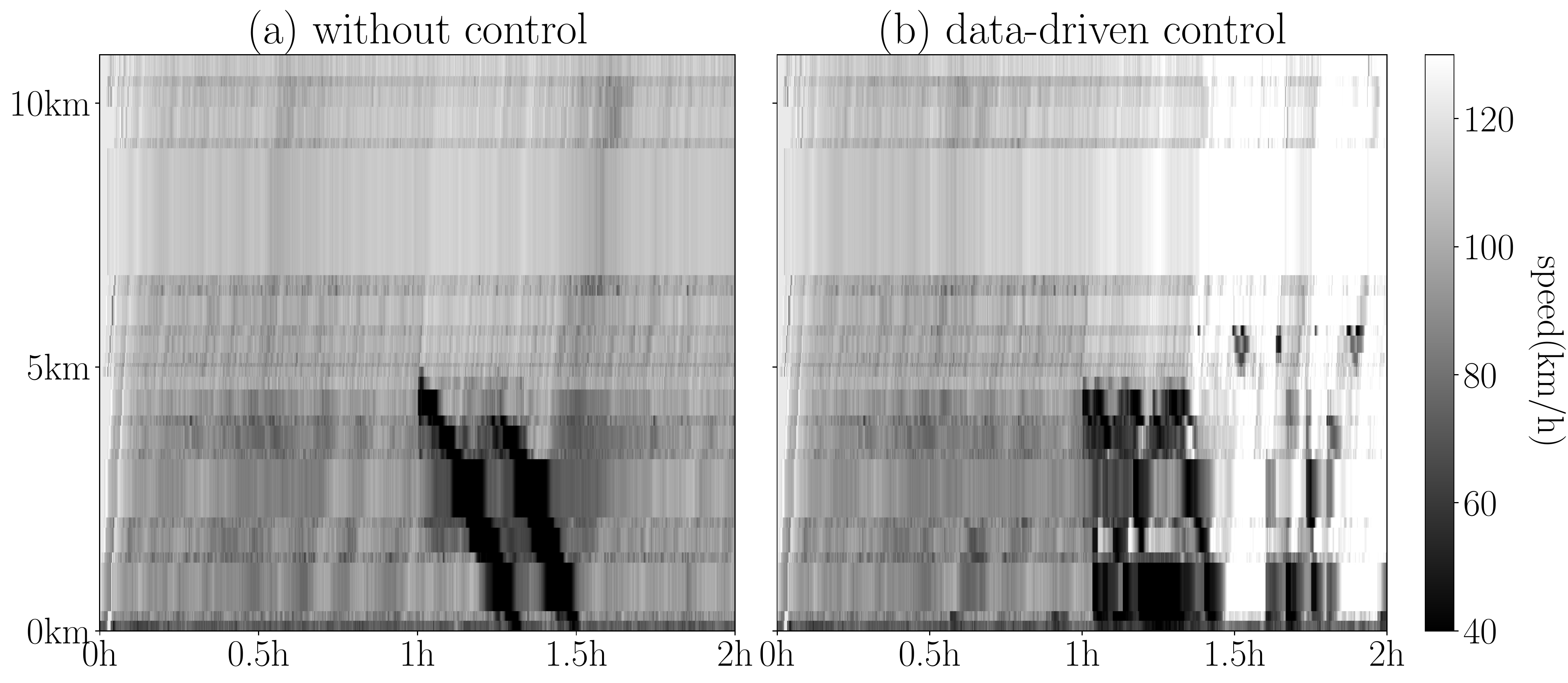}
  \caption{\small Time-space profile of traffic average speed.}%
  \label{fig:sd_speed}%
  \end{figure}
  \begin{figure}[tp]%
  \centering
  \includegraphics[width=0.49\textwidth]{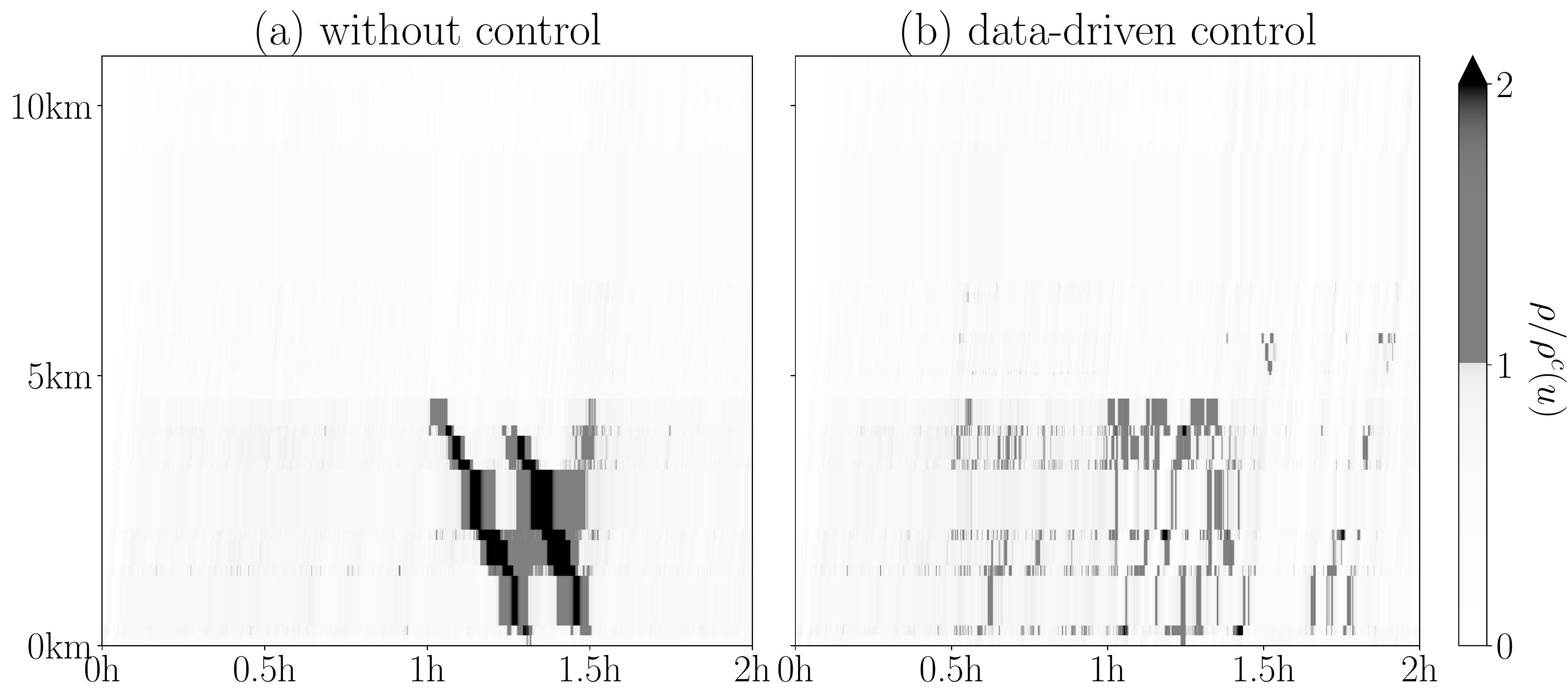}
  \caption{\small Time-space profile of congestion ratio $\rho(t)/\rho^{c}(u(t))$.}%
  \label{fig:sd_congestion_ratio}%
  \end{figure}
\reviseone{\noindent \textbf{Performance analysis:}
  Fig.~\ref{fig:sd_density} shows the traffic density profile in time
  and space, where the origin indicates the entry of the highway at
  the initial time, i.e., Encinitas at time $12$pm. The x-axis
  indicates the time (number of hours) passed from $12$pm and the
  y-axis is the distance to the highway entry. Similarly,
  Fig.~\ref{fig:sd_flow} and Fig.~\ref{fig:sd_speed} demonstrate the
  flow and average speed profiles of those, respectively. For
  comparison, subfigures (a)
  are profiles without control and subfigures (b) are those with the
  data-driven control. It can be observed that, during the period $0$h
  to $0.5$h ($12$pm to $12:30$pm), the two profiles were statistically
  identical. During period $0.5$h to $1$h ($12:30$pm to $1$pm), as
  the mainstream flows were moderate, the speed-limit control
  assigned a higher speed limit ($130$ km/h) than the default $105$
  km/h on the majority of the segments.
  See %
  Fig.~\ref{fig:sd_speed_limit} for the speed limit
  profile. After $1$h ($1$pm), a significant capacity drop occurs on
  the middle of the highway due to a temporary lane closure. This
  event leads to congestion on the preceding segment and backward
  congestion waves start to transmit on the highway, see, e.g.,
  subfigures (a). The darker parts on the profiles indicate congestion
  and notice how the congestion was transmitted to the entry of the
  highway over time. In addition, during the lane closure period
  ($1$pm to $1:20$pm), the data-driven control then took effect to
  cancel/reducing the congestion transmission by dynamically assigning
  low speed limits on the upper stream of the highway. See, e.g.,
  Fig.~\ref{fig:sd_speed_limit} on the speed limit assignment
  during $1$h to $1.4$h, which regulates the highway average speed as
  in Fig.~\ref{fig:sd_speed}(b). These actions eliminated the
  congestion waves and reduced the effect of the random events.
  See, e.g., Fig.~\ref{fig:sd_density} and~\ref{fig:sd_flow} for comparison of the effect of the congestion elimination,
  Fig.~\ref{fig:sd_congestion_ratio} for the significant reduction of the congestion ratio $\rho(t)/\rho^{c}(t)$, and Fig.~\ref{fig:sd_speed_limit} for the assignment of the variable
  speed limit.
  When the random event ended, the speed-limit control then resumed to
  normal operation, which, during $1.4$h to $2$h, dynamically assigned
  speed limits to account for uncertainties on random ramp flows.
  See, e.g., the scattered speed restrictions in
  Fig.~\ref{fig:sd_speed_limit}. Furthermore, notice in
  Fig.~\ref{fig:sd_congestion_ratio}(b) how the actual highway density
  $\rho(t)$ violated the prediction-and-assignment critical density
  $\rho^{c}(u)$ via speed limits $u$. These are the original driving
  forces to update speed limits. Notice that, during this whole
  scenario, the data-driven control problem~\eqref{eq:P} is
  feasible. Otherwise, the speed limit would be set to the default
  $105$ km/h at some time later than $0.5$h in
  Fig.~\ref{fig:sd_speed_limit}. This indicates a successful
  containment of the congestion. When the congestion is too
  heavy,~\eqref{eq:P} could be infeasible for some time period due to
  the extreme-high density on some of the segments.  In those
  scenarios, to handle the congestion, small enough candidates can be
  added to the candidate speed limit set $\Gamma$ in order for larger,
  admissible operation zone of~\eqref{eq:P}. Otherwise, the congestion
  is inevitable as~\eqref{eq:P} is infeasible, and we simply select
  the default, pre-selected speed limits.  At last, notice that the
  control performance relies heavily on the selection of the objective
  function of~\eqref{eq:P} as well as on the available information on
  the flow and random events data. We leave the questions regarding
  other traffic performance metrics and the improvement of the
  controller employing more accurate traffic models for the future
  work.  }
\section{Conclusions}\label{sec:Conclude}
In this paper, we formulate a data-driven predictive control problem with
probabilistic performance guarantees as a distributionally optimization problem. We
equivalently reformulate this intractable Problem~\eqref{eq:P0} %
into a
non-convex but tractable Problem~\eqref{eq:P4}, or Mixed Integer Second
Order Cone Problem~\eqref{eq:P5}. This is achieved by: 1) extending
Distributionally Robust Optimization theory to account for system
dynamics; 2) reformulation and relaxation techniques. Finally, we
adapt the idea of decomposition and propose an integer-solution search
algorithm for efficient solutions of Problem~\eqref{eq:P4},
or commercial solvers for that of
Problem~\eqref{eq:P5}. The proposed approaches provide performance guarantee~\eqref{eq:perf} for Problem~\eqref{eq:P}.
To explicitly demonstrate the proposed approach, we consider a highway speed-limit control problem that accounts for random inflows, outflows and events.
Finally, we demonstrate the theoretical effectiveness of this work via simulations. In particular, we simulated traffic flows on a highway in San Diego, and show the effectiveness of the designed speed limits numerically.

\appendix \label{appx:whole}

\begin{proof}{Lemma~\ref{lemma:RhoGener}}
Continuous functions of i.i.d. random variables generate
i.i.d. random variables. Thus, the admissible sample trajectories
$\{{\rho}^{(l)}\}_{l\in\mathcal{L}}$ generated
by~\eqref{eq:Rhohat} are i.i.d.~realizations~of~$\prob(u)$.
\qed \end{proof}

\begin{proof}{Lemma~\ref{lemma:hatRho}}
  The proof consists of two steps. First, we explore the boundedness
  property of $\rho_e(t)$ w.r.t. $\varpi$, for all $e \in \edges$ and
  $t\in \mathcal{T}$. Second, we claim that there exists an $a >1$
  such that $\mathbb{E}_{\prob(u)}[\exp(\Norm{{\rho}}^a)] < \infty$.

  \noindent \textbf{Step 1: (Boundedness of $\rho$)} Consider for each
  speed limit $u$ and time $t \in
  \mathcal{T}\setminus \{ 0\}$. Let $A(u,t)$ denote the matrix that is
  consistent with the highway system $\mathcal{G}$,  let $B(t)$ denote
  the column vector that encodes the mainstream flow $\omega$, and
  write the density $\rho(t)$ in the following compact form
\begin{equation*}
  \rho(t)= \Phi(t,0)\rho(0) + \sum\limits_{\tau=0}^{t-1} \Phi(t,
  \tau+1) B(\tau), \\
  \label{eq:Rhohat2}
\end{equation*}
where
\begin{equation*} \small
\begin{aligned}
&\Phi(t,\tau):=\begin{cases}
          I_n , & \; {\textrm{ if }} t=\tau, \\
          A(u,t-1)A(u,t-2)\cdots A(u,\tau), & \; {\textrm{ if } t > \tau,}
	\end{cases}
	\end{aligned}
\end{equation*}
\[ \tiny
A(u,t)= \begin{bmatrix}
    \reviseone{m_1(t)} &       &      &       \\
    n_1(t)  & \ddots &          & \\
        &   \ddots & \ddots &  &  \\
          &        &  &  &  \\
          &    &    n_{n-1}(t) & \reviseone{m_n(t)}
\end{bmatrix},  B(t)= \begin{bmatrix}
\reviseone{h_{(0,1)}} \omega(t) \\ 0 \\ \vdots \\ 0
\end{bmatrix},
\]
with
\reviseone{
\begin{equation*} \small
  \begin{aligned}
  & m_i(t)=1-\reviseone{h_e} u_{e}(t), \quad \forall  i \in \untilinterval{1}{n}, e=(i-1,i), \\
  &  n_j(t)= \min \left\{\frac{1-\supscr{r}{o}_{s}(t)}{1-\supscr{r}{in}_{e}(t)}
  \reviseone{h_e} u_{s}(t), \frac{\overline{f}_{e}}{\rho_{s}(t)}, \tau_e\overline{u}_e \frac{ \overline{\rho}_e - \rho_{e}(t)} {\rho_{s}(t)} \right\}, \\
  & \hspace{1.5cm} \forall j \in \untilinterval{1}{n-1}, s=(i-1,i), e=(i,i+1).
  \end{aligned}
\end{equation*}
}

Given that each component of $A(u,t)$ is bounded for fixed $u$ and
$t$, the induced norm of $A(u,t)$ is also bounded and we denote their
universal bound by $A_0$. \reviseone{Similarly, we denote by $\overline{h}$ the upper bound on $h_e$, $e\in \edges$.}
Then for every $t$ we can bound the infinity
norm of $\rho(t)$ as follows
\begin{equation*}
  \begin{aligned}
    \Norm{\rho(t)}_{\infty} &\leq \Norm{ \Phi(t,0) }_{\infty}
    \Norm{\rho(0)}_{\infty} \\
    & \hspace{2.5cm} +  \reviseone{\overline{h}}\sum\limits_{\tau=0}^{t-1} \Norm{\Phi(t, \tau+1)}_{\infty} \Norm{\omega(\tau)}_{\infty}, \\
    &\leq A_0^{t} \Norm{\rho(0)}_{\infty}
    +   \reviseone{\overline{h}} \sum\limits_{\tau=0}^{t-1} A_0^{t-\tau-1} \Norm{\omega(\tau)}_{\infty}, \\
    & \leq M_1(t)\left(\Norm{\rho(0)}_{\infty}
      +\sum\limits_{\tau=0}^{t-1}
      \Norm{\omega(\tau)}_{\infty} \right), \\
    &\leq (t+1)M_1(t)\Norm{\varpi}_{\infty}, \leq
    (t+1)M_1(t)\Norm{\varpi},
  \end{aligned}
\end{equation*}
where $M_1(t):=\max\{A_0^{t},\;  \reviseone{\overline{h}},\;  \reviseone{\overline{h}} A_0^{t} \}$.

\noindent \textbf{Step 2: (Light-tailed distribution)}
\reviseone{Given an $u$}, consider a time $t^{\star} \in \argmax_{t \in \mathcal{T}\setminus
  \{ 0\} } \{ \Norm{\rho(t)}_{\infty} \}$. Then by
norm equivalence we have
  \begin{equation*}
          \begin{aligned}
            \Norm{{\rho}} &\leq n \Norm{{\rho}}_{\infty}
            = n \max_{t \in \mathcal{T}\setminus \{ 0\}} \{ \Norm{\rho(t)}_{\infty} \}, \\
            &= n \Norm{\rho(t^{\star})}_{\infty}, \leq n
            (t^{\star}+1)M_1(t^{\star})\Norm{\varpi}.
          \end{aligned}
\end{equation*}
Then for any $a>1$ with
$\mathbb{E}_{\prob_{\varpi}}[\exp(\Norm{\varpi}^a)] < \infty$, let
$M_2=[{n(t^{\star}+1)M_1(t^{\star})}]^a < \infty$ and we have
\begin{equation*}
        \begin{aligned}
          \mathbb{E}_{\prob(u)}
          [\exp(\Norm{{\rho}}^a)] & \leq \mathbb{E}_{\prob_{\varpi}}
          [\exp( M_2 \Norm{\varpi}^a)], \\
          &= \exp( M_2)^a \cdot \mathbb{E}_{\prob_{\varpi}}
          [\exp(\Norm{\varpi}^a)] < \infty.
        \end{aligned}
\end{equation*}
That is, $\prob(u)$ is light tailed.
\qed \end{proof}

\begin{proof}{Theorem~\ref{thm:cert}}
  We prove \reviseone{the result} in two steps. First, we show the
  proposed set $\mathcal{P}(u)$ contains the unknown distribution
  $\prob(u)$ with high probability. Then, we show that the proposed
  certificate ${J}({u})$ provides \reviseone{performance} guarantees
  as in~\eqref{eq:perfnew}.

  \noindent \textbf{Step 1: (Tractable set containing $\prob(u)$)}
  Under Assumption~\ref{assump:1} and using $N$ i.i.d.~samples of
  ${\varpi}$, we obtain i.i.d.~samples
  $\{{\rho}^{(l)}\}_{l\in\mathcal{L}}$ of $\prob(u)$ via
  Lemma~\ref{lemma:RhoGener}. Then, with the distribution
\begin{equation*}
  \hat{\prob}(u):=
  \frac{1}{N}\sum_{l\in\mathcal{L}} \delta_{\{{\rho}^{(l)}\}},
\end{equation*}
we quantify the relation between $\hat{\prob}(u)$ and $\prob(u)$ via
Lemma~\ref{lemma:hatRho} and the following theorem:
\begin{theorem}[Measure concentration~{\cite[Theorem~2]{NF-AG:15}}]
  If $\prob(u) \in
  \subscr{\mathcal{M}}{lt}(\mathcal{Z}{(u)})$, then
\begin{equation*} \small
  \Prob^N\{ d_W(\prob(u),\hat{\prob}(u)) > \epsilon \}
  < \left\{ {\begin{array}{*{20}{l}}
        c_1 e^{-c_2N \epsilon^{\max \{2,\ell\}}},  & \textrm{if} \;
        \epsilon \leq 1,    \\
        c_1 e^{-c_2N \epsilon^a} ,  &  \textrm{if} \; \epsilon > 1,
\end{array}} \right.
\label{eq:radB}
\end{equation*}
for all $N \geq 1$, $\ell \ne 2$, and $\epsilon >0$, where the
parameter $\ell$ is the dimension of $\rho$, and parameters $c_1$,
$c_2$ are positive constants that only depend on $\ell$, $a$ and $b$
as in Assumption~\ref{assump:1}. \hfill $\square$
\label{thm:MeasConc}
\end{theorem}
Let us select $\epsilon:=\epsilon(\beta)$ to be the following
\begin{equation*}
  \epsilon(\beta) := \left\{ {\begin{array}{*{20}{l}}
        \left( \frac{\log(c_1 \beta^{-1})}{c_2 N}
        \right)^{1/\max\{2, \ell \}},
        \; & \textrm{if} \; N \geq \frac{\log(c_1 \beta^{-1})}{c_2}, \\
        \left( \frac{\log(c_1 \beta^{-1})}{c_2 N} \right)^{1/a},
        \; & \textrm{if} \; N < \frac{\log(c_1 \beta^{-1})}{c_2},
\end{array}} \right.
\end{equation*}
then Theorem~\ref{thm:MeasConc} leads to
\begin{equation*}
  {\Prob^N}\left( d_W(\prob(u),\hat{\prob}(u)) > \epsilon(\beta)  \right) < \beta,
\end{equation*}
or, equivalently,
\begin{equation*}
  {\Prob^N}\left( d_W(\prob(u),\hat{\prob}(u)) \leq \epsilon(\beta)  \right) \geq 1- \beta.
\end{equation*}
From the definition of the Wasserstein ball
$\mathbb{B}_{\epsilon(\beta)}(\hat{\prob}({u}))$ and the fact that
$\prob(u) \in \subscr{\mathcal{M}}{lt}(\mathcal{Z}{(u)})$, we have
\begin{equation*}
\begin{aligned}
  {\Prob^N} & \left( \prob(u) \in  \mathbb{B}_{\epsilon(\beta)}(\hat{\prob}({u}))  \right) \geq 1- \beta, \\
  {\Prob^N} & \left( \prob(u) \in
    \subscr{\mathcal{M}}{lt}(\mathcal{Z}{(u)}) \right) =1,
\end{aligned}
\end{equation*}
which results in
\begin{equation*}
  {\Prob^N}\left( \prob(u) \in \mathcal{P}(u) \right)\geq 1- \beta,
\end{equation*}
with
\begin{equation*}
  \mathcal{P}({u}) := \mathbb{B}_{\epsilon(\beta)}(\hat{\prob}({u})) \cap  \subscr{\mathcal{M}}{lt}(\mathcal{Z}{(u)}).
\end{equation*}
\noindent \textbf{Step 2: (Certificate of~\eqref{eq:P1})} From the
previous reasoning, for a given $u$ we have $\prob(u) \in
\mathcal{P}(u)$ with probability at least $1-\beta$. Then, with
probability at least $1-\beta$ the objective value of~\eqref{eq:P1}
satisfies the following
\begin{equation*}
  \mathbb{E}_{\prob(u)}\left[ H(u;{\rho}) \right]
  \geq \inf_{\mathbb{Q} \in \mathcal{P}({u})}
  {\mathbb{E}_{\mathbb{Q}} \left[H(u;{\rho})\right]}.
\end{equation*}
Let ${J}({u})$ be
\begin{equation*}
  {J}({u}):= \inf_{\mathbb{Q} \in \mathcal{P}({u})}
  {\mathbb{E}_{\mathbb{Q}} \left[H(u;{\rho})\right]},
\end{equation*}
which results in
\begin{equation*}
  {\Prob^N}\left(\mathbb{E}_{\prob(u)}
    {\left[ H(u;{\rho}) \right] }\geq
    {J}({u}) \right)\geq 1- \beta.
\end{equation*}
Then the proposed certificate ${J}({u})$ of~\eqref{eq:P1} satisfies
guarantee~\eqref{eq:perfnew}, which completes the proof. %
\qed \end{proof}

\begin{proof}{Theorem~\ref{thm:P2}}
  We achieve it by three steps. First, we show that
  Problem~\eqref{eq:P2} can be equivalently reduced to a
  finite-dimensional problem. Then we equivalently reformulate the
  resulting problem into a maximization problem. Finally, we show the
  performance guarantees of~\eqref{eq:P3}.

\noindent \textbf{Step 1: (Finite-dimensional reduction of~\eqref{eq:P2})}
  We express $J(u)$ as follows
 \begin{equation*} \small
   \begin{aligned}
     J(u) &= \begin{cases}
       \inf\limits_{\mathbb{Q}} & \int_{\mathcal{Z}(u)} H(u;{\rho}) \mathbb{Q}(d \rho), \\
       \st &  \mathbb{Q} \in \subscr{\mathcal{M}}{lt}(\mathcal{Z}{(u)}), \; d_W(\mathbb{Q}, \hat{\prob}(u)) \leq \epsilon(\beta), \\
     \end{cases}  \\
     &= \begin{cases}
       \inf\limits_{\mathbb{Q}, \Pi} & \int_{\mathcal{Z}(u)} H(u;{\rho}) \mathbb{Q}(d \rho), \\
       \st &   \int_{\mathcal{Z}(u) \times \mathcal{Z}(u)} \Norm{\rho -\rho^{\prime}} \Pi(d \rho, d \rho^{\prime}) \leq \epsilon(\beta), \\
       & \Pi \textrm{ is a distribution of } \rho \textrm{ and } \rho^{\prime} \\
       & \quad \quad \textrm{ with marginals } \mathbb{Q} \textrm{ and
       } \hat{\prob}(u) \in
       \subscr{\mathcal{M}}{lt}(\mathcal{Z}{(u)}),
     \end{cases}  \\
     &= \begin{cases}
       \inf\limits_{\mathbb{Q}^{(l)}, \; l \in \mathcal{L}} & \frac{1}{N} \sum\limits_{l\in\mathcal{L}}  \int_{\mathcal{Z}(u)} H(u;{\rho}) \mathbb{Q}^{(l)}(d \rho), \\
       \st & \frac{1}{N} \sum\limits_{l\in\mathcal{L}} \int_{\mathcal{Z}(u)} \Norm{\rho -\rho^{(l)}} \mathbb{Q}^{(l)}(d \rho) \leq \epsilon(\beta), \\
       & \eqref{eq:Rhohat}, \; \mathbb{Q}^{(l)} \in \subscr{\mathcal{M}}{lt}(\mathcal{Z}{(u)}), \; \forall l \in \mathcal{L},  \\
     \end{cases}  \\
   \end{aligned}
  \end{equation*}
     \begin{equation*} \small
       \begin{aligned}
     &= \begin{cases}
       \inf\limits_{\mathbb{Q}^{(l)}, \; l \in \mathcal{L}} \sup\limits_{\lambda \geq 0} \frac{1}{N} \sum\limits_{l\in\mathcal{L}} \int_{\mathcal{Z}(u)} H(u;{\rho}) \mathbb{Q}^{(l)}(d \rho) \\
       \quad \qquad  + \lambda \left( \frac{1}{N} \sum\limits_{l\in\mathcal{L}} \int_{\mathcal{Z}(u)} \Norm{\rho -\rho^{(l)}} \mathbb{Q}^{(l)}(d \rho) -  \epsilon(\beta)   \right), \\
       \st \quad  \eqref{eq:Rhohat}, \; \mathbb{Q}^{(l)} \in \subscr{\mathcal{M}}{lt}(\mathcal{Z}{(u)}), \; \forall l \in \mathcal{L},  \\
\end{cases}
 	\end{aligned}
 \end{equation*}
 where the first equality applies the definition of the expectation
 operation; the second equality uses the definition of Wasserstein
 metric; the third equality exploits the fact that the joint distribution
 $\Pi$ can be characterized by the marginal distribution $\hat{\prob}(u)$ of $\rho^{\prime}$ and the conditional distributions $\mathbb{Q}^{(l)}$ of $\rho$ given $\rho^{\prime}={\rho}^{(l)}$, $l \in \mathcal{L}$, written as
\begin{equation*}
  \Pi:=\frac{1}{N} \sum\limits_{l\in\mathcal{L}} \delta_{\{{\rho}^{(l)}\}} \otimes \mathbb{Q}^{(l)},
\end{equation*}
where admissible sample trajectories $\{{\rho}^{(l)}\}_{l\in\mathcal{L}}$ come
from~\eqref{eq:Rhohat} and each conditional distribution $\mathbb{Q}^{(l)}$ is supported on $\subscr{\mathcal{M}}{lt}(\mathcal{Z}{(u)})$; %
and, on the other hand, the fourth equality
applies the Lagrangian representation of the %
problem.

Then, with an extended version of the strong duality results for the
moment problem~\cite[Lemma~3.4]{AS-MG-MAL:01}, the order of the
$\inf$-$\sup$ operator in the resulting representation of $J(u)$ can
be switched, resulting in the following expression
\begin{equation*} \small
  \begin{aligned}
    J(u) &= \begin{cases} \sup\limits_{\lambda \geq 0}
      \inf\limits_{\mathbb{Q}^{(l)}, \; l \in \mathcal{L}  }
      - \lambda\epsilon(\beta) \\
      \quad \; +\frac{1}{N} \sum\limits_{l\in\mathcal{L}}
      \int_{\mathcal{Z}(u)}
      \left(  \lambda \Norm{{\rho}-{\rho}^{(l)}} +  H(u;{\rho}) \right) \mathbb{Q}^{(l)}(d \rho), \\
      \st \quad \eqref{eq:Rhohat},
      \; \mathbb{Q}^{(l)} \in \subscr{\mathcal{M}}{lt}(\mathcal{Z}{(u)}), \; \forall l \in \mathcal{L},  \\
\end{cases}
 \end{aligned}
\end{equation*}
Move the $\inf$ operator into the sum operator, we have
\begin{equation*} \small
  \begin{aligned}
    J(u) &= \begin{cases} \sup\limits_{\lambda \geq 0}
      - \lambda\epsilon(\beta) \\
      +\frac{1}{N} \sum\limits_{l\in\mathcal{L}} \left\{ \inf\limits_{ \mathbb{Q}^{(l)} \in \subscr{\mathcal{M}}{lt}(\mathcal{Z}{(u)})}
{\mathbb{E}_{\mathbb{Q}^{(l)}} \hspace{-1ex} \left[  \lambda \Norm{{\rho}-{\rho}^{(l)}} +  H(u;{\rho}) \right]}  \right\} \hspace{-0.5ex} , \\
      \st \quad  \eqref{eq:Rhohat}.  \\
\end{cases}
 \end{aligned}
\end{equation*}
For each $l \in \mathcal{L}$, we claim that
\begin{equation*}
  \begin{aligned}
  &  \inf\limits_{\mathbb{Q}^{(l)} \in \subscr{\mathcal{M}}{lt}(\mathcal{Z}{(u)})}
  {\mathbb{E}_{\mathbb{Q}^{(l)}} \left[  \lambda \Norm{{\rho}-{\rho}^{(l)}} +  H(u;{\rho}) \right]} \\
  & \hspace{3.43cm} = \inf\limits_{{\rho} \in \mathcal{Z}(u) }
  \left(\lambda \Norm{{\rho}-{\rho}^{(l)}}+
    {H(u;{\rho})} \right).
  \end{aligned}
\end{equation*}
The above claim can be clarified as the following
\begin{itemize}
  \item[(a)]
  Let $p^{\star}$ denote the value of the second term. Then, for any ${\rho} \in \mathcal{Z}(u)$, we have
  \begin{equation*}
    p^{\star} \leq \lambda \Norm{{\rho}-{\rho}^{(l)}}+
      {H(u;{\rho})},
  \end{equation*}
implying
\begin{equation*}
  p^{\star} \leq {\mathbb{E}_{\mathbb{Q}^{(l)}} \left[  \lambda \Norm{{\rho}-{\rho}^{(l)}} +  H(u;{\rho}) \right]},
\end{equation*}
holds for any probability distribution $\mathbb{Q}^{(l)} \in \mathcal{M}(\mathcal{Z}{(u)})$, so does for $\mathbb{Q}^{(l)} \in \subscr{\mathcal{M}}{lt}(\mathcal{Z}{(u)})$. Therefore, we have
\begin{equation*}
  p^{\star} \leq  \inf\limits_{\mathbb{Q}^{(l)} \in \subscr{\mathcal{M}}{lt}(\mathcal{Z}{(u)})}
  {\mathbb{E}_{\mathbb{Q}^{(l)}} \left[  \lambda \Norm{{\rho}-{\rho}^{(l)}} +  H(u;{\rho}) \right]}.
\end{equation*}

\item[(b)]
Next, we claim that the set $\subscr{\mathcal{M}}{lt}(\mathcal{Z}{(u)})$ contains all the Dirac distributions supported on $\mathcal{Z}(u)$. Then, for any ${\rho}^{\prime} \in \mathcal{Z}(u)$, we have
\begin{equation*}
  \delta_{\{\rho^{\prime} \}} \in \subscr{\mathcal{M}}{lt}(\mathcal{Z}{(u)}),
\end{equation*}
resulting in
\begin{equation*}
  \begin{aligned}
  &  \inf\limits_{\mathbb{Q}^{(l)} \in \subscr{\mathcal{M}}{lt}(\mathcal{Z}{(u)})}
  {\mathbb{E}_{\mathbb{Q}^{(l)}} \left[  \lambda \Norm{{\rho}-{\rho}^{(l)}} +  H(u;{\rho}) \right]} \\
  & \hspace{1.5cm} \leq \lambda \Norm{{\rho}^{\prime}-{\rho}^{(l)}}+
    {H(u;{\rho}^{\prime})}, \; \forall \rho^{\prime} \in  \mathcal{Z}(u),
  \end{aligned}
\end{equation*}
which implies
\begin{equation*}
  \begin{aligned}
  &  \inf\limits_{\mathbb{Q}^{(l)} \in \subscr{\mathcal{M}}{lt}(\mathcal{Z}{(u)})}
  {\mathbb{E}_{\mathbb{Q}^{(l)}} \left[  \lambda \Norm{{\rho}-{\rho}^{(l)}} +  H(u;{\rho}) \right]} \leq p^{\star}.
  \end{aligned}
\end{equation*}
\end{itemize}
Therefore, we equivalently write $J(u)$ as
\begin{equation*} \small
  \begin{aligned}
J(u)=
\sup\limits_{\lambda \geq 0} & - \lambda\epsilon(\beta) + \frac{1}{N}
\sum\limits_{l\in\mathcal{L}} \inf\limits_{{\rho} \in \mathcal{Z}(u) }
\left(\lambda \Norm{{\rho}-{\rho}^{(l)}}+
  {H(u;{\rho})} \right) ,\\
\st &  \;  \eqref{eq:Rhohat}. \\
 \end{aligned}
\end{equation*}
  Finally, we write Problem~\eqref{eq:P2} as follows
{\leqnomode \begin{align} \small
    \label{eq:P2A} \tag{$\mathbf{P2^{\prime}}$}~\sup_{u,\lambda \geq 0} \;& - \lambda\epsilon(\beta) + \frac{1}{N}
    \sum_{l\in\mathcal{L}} \inf_{{\rho} \in \mathcal{Z}(u) }
    \left(\lambda \Norm{{\rho}-{\rho}^{(l)}}+
      {H(u;{\rho})} \right) , \nonumber \\
    \st &  \; \eqref{eq:ue},\; \eqref{eq:Rhohat}. \nonumber
	\end{align}}

\noindent \textbf{Step 2: (Equivalent reformulation of~\eqref{eq:P2A})}
Using the definition of the dual norm and moving its $\sup$ operator
we can write Problem~\eqref{eq:P2A} as
\begin{equation*}
	\begin{aligned}
          \sup_{u,\lambda \geq 0} \;& - \lambda\epsilon(\beta) +
          \frac{1}{N} \sum_{l\in\mathcal{L}}
          \inf_{{\rho} \in \mathcal{Z}(u) } \sup_{\Norm{\nu^{(l)}}_{\star}\leq \lambda}   \\
          & \hspace{2cm} \left({\left\langle \nu^{(l)} ,
                {\rho}-{\rho}^{(l)}\right\rangle}
            +{H(u;{\rho})} \right),\\
          \st &  \; \eqref{eq:ue},\; \eqref{eq:Rhohat}. \\
	\end{aligned}
\end{equation*}
Given $\lambda \geq 0$, the sets $\setdef{\nu^{(l)} \in
  \real^{nT}}{{\Norm{\nu^{(l)}}_{\star}\leq \lambda}}$ are compact for
all $l\in \mathcal{L}$. We then apply the minmax
  theorem between $\inf$ and the second $\sup$
operators. This results in the switch of the operators, and by
combining the two $\sup$ operators we have
\begin{equation*}
	\begin{aligned}
          \sup_{u,\lambda, \nu} & - \lambda\epsilon(\beta) +
          \frac{1}{N} \sum_{l\in\mathcal{L}} \inf_{{\rho} \in
            \mathcal{Z}(u) } \left({\left\langle \nu^{(l)} ,
                {\rho}-{\rho}^{(l)}
              \right\rangle} +{H(u;{\rho})} \right),\\
          \st &  \; \eqref{eq:ue},\; \eqref{eq:Rhohat},\; \lambda \geq 0, \\
          & \Norm{\nu^{(l)}}_{\star}\leq \lambda, \; \forall l \in
          \mathcal{L}.
	\end{aligned}
\end{equation*}
The objective function can be simplified as follows
\begin{equation*}
	\begin{aligned}
          - \lambda\epsilon(\beta) +
          \frac{1}{N}\sum_{l\in\mathcal{L}}{\left\langle -\nu^{(l)} ,
              {\rho}^{(l)}\right\rangle} + \frac{1}{N}
          \sum_{l\in\mathcal{L}} h^{(l)}(u),\\
	\end{aligned}
\end{equation*}
where
\begin{equation*}
	\begin{aligned}
          h^{(l)}(u):=&\inf_{{\rho} \in \mathcal{Z}(u) }
          \left({\left\langle \nu^{(l)} , {\rho}\right\rangle}
            +{H(u;{\rho})} \right), \; \forall l\in \mathcal{L}.
	\end{aligned}
\end{equation*}
For each $l\in \mathcal{L}$, we rewrite $h^{(l)}(u)$ by firstly taking
a minus sign out of the $\inf$ operator, then exploiting the
equivalent representation of $\sup$ operation, and finally using the
definition of conjugate functions. The function $h^{(l)}(u)$ results
in the following form
	\begin{equation*}
	\begin{aligned}
h^{(l)}(u)=&-\sup_{{\rho} \in \mathcal{Z}(u) }
          \left({\left\langle -\nu^{(l)} , {\rho}\right\rangle}
            -{H(u;{\rho})} \right), \\
          =& -\sup_{{\rho} } \left({\left\langle -\nu^{(l)} ,
                {\rho}\right\rangle} -{H(u;{\rho})} -
            \chi_{\mathcal{Z}(u)}({\rho}) \right), \\
          =& - \left[ {H(u;\cdot)} + \chi_{\mathcal{Z}(u)}(\cdot)
          \right]^{\star}(-\nu^{(l)}).
	\end{aligned}
\end{equation*}
Further, we apply the property of the inf-convolution operation and
push the minus sign back into the $\inf$ operator, for each
$h^{(l)}(u)$, $l\in \mathcal{L}$. The representation of $h^{(l)}(u)$
results in the following relation
\begin{equation*}
	\begin{aligned}
          h^{(l)}(u)=& - \inf_{\mu} \left( \left[
              {H(u;\cdot)}\right]^{\star}
            (-\mu^{(l)}-\nu^{(l)}) \right. \\
          & \hspace{2cm} \left. +\left[ \chi_{\mathcal{Z}(u)}(\cdot)
            \right]^{\star}(\mu^{(l)}) \right), \\
          =& \sup_{\mu} \left( -\left[ {H(u;\cdot)}\right]^{\star}
            (-\mu^{(l)}-\nu^{(l)}) \right. \\
          & \hspace{2cm} \left.-\left[
              \chi_{\mathcal{Z}(u)}(\cdot)\right]^{\star}
            (\mu^{(l)}) \right). \\
	\end{aligned}
\end{equation*}
By substituting $-\nu^{(l)}$ by $\nu^{(l)}$, the resulting
optimization problem has the following form
\begin{equation*}
	\begin{aligned}
          \sup_{u,\lambda,\mu, \nu} \; & - \lambda\epsilon(\beta)
          -\frac{1}{N} \sum_{l\in\mathcal{L}} \left( \left[
              {H(u;\cdot)}\right]^{\star}
            (-\mu^{(l)}+\nu^{(l)}) \right. \\
          & \hspace{1cm} \left. +\left[ \chi_{\mathcal{Z}(u)}(\cdot)
            \right]^{\star}(\mu^{(l)}) -{\left\langle \nu^{(l)} ,
                {\rho}^{(l)}\right\rangle} \right), \\
          \st &  \; \eqref{eq:ue},\; \eqref{eq:Rhohat}, \; \lambda \geq 0,  \\
          & \Norm{\nu^{(l)}}_{\star}\leq \lambda, \; \forall
          l \in \mathcal{L}.
	\end{aligned}
\end{equation*}
Given $u$, the strong duality of linear programs is applicable for
the conjugate of the function ${H(u;\cdot)}$ and the support function
$\sigma_{\mathcal{Z}(u)}(\mu^{(l)})$. Using the strong duality and the
definition of the support function, we compute, for each $l\in\mathcal{L}$, the following
\begin{equation*}
  \begin{aligned}
    \left[ {H(u;\cdot)}\right]^{\star}
    &(\nu^{(l)}-\mu^{(l)})
    :=\begin{cases}
      0, & \nu^{(l)}=\mu^{(l)}+ \frac{1}{T}u %
      , %
      \\
      \infty, & {\rm{o.w.}}, \\
          \end{cases}
	\end{aligned}
      \end{equation*}
      and
      \begin{equation*}
        \begin{aligned}
          &\left[ \chi_{\mathcal{Z}(u)}(\cdot)\right]^{\star}(\mu^{(l)})=
          \sigma_{\mathcal{Z}(u)}(\mu^{(l)}) \\
          &= \begin{cases}
            \sup\limits_{\xi^{(l)}} & {\left\langle \mu^{(l)} ,
                \xi^{(l)} \right\rangle}, \\
            \st &\;  0 \leq \xi^{(l)}_e(t) \leq \rho^{c}_{e}(\reviseone{u_e(t)}),~\forall e\in\edges, \; \forall t\in\mathcal{T}, \\
          \end{cases}  \\
          &= \begin{cases}
            \inf\limits_{\eta^{(l)}} &
            \sum\limits_{e\in\edges,t\in\mathcal{T}}{
              \overline{f}_e\overline{\rho}_e\eta_e^{(l)}(t)  }, \\
            \st &  \reviseone{ [\overline{f} \otimes
            \vectorones{T} +
            (\overline{\rho}-\rho^{c}(\overline{u})) \otimes
            \vectorones{T} \circ u  ] \circ \eta^{(l)} } \\
            & \hspace{5cm} -\mu^{(l)} \geq \vectorzeros{nT},
            \\
            &  \eta^{(l)} \geq \vectorzeros{nT},
          \end{cases}
	\end{aligned}
\end{equation*}
We substitute these parts for that in the objective function, take the above $\inf$ operator out of a minus sign, and obtain~\eqref{eq:P3}.

Given that all the reformulations in this step hold with equalities,
we therefore claim that the above problem is equivalent
to~\eqref{eq:P2}. Finally, we claim that the $\sup$ operation is
indeed achievable, \reviseone{because} 1) \reviseone{each component
  of} the variable $u$ is in a finite set $\Gamma$ and 2), for any $u$
that is feasible to the above problem, the above problem with that
fixed $u$ satisfies the Slater's condition, which implies that the
above problem is achievable. We therefore claim~\eqref{eq:P2} is
equivalent to~\eqref{eq:P3}.

\noindent \textbf{Step 3: (Performance guarantees of~\eqref{eq:P3})}
Given any feasible point $(u,{\rho},\lambda,\mu,\nu,\eta)$
of~\eqref{eq:P3}, we denote its objective value by $\hat{J}(u)$.  The
value $\hat{J}(u)$ is a lower bound of~\eqref{eq:P3} and therefore a
lower bound for~\eqref{eq:P2}, i.e., $\hat{J}(u) \leq {J}(u)$. Thus
$\hat{J}(u)$ is an estimate of the certificate for the performance
guarantee~\eqref{eq:perfnew}. Therefore, $(u,\;\hat{J}(u) )$ is a
data-driven solution and certificate pair for~\eqref{eq:P1}. \qed
\end{proof}

\begin{proof}{Lemma~\ref{lemma:linearization}}
  The proof follows by the application of the following proposition on
  each bi-linear term in~\eqref{eq:bilinearNew}:
  \begin{proposition}[Equivalent reformulation of bi-linear terms
    {\cite[Section 2]{FG:75}}]\label{glov} Let
    $\mathcal{Y}\subset\mathbb{R}$ be a compact set. Given a binary
    variable $x$ and a linear function $g(y)$ in a continuous variable
    $y\in \mathcal{Y}$, $z$ equals the quadratic function $xg(y)$ if and
    only if
  \begin{align}
    &\underline{g}x \le z \le \overline{g}x,\nonumber\\
    &g(y)-\overline{g}\cdot(1-x) \le z\le
    g(y)-\underline{g}\cdot(1-x),\nonumber
  \end{align}
  where $\underline{g}=\min_{y\in \mathcal{Y}}\{g(y)\}$ and
  $\overline{g}=\max_{y\in \mathcal{Y}}\{g(y)\}$.  \hfill $\square$
  \end{proposition}
\begin{remark}[Regularization technique]{\rm
  In a later program, we add the following extra constraints to speed up
  the internal computation of solvers
\begin{equation*}
\begin{array}{l}
  \sum_{i\in \mathcal{O}}z^{(l)}_{e,i}(t)=\eta^{(l)}_e(t), \; \\
  \hspace{2cm}\;
  \forall e\in\edges, \; t\in\mathcal{T}\setminus\{0\},\;
  l\in\mathcal{L}, \\
  \sum_{i\in \mathcal{O}}y^{(l)}_{e,i}(t)={\rho}^{(l)}_e(t), \\
  \hspace{2cm}\; \forall e\in\edges, \;
  t\in\mathcal{T}\setminus\{0\},\; l\in\mathcal{L}. \\
\end{array}
\end{equation*}
These are adapted from the binary representation~\eqref{eq:ueNew} and
the definition of $y^{(l)}_{e,i}(t)$ and $z^{(l)}_{e,i}(t)$.}
\qed \end{remark}
 \end{proof}

\begin{proof}{Proposition~\ref{prop:cone}}
  Knowing that ${\rho}^{\star} \geq \vectorzeros{nTN}$ by
  constraints~\eqref{eq:y}, we only need to show $\nu^{\star} \geq
  \vectorzeros{nTN}$. To prove this, let us assume there exists an
  optimizer $\sol^{\star}$ such that, for at least one $\varepsilon \in
  \edges$, $\tau \in \mathcal{T}$ and $\ell \in\mathcal{L}$, it holds
  $\nu^{(\ell),\star}_{\varepsilon}(\tau) <0$. Then, using
  constraint~\eqref{eq:dual2}, we claim
  $\mu^{(\ell),\star}_{\varepsilon}(\tau) <0$. Next, we show the
  contradiction to an optimizer by constructing a feasible solution
  that gives us higher objective value than that resulted from
  $\sol^{\star}$. To achieve this, we perturb variables
  $\lambda^{\star}$, $\mu^{(\ell),\star}_{\varepsilon}(\tau)$ and
  $\nu^{(\ell),\star}_{\varepsilon}(\tau)$, and leave other components
  the same as that in $\sol^{\star}$. With such perturbation, only
  constraints~\eqref{eq:dual1},~\eqref{eq:dual2}, and~\eqref{eq:dual3}
  are varied.

  Let $\sol:=(x^{\star},y^{\star},z^{\star}, {\rho}^{\star},
  \hat{\lambda},\hat{\mu},\hat{\nu},\eta^{\star})$ denote the feasible
  solution we are to construct. We denote by
  $\hbar^{\star}:=\sum_{i\in \mathcal{O}} \gamma^{(i)}
  (\overline{\rho}-\rho^{c}(\overline{u}))\otimes \vectorones{T} \circ
  z^{(l),\star}_i +\overline{f}\otimes \vectorones{T} \circ
  \eta^{(l),\star}$ the unperturbed part in
  constraint~\eqref{eq:dual1} and construct $\hat{\mu}$ as follows
\begin{equation*}
\small
 \hat{\mu}^{(l)}_e(t)=
	\begin{cases}
          \mu^{(l),\star}_e(t), & \; {\rm{if}} \; (e, t, l) \ne (\varepsilon,\tau, \ell), \\
          \min \{ \hbar^{(\ell),\star}_{\varepsilon}(\tau), \; -
          \mu^{(\ell),\star}_{\varepsilon}(\tau) \} , & \; {\rm{o.w}.}
	\end{cases}
\end{equation*}
The above construction ensures the feasibility of
constraints~\eqref{eq:dual1} and furthermore, because
$\hbar^{(\ell),\star}_{\varepsilon}(\tau) \geq 0$ and
$\mu^{(\ell),\star}_{\varepsilon}(\tau) <0$, we claim
$\hat{\mu}^{(\ell)}_{\varepsilon}(\tau) \geq 0$.
 Then let us denote by $g^{\star}:=\frac{1}{T}\sum_{i\in
  \mathcal{O}} \gamma^{(i)}x_{i}^{\star} %
  $ the
unperturbed part of constraints~\eqref{eq:dual2} and construct
variable $\hat{\nu}$ as follows
\begin{equation*}
\small
 \hat{\nu}^{(l)}_e(t)=
	\begin{cases}
          \nu^{(l),\star}_e(t), & \; {\rm{if}} \; (e, t, l) \ne (\varepsilon,\tau, \ell), \\
          \hat{\mu}^{(\ell)}_{\varepsilon}(\tau) + g^{\star}_{\varepsilon}(\tau), & \; {\rm{o.w}.}
	\end{cases}
\end{equation*}
Again,
we have $\hat{\nu}^{(\ell)}_{\varepsilon}(\tau) \geq 0$. Then by
letting $\hat{\lambda}:= \max\{ \lambda^{\star}, \;
\hat{\nu}^{(\ell)}_{\varepsilon}(\tau) \}$,
constraints~\eqref{eq:dual3} are satisfied. In this way, a feasible
solution $\sol$ is constructed.

Next, we evaluate the difference of the objective values
of~\eqref{eq:P4} resulting from $\sol$ and $\sol^{\star}$ in the
following
\begin{equation*}
\begin{aligned}
& \textrm{objective}(\sol)- \textrm{objective}(\sol^{\star})  \\
& \quad = \left(- \hat{\lambda} +\lambda^{\star}\right) \epsilon(\beta)  + \left(\hat{\nu}^{(\ell)}_{\varepsilon}(\tau) - {\nu}^{(\ell),\star}_{\varepsilon}(\tau) \right) {\rho}^{(\ell),\star}_{\varepsilon}(\tau), \\
& \quad \geq \left( -\hat{\lambda} +\lambda^{\star} +\hat{\nu}^{(\ell)}_{\varepsilon}(\tau) - {\nu}^{(\ell),\star}_{\varepsilon}(\tau) \right) \epsilon(\beta), \\
& \quad = \left( \min\{ -\lambda^{\star}, \; -\hat{\nu}^{(\ell)}_{\varepsilon}(\tau) \}+\lambda^{\star} +\hat{\nu}^{(\ell)}_{\varepsilon}(\tau) \right) \epsilon(\beta) \\
&\qquad\qquad\qquad\qquad \qquad\qquad\qquad\qquad- {\nu}^{(\ell),\star}_{\varepsilon}(\tau)\epsilon(\beta), \\
& \quad >0,
\end{aligned}
\end{equation*}
where the first equality cancels out unperturbed terms; the second
inequality applies Assumption~\ref{assump:3} and the fact that
$\hat{\nu}^{(\ell)}_{\varepsilon}(\tau) \geq 0$ and
${\nu}^{(\ell),\star}_{\varepsilon}(\tau) <0$; the third equality
applies construction of $\hat{\lambda}$; and the last inequality is
achieved by summing the nonnegative first term and the strict positive
second term.

By the above computation, we constructed a feasible solution $\sol$
with a higher objective value than that of $\sol^{\star}$, contradicting
the assumption that $\sol^{\star}$ is an optimizer.
\qed \end{proof}

\begin{proof}{Lemma~\ref{lemma:coneP4A}}
  Let us denote by (P4$''$) the Problem~\eqref{eq:P4} with an extra
  set of constraints $\nu \geq \vectorzeros{nTN}$. We prove the lemma
  in two steps.

\noindent \textbf{Step 1: (Equivalence of optimizers sets)} First, we use Proposition~\ref{prop:cone} to claim
  that the set of optimizers of~\eqref{eq:P4} is the same as that of
  (P4$''$). Second, we claim that for any optimizer of~\eqref{eq:P4A},
  all the constraints in~\eqref{eq:cone1} are active. This means that
  the set of optimizers of (P4$''$) are the same as the projection of
  that of~\eqref{eq:P4A}. Therefore, the optimizers set of~\eqref{eq:P4} and~\eqref{eq:P4A} are equivalent.

\noindent \textbf{Step 2: (Performance guarantees)} First, any feasible solution of~\eqref{eq:P4A} correspond to a feasible solution of~\eqref{eq:P4}. This holds because any feasible solution of~\eqref{eq:P4A} satisfies all the constraints of~\eqref{eq:P4}. Next, the objective value of~\eqref{eq:P4A} gives a lower bounds of that of~\eqref{eq:P4}. This can be verified using constraints~\eqref{eq:cone1}. Finally, the performance guarantees~\eqref{eq:perfnew} of feasible solution of~\eqref{eq:P4A} can be derived from that of~\eqref{eq:P4} as in Remark~\ref{remark:P4}.
\qed \end{proof}

\begin{proof}{Lemma~\ref{lemma:vartheta}}
  We construct $\overline{\vartheta}$ by showing boundedness of $\nu
  \circ {\rho}$. It's known for each $e\in\edges$, $t\in\mathcal{T}$
  and $l\in\mathcal{L}$, the density $\rho^{(l)}_e(t)$ is nonnegative
  and bounded above by $\max_{e\in\edges}\{\overline{\rho}_e \}$. Then
  we only need to find the upper bound of $\nu^{(l)}_e(t)$. By
  Assumption~\ref{assump:2}, the variable $\eta^{(l)}_e(t)$ is
  bounded. Then computations via constraints~\eqref{eq:dual1}
  and~\eqref{eq:dual2} result in upper bound of $\nu^{(l)}_e(t)$ as
  the following
\begin{equation*}
\begin{aligned}
\nu^{(l)}_e(t) \leq & \max_{e\in\edges, \reviseone{u_e(t)} \in \Gamma} \left\{ {\left(\overline{f}_e +\reviseone{u_e(t)}(\overline{\rho}_e- \rho^{c}_e(\overline{u}_e))\right) \overline{\eta} + \frac{1}{T}\reviseone{u_e(t)}}\right\} ,\\
=& \max_{e\in\edges} \left\{ { \overline{u}_e\overline{\rho}_e \overline{\eta} + \frac{1}{T}\overline{u}_e}\right\}. \\
\end{aligned}
\end{equation*}
By letting $\overline{\vartheta}=\sqrt{\max_{e\in\edges} \left\{ { \overline{u}_e\overline{\rho}^2_e \overline{\eta} + \frac{1}{T}\overline{u}_e}\overline{\rho}_e\right\}}$, we complete the proof.
\qed \end{proof}

\bibliographystyle{IEEEtran}

\bibliography{../../../bib/alias,../../../bib/SMD-add,../../../bib/SM,../../../bib/JC,../../../bib/Main,../../../bib/Main-sonia}

\begin{IEEEbiography}[{\includegraphics[width=1in,height=1.25in,clip,keepaspectratio]{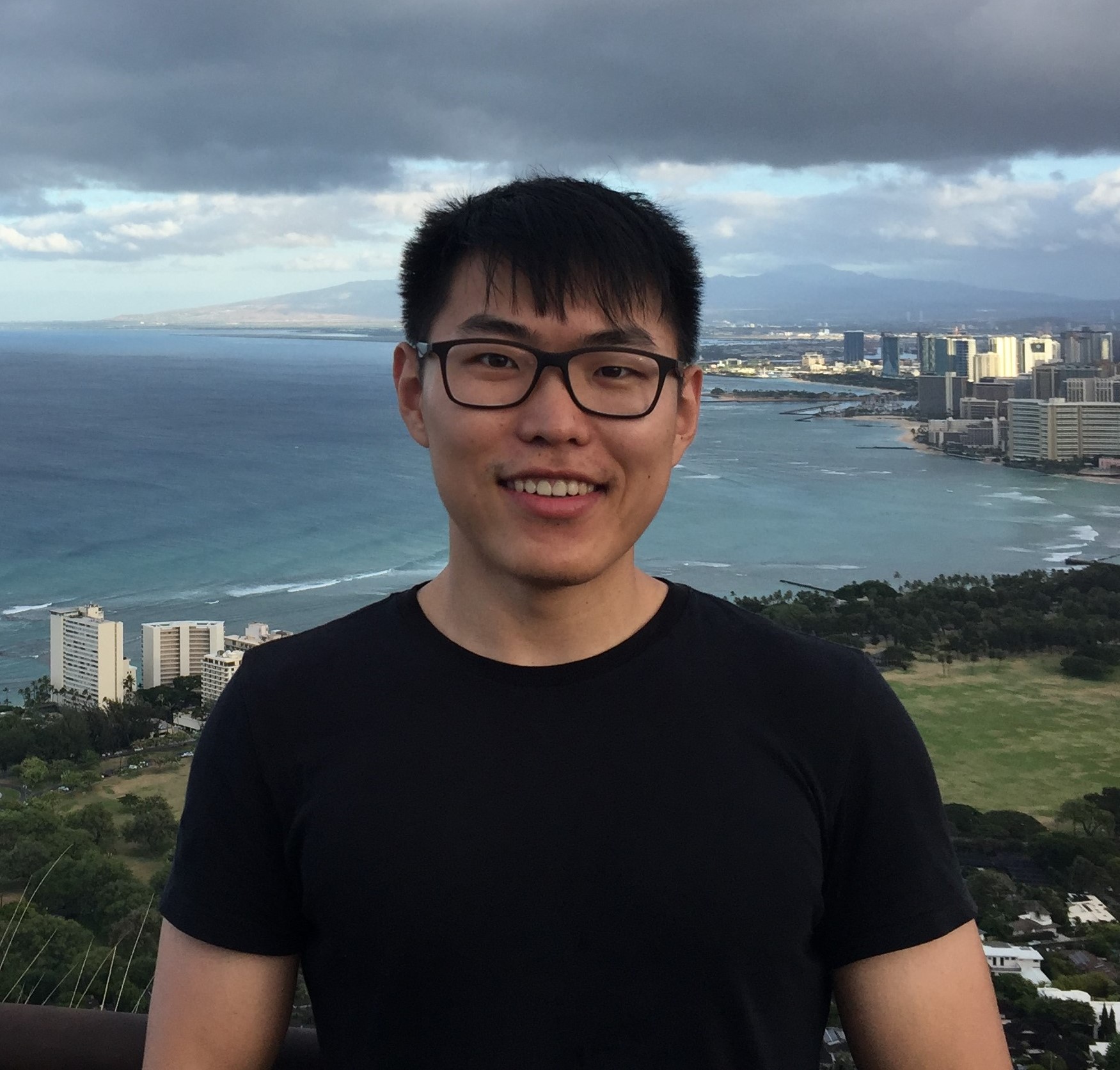}}]{Dan
    Li} received the B.E. degree in automation from the Zhejiang
  University, Hangzhou, China, in 2013, the M.Sc. degree in chemical
  engineering from Queen's University, Kingston, Canada, in 2016. He
  is currently a Ph.D. student at University of California, San Diego,
  CA, USA. His current research interests include data-driven systems
  and optimization, dynamical systems and control, optimization
  algorithms, applied computational methods, and stochastic systems. He
  received Outstanding Student Award from Zhejiang University in 2012,
  Graduate Student Award from Queen's University in 2014, and
  Fellowship Award from University of California, San Diego, in 2016.
\end{IEEEbiography}
\begin{IEEEbiography}[{\includegraphics[width=1in,height=1.25in,clip,keepaspectratio]{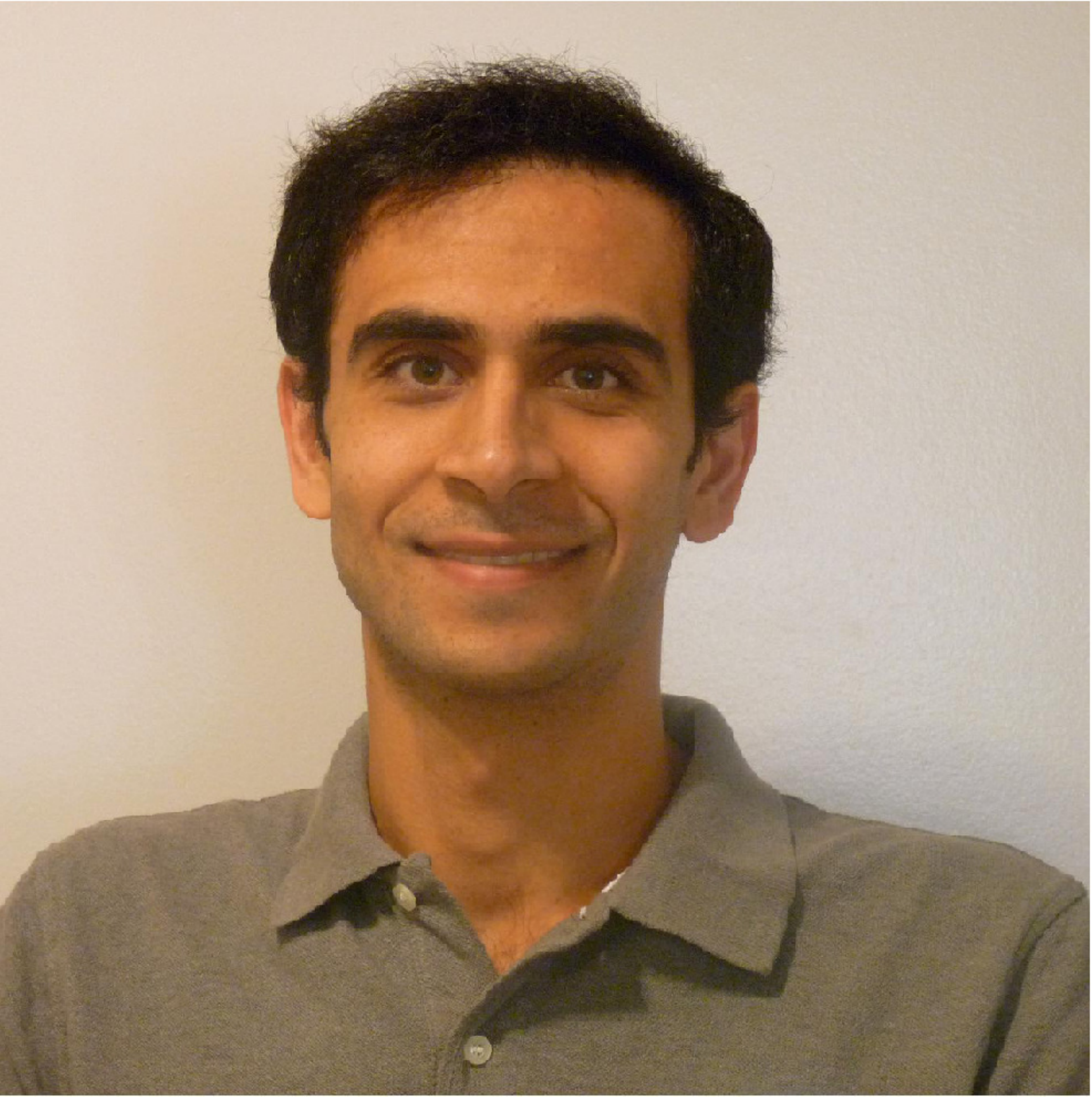}}]{Dariush Fooladivanda}
received the Ph.D. degree from the University of Waterloo, in 2014, and a B.S. degree from the Isfahan University of Technology, all in electrical engineering. He is currently a Postdoctoral Research Associate in the Department of Electrical Engineering and Computer Sciences at the University of California Berkeley. His research interests include theory and applications of control and optimization in large scale dynamical systems.
\end{IEEEbiography}
\begin{IEEEbiography}[{\includegraphics[width=1in,height=1.25in,clip,keepaspectratio]{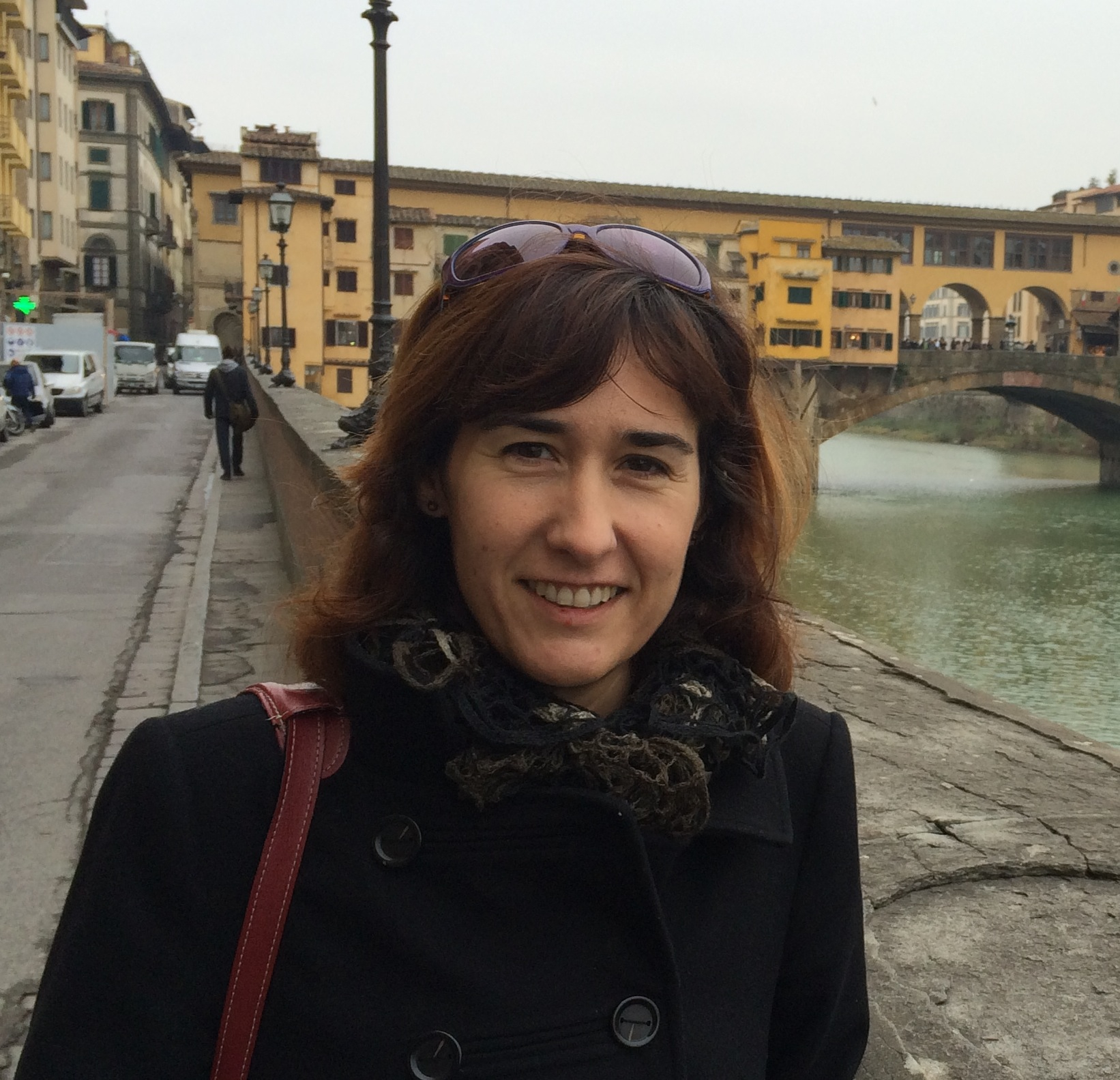}}]{Sonia
   Martínez} is a Professor at the Department of Mechanical and
 Aerospace Engineering at the University of California, San
 Diego. She received her Ph.D. degree in Engineering Mathematics from
 the Universidad Carlos III de Madrid, Spain, in May 2002. Following
 a year as a Visiting Assistant Professor of Applied Mathematics at
 the Technical University of Catalonia, Spain, she obtained a
 Postdoctoral Fulbright Fellowship and held appointments at the
 Coordinated Science Laboratory of the University of Illinois,
 Urbana-Champaign during 2004, and at the Center for Control,
 Dynamical systems and Computation (CCDC) of the University of
 California, Santa Barbara during 2005.

 Her research interests include networked control systems,
 multi-agent systems, and nonlinear control theory with applications
 to robotics and cyber-physical systems. For her work on the control
 of underactuated mechanical systems she received the Best Student
 Paper award at the 2002 IEEE Conference on Decision and Control. She
 co-authored with Jorge Cortés and Francesco Bullo "Motion
 coordination with Distributed Information" for which they received
 the 2008 Control Systems Magazine Outstanding Paper Award. She is a
 Senior Editor of the IEEE Transactions on Control of Networked
 Systems and an IEEE Fellow.
\end{IEEEbiography}
\end{document}